\documentclass[11pt,a4paper]{article}
\usepackage[margin=27mm]{geometry}
\usepackage{amsmath,amssymb,amsthm,mathtools,bm,booktabs,longtable,array}
\usepackage[T1]{fontenc}
\usepackage{lmodern,microtype}
\usepackage[colorlinks=true,linkcolor=blue,citecolor=blue,urlcolor=blue]{hyperref}
\usepackage{enumitem,listings}
\allowdisplaybreaks[2]
\numberwithin{equation}{section}
\newtheorem{theorem}{Theorem}[section]
\newtheorem{proposition}[theorem]{Proposition}
\newtheorem{lemma}[theorem]{Lemma}
\newtheorem{corollary}[theorem]{Corollary}
\theoremstyle{definition}\newtheorem{definition}[theorem]{Definition}

\theoremstyle{remark}\newtheorem{remark}[theorem]{Remark}
\newcommand{\Wres}{\operatorname{Wres}}
\newcommand{\Res}{\operatorname{Res}}
\newcommand{\tr}{\operatorname{tr}}
\newcommand{\Ric}{\operatorname{Ric}}
\newcommand{\Id}{\operatorname{Id}}
\newcommand{\hc}{\widehat c}
\newcommand{\dd}{\mathrm d}
\newcommand{\vol}{\mathrm{vol}}

\newcommand{\A}{\mathcal A}
\newcommand{\HH}{\mathcal H}
\newcommand{\RR}{\mathcal R}
\newcommand{\fres}{\mathfrak f}

\title{Riemannian Curvature Functionals for de Rham Spectral Triples}
\author{Yong Wang\\[4pt]
\small School of Mathematics and Statistics,\\
\small Northeast Normal University,\\
\small Changchun Jilin, 130024, China\\[3pt]
\small E-mail: \href{mailto:wangy581@nenu.edu.cn}{\texttt{wangy581@nenu.edu.cn}}}
\date{}
\begin{document}
\maketitle
\begin{abstract}
We compute mixed Clifford residue functionals for the Hodge--de Rham operator in even dimensions. A symbol-composition calculation recovers the Riemann and Einstein tensors and yields two ordered curvature functionals with the same classical limit. For specified truncation products on warped collars, we evaluate all five boundary cases and prove that the curvature-recovering antisymmetrizations have no boundary correction. We extend the constructions to spectral data equipped with two one-form actions. Explicit noncommutative torus and finite-direction examples distinguish the functionals and exhibit failures of module balancing and Riemann symmetries. A conformal example shows that a naive exterior lift need not recover classical curvature. The additional structures and limitations of each extension are stated explicitly.
\end{abstract}
\noindent\textbf{Keywords.} Hodge--de Rham operator; noncommutative residue; Riemann curvature; boundary symbols; noncommutative torus; spectral tensoriality.\par
\noindent\textbf{2020 Mathematics Subject Classification.} 58J42, 58B34, 53B20, 46L87.

\clearpage
\setcounter{tocdepth}{2}
\begingroup
\small
\tableofcontents
\endgroup
\clearpage

\section{Introduction and scope}
\subsection{From scalar spectral invariants to the full curvature tensor}
The noncommutative residue turns the homogeneous part of a pseudodifferential symbol into a local geometric density. Its importance lies not only in its trace property, but also in the possibility of recovering geometric information without introducing that information explicitly into an operator insertion. The work of Wodzicki and Guillemin established the residue framework \cite{Wodzicki,Guillemin}; the gravitational interpretation of inverse Dirac powers was developed in the Kastler--Kalau--Walze theorem \cite{Kastler,KalauWalze}. In that theorem the recovered invariant is scalar curvature. A natural next question is whether more selective insertions recover tensors before their indices are contracted.

The spectral metric and Einstein functionals of D\k{a}browski, Sitarz and Zalecki \cite{DSZmetric} provide a systematic realization of this idea. Their spectral torsion functional \cite{DSZtorsion} addresses a different geometric order. The spectral forms of Wang, Wang and Liu \cite{WWL} further motivate multilinear expressions built from the Hodge--de Rham operator. These constructions suggest that the choice of representation and the ordering of inserted one-forms are essential parts of a spectral tensor formula.

On the full exterior bundle there are two natural actions, $c(X)=\varepsilon(X^\flat)-\iota_X$ and $\hc(X)=\varepsilon(X^\flat)+\iota_X$. They have opposite Clifford signs and anticommute with each other. Their coexistence makes the de Rham representation particularly suitable for separating curvature components. Indeed, the Weitzenbock endomorphism of $D=\dd+\delta$ contains the mixed term $-R_{ijkl}c_ic_j\hc_k\hc_l/8$. The two unhatted and two hatted slots can test the two skew index pairs of $R$. Replacing this endomorphism by the spinorial scalar potential would remove precisely the term that is needed for full curvature recovery.

We study four ordered residues $F_\ell=\Wres(A_\ell D B_\ell D^{1-n})$ on manifolds of even dimension $n=2m$. Their interior densities contain scalar, Ricci, Einstein and Riemann contributions. In particular, the first-pair antisymmetrizations satisfy
\[
 I_2(X,Y;Z,W)-I_2(Y,X;Z,W)=(m-1)R(X,Y,Z,W),
\]
and the analogous difference for $I_3$ is its negative. Thus a single antisymmetrization suffices for recovery when $n>2$. Inserting arbitrary scalar weights on the left determines the local tensor, rather than only its unweighted integral. The exceptional dimension two is retained explicitly: these particular residues vanish there, and dividing by $m-1$ would be meaningless.

\subsection{Why a boundary calculation requires a separate theorem}
On a manifold with boundary, an interior pseudodifferential expression does not uniquely determine a boundary functional. A product of truncated factors acquires a singular Green contribution. The residue theorem of Fedosov, Golse, Leichtnam and Schrohe \cite{FGLS}, together with the boundary analysis discussed by Schrohe \cite{Schrohe}, provides the appropriate trace. Yong Wang's gravitational calculations \cite{WangBoundary,WangGravity} show how the homogeneous-symbol formula can be evaluated by normal contour integration and tangential spherical moments.

The six-dimensional spectral Einstein calculations of Jian Wang and Yong Wang \cite{WWJMP} are especially relevant to the presentation adopted here. They distinguish different truncation factorizations and organize the boundary calculation into five cases before rewriting the answer in geometric quantities. General-dimensional methods also appear in \cite{WW6,WangGeneral,WuWangGeneral}. These papers concern different operator insertions. Their coefficients cannot simply be transferred to the mixed exterior representation used here.

We fix the factorization $\pi^+(A_\ell D B_\ell D^{-2})\circ\pi^+(D^{3-n})$ and a warped collar $g=h(x_n)^{-1}g_{\partial M}+\dd x_n^2$. For this choice there is exactly one loss of homogeneous degree in the boundary order equation, hence five cases. We retain all first derivatives of the inserted fields. The final answer is expressed through inner products, normal components, tangential divergences, covariant derivatives and the mean curvature $K_{\partial}=-(n-1)h'(0)/2$ in our inward-normal convention. All dimension-dependent coefficients are explicit rational functions times a central binomial coefficient. This formulation avoids treating a finite symbol recipe as the final geometric answer.

An additional feature emerges after the five cases are summed. The boundary densities of the second and third insertions are symmetric in $X,Y$, even though individual cases of the third insertion need not be. Their first-pair antisymmetrizations therefore have no boundary correction for this factorization. This yields a boundary version of the curvature-recovery corollary. It is a property of the specified extension, not a claim that all factorizations of a curvature residue share the same boundary behavior.

The inserted-Clifford calculations of Wu and Wang \cite{WuWangInserted} offer a complementary reference: they show explicitly how products of Clifford actions produce metric contractions and how normal derivatives of those contractions enter boundary residues. Here the mixed hatted action introduces additional contractions and requires the full exterior trace. Throughout the paper this is the ordinary trace on all degrees, of rank $2^n$, rather than a spinor trace or a supertrace.

\subsection{Noncommutative extensions and their limitations}
Classical recovery and noncommutative tensoriality are separate mathematical questions. A residue expression can be meaningful on spectral data without descending to a balanced multilinear map on a bimodule. We introduce differential de Rham data with a nilpotent Hilbert-complex differential and distinguish them from weaker data carrying two one-form actions. Under the stated regularity and residue assumptions, the mixed commutator insertion defines a curvature functional. Module balancing, pair interchange and the first Bianchi identity are additional properties to be checked.

On the smooth noncommutative torus, the flat exterior operator gives zero curvature functional for arbitrary algebra-valued inputs. For zero-order perturbations we compute the relevant ordered potential and obtain explicit failures of balancing and of Riemann symmetries. These examples identify limitations of the proposed functional; they do not establish the nonexistence of other noncommutative curvature constructions. A conformal exterior amplification supplies a further warning: vanishing numerical-frame components can reflect an inadequate second action rather than a flat metric. The gauge and two-sheet examples are treated only after their additional exterior data have been specified.

\subsection{Conventions and organization}
Our residue includes $(2\pi)^{-n}$. To compare a manifold formula with an unnormalized boundary residue convention, multiply the entire formula, including its boundary term, by $(2\pi)^n$. Negative powers are parametrices modulo smoothing operators. No APS, absolute or relative boundary realization is being imposed: the boundary object belongs to the transmission calculus.

Sections~\ref{sec:conventions}--\ref{sec:four} establish the exterior trace identities and the four local residues. Section~\ref{sec:curvature} explains weighted curvature recovery. Section~\ref{sec:boundary} derives all five boundary cases, gives their geometric sums in every even dimension $n\ge4$, and proves the boundary recovery corollary. Sections~\ref{sec:data}--\ref{sec:examples} develop the proposed spectral data, the torus calculations and other examples. The scope of the boundary theorem is the stated warped collar, rather than an arbitrary second fundamental form. No uniqueness or reconstruction theorem for general noncommutative de Rham data is asserted.

\section{Exterior calculus and normalizations}\label{sec:conventions}
Let $(M,g)$ be an oriented closed Riemannian manifold, $n=2m$, and let
\begin{equation}
 \mathcal E=\Lambda^*T^*M\otimes\mathbb C,\qquad r=2^n,\qquad
 D=\dd+\delta=\sum_i c_i\nabla_i.
\end{equation}
A spin structure is not required. For a real vector field $X$, put
\begin{equation}
 c(X)=\varepsilon(X^\flat)-\iota_X,\qquad
 \hc(X)=\varepsilon(X^\flat)+\iota_X.
\end{equation}
Thus $c(X)^*=-c(X)$, $\hc(X)^*=\hc(X)$, and
\begin{equation}\label{eq:CAR}
 \{c(X),c(Y)\}=-2g(X,Y),\quad
 \{\hc(X),\hc(Y)\}=2g(X,Y),\quad
 \{c(X),\hc(Y)\}=0.
\end{equation}
We use
\begin{equation}
 R(X,Y)=[\nabla_X,\nabla_Y]-\nabla_{[X,Y]},\quad
 R_{ijkl}=g(R(e_i,e_j)e_k,e_l),\quad
 \Ric_{jk}=\sum_iR_{ijki}.
\end{equation}
Consequently a round sphere of curvature $K$ has
$R_{ijkl}=K(\delta_{il}\delta_{jk}-\delta_{ik}\delta_{jl})$ and positive scalar curvature $s=n(n-1)K$.

The left-symbol convention is
\begin{equation}\label{eq:star}
 p\#q\sim\sum_\alpha\frac1{\alpha!}
 \partial_\xi^\alpha p\,D_x^\alpha q,
 \qquad D_{x_j}=\frac1i\partial_{x_j}.
\end{equation}
The residue is
\begin{equation}\label{eq:Wres}
 \Wres(P)=(2\pi)^{-n}\int_M\int_{|\xi|=1}
 \tr_{\mathcal E}\sigma_{-n}(P)(x,\xi)\,\dd S(\xi)\,\dd x.
\end{equation}
Let $\nu_{n-1}=\vol(S^{n-1})$. We will repeatedly use
\begin{equation}\label{eq:Cm}
 C_m=\frac{r\nu_{n-1}}{(2\pi)^n}
 =\frac{2^{2m+1}}{(4\pi)^m\Gamma(m)}
 =\frac2{\pi^m\Gamma(m)}.
\end{equation}
Negative powers mean classical parametrices, or generalized inverse powers with kernel projection removed. The difference is smoothing and has no residue.

\subsection{Curvature of the exterior connection}
The induced connection has curvature
\begin{equation}\label{eq:Omega}
 \Omega_{ij}=\frac14R_{ijkl}(c_kc_l-\hc_k\hc_l).
\end{equation}
Indeed both sides have commutator $c(R(e_i,e_j)V)$ with $c(V)$ and commutator $\hc(R(e_i,e_j)V)$ with $\hc(V)$; both annihilate constants as the induced action on exterior forms. Squaring $D$ gives
\begin{align}
 L:=D^2&=\nabla^*\nabla+E,\label{eq:L}\\
 E&=\frac12c_ic_j\Omega_{ij}
 =\frac s4-\frac18R_{ijkl}c_ic_j\hc_k\hc_l.\label{eq:E}
\end{align}
For the last equality the first Bianchi identity eliminates the degree-four part of $R_{ijkl}c_ic_jc_kc_l$ and its scalar part is $2s$. Formula~\eqref{eq:E} is the source of the mixed Riemann term. Replacing $E$ by $s/4$ would be the spinor calculation and would give a different answer.

\paragraph{Verification on exterior generators.}
Let $\varepsilon_k$ and $\iota_k$ be wedge multiplication and contraction.
Their relations give
\[
 \tfrac14R_{ijkl}(c_kc_l-\hc_k\hc_l)
 =-R_{ijkl}\varepsilon_k\iota_l.
\]
On a one-form $e^a$ this is $-R_{ijka}e^k=R_{ijak}e^k$, which is
exactly the induced curvature action under the metric identification.
Both sides extend as degree-zero derivations of the exterior algebra,
so agreement on scalars and one-forms proves~\eqref{eq:Omega} on every
degree. At a geodesic center, splitting $c_ic_j\nabla_i\nabla_j$ into
its diagonal and skew parts gives $\nabla^*\nabla+c_ic_j\Omega_{ij}/2$.
In the pure $c$ contraction the four-generator part vanishes by Bianchi,
the two-generator part is a symmetric Ricci tensor contracted with a
skew Clifford bivector, and the scalar part is
\[
 R_{ijkl}(\delta_{ij}\delta_{kl}-\delta_{ik}\delta_{jl}
                       +\delta_{il}\delta_{jk})=2s.
\]
This proves~\eqref{eq:E} with the stated curvature sign.

\subsection{Trace identities}
Write $\tau_\Lambda=r^{-1}\tr_{\mathcal E}$. The basic contractions are
\begin{align}
 \tau_\Lambda(c(X)c(Y))&=-g(X,Y),&
 \tau_\Lambda(\hc(X)\hc(Y))&=g(X,Y),\\
 \tau_\Lambda(c(X)c(Y)c(Z)c(W))&=K(X,Y,Z,W),\label{eq:fourtrace}\\
 \tau_\Lambda(\hc(X)\hc(Y)\hc(Z)\hc(W))&=K(X,Y,Z,W),\\
 \tau_\Lambda(\hc(X)\hc(Y)c(Z)c(W))&=-g(X,Y)g(Z,W),
\end{align}
where
\begin{equation}
 K(X,Y,Z,W)=g(X,Y)g(Z,W)-g(X,Z)g(Y,W)+g(X,W)g(Y,Z).
\end{equation}
These identities follow by moving generators using~\eqref{eq:CAR} and keeping the scalar Clifford component. On the full exterior representation every nonempty ordered monomial in the $2n$ distinct generators has zero trace. This provides an exact finite algorithm for all traces in the paper.

\paragraph{Why scalar Clifford extraction computes the ordinary trace.}
The $2n$ matrices $c_i,\hc_j$ mutually anticommute when their labels
are distinct, and each is invertible. For a nonempty ordered monomial
of even length, conjugation by a generator occurring in it changes its
sign. For odd length, conjugation by a generator not occurring in it
changes its sign; such a generator exists because $2n$ is even.
Trace invariance under conjugation therefore makes every such monomial
traceless. Only the scalar component remains, with trace $r=2^n$.
In particular
\[
 \tau_\Lambda(P(c)Q(\hc))=\tau_\Lambda(P(c))\tau_\Lambda(Q(\hc)).
\]
Moving generators through one another gives the three pairings
\[
 g(X,Y)g(Z,W),\qquad -g(X,Z)g(Y,W),\qquad g(X,W)g(Y,Z)
\]
in the four-generator trace. For matrix-valued coefficients this argument applies to the
Clifford matrices while retaining every coefficient in its original
order. No supertrace is used.

\section{Symbol calculation at curvature order}\label{sec:symbol}
We record enough of the recursion to fix signs and to permit an independent implementation. If $l_2+l_1+l_0$ is the symbol of $L$, the first inverse symbols are
\begin{align}
 b_{-2}&=l_2^{-1},\\
 b_{-3}&=-l_2^{-1}\left(l_1b_{-2}+\partial_{\xi_j}l_2D_{x_j}b_{-2}\right),\\
 b_{-4}&=-l_2^{-1}\bigg(l_0b_{-2}+l_1b_{-3}
 +\partial_{\xi_j}l_2D_{x_j}b_{-3}
 +\partial_{\xi_j}l_1D_{x_j}b_{-2}\notag\\
 &\hspace{35mm}+\frac12\partial_{\xi_j}\partial_{\xi_k}l_2
 D_{x_j}D_{x_k}b_{-2}\bigg).\label{eq:inverse-rec}
\end{align}
For integer $a\geq1$, the three needed symbols of $L^{-a}$ are the terms of degrees $-2a,-2a-1,-2a-2$ in the $a$-fold star product of $b$. Thus all recursions below terminate after two losses of homogeneous degree.

At a normal-coordinate center with a radial parallel bundle frame,
\begin{equation}
 g_{ij}=\delta_{ij},\quad\partial_kg_{ij}=0,\quad
 \omega_i=0,\quad\partial_i\omega_j=\tfrac12\Omega_{ij}.
\end{equation}
The metric and density jets, with our curvature ordering, are
\begin{align}
 g_{ij}(x)&=\delta_{ij}+\tfrac13R_{ikjl}x^kx^l+O(|x|^3),\\
 g^{ij}(x)&=\delta_{ij}-\tfrac13R_{ikjl}x^kx^l+O(|x|^3),\\
 \sqrt{\det g}(x)&=1-\tfrac16\Ric_{kl}x^kx^l+O(|x|^3).
\end{align}
The angular rules required for the recursion are
\begin{align}
 \int_{S^{n-1}}\xi_i\xi_j\,\dd S&=\frac{\nu_{n-1}}n\delta_{ij},\\
 \int_{S^{n-1}}\xi_i\xi_j\xi_k\xi_l\,\dd S
 &=\frac{\nu_{n-1}}{n(n+2)}
 (\delta_{ij}\delta_{kl}+\delta_{ik}\delta_{jl}+\delta_{il}\delta_{jk}).
\end{align}

\begin{lemma}[The two universal contractions]\label{lem:contractions}
For a Laplace-type operator $L=\nabla^*\nabla+E$ and a zero-order endomorphism $T$, the local residue densities satisfy
\begin{align}
 \operatorname{res}_x(TL^{1-m})
 &=\frac{\nu_{n-1}}{(2\pi)^n}(m-1)
 \tr\left[T\left(\frac s6-E\right)\right]\dd V_g,\label{eq:static}\\
 \operatorname{res}_x(T\nabla_i\nabla_jL^{-m})
 &=\frac{\nu_{n-1}}{(2\pi)^n}\tr(T H_{ij})\dd V_g,\label{eq:H}\\
 H_{ij}&=\frac16\Ric_{ij}-\frac{s}{12}\delta_{ij}
 +\frac12\delta_{ij}E+\frac12\Omega_{ij}.
\end{align}
Here $\nabla_i\nabla_j$ in~\eqref{eq:H} is evaluated at the normal-coordinate center. Also
$\operatorname{res}_x(T\nabla_iL^{-m})=0$.
\end{lemma}
\begin{proof}
We use only the left-symbol composition rule~\eqref{eq:star}. All the
following identities are evaluated at a fixed normal-coordinate center
$x_0=0$, in a radial parallel frame of the bundle. Put $q=g^{ab}\xi_a\xi_b$
and $q(0,\xi)=|\xi|^2$. Coordinate indices in this proof are raised with
$\delta$ at the center.

\paragraph{Metric, Christoffel and connection derivatives.}
The curvature conventions are
$\Ric_{ij}=\sum_aR_{aija}$ and $s=\sum_i\Ric_{ii}$.
Differentiating the metric expansion gives
\begin{align*}
 \partial_r\partial_s g_{ab}(0)&=\tfrac13(R_{arbs}+R_{asbr}),\\
 \partial_r\partial_s g^{ab}(0)&=-\tfrac13(R_{arbs}+R_{asbr}),\\
 \partial_r\Gamma^c_{ab}(0)&=-\tfrac13(R_{brac}+R_{arbc}),\\
 \partial_r(g^{ab}\Gamma^c_{ab})(0)&=\tfrac23\Ric_{rc}.
\end{align*}
For example, the third line follows by differentiating
$\Gamma^c_{ab}=g^{cd}(\partial_ag_{bd}+\partial_bg_{ad}-\partial_dg_{ab})/2$
and inserting the first line; the terms involving first derivatives of
$g$ vanish at $0$. The skew symmetries and the first Bianchi identity
reduce the resulting six curvature terms to the displayed expression.
The last line follows by setting $a=b$ and summing, using
$\sum_aR_{arac}=-\Ric_{rc}$. Thus these formulas fix the signs of both
the Christoffel derivative and the Ricci contraction.

Write $\nabla_a=\partial_a+\omega_a$. Radial parallel transport gives
$x^a\omega_a(x)=0$. Its quadratic Taylor coefficient implies
$\partial_r\omega_s+\partial_s\omega_r=0$ at $0$, whereas the curvature
equation gives $\partial_r\omega_s-\partial_s\omega_r=\Omega_{rs}$.
Consequently $\partial_r\omega_s=\Omega_{rs}/2$ and
$\sum_a\partial_a\omega_a=0$.

Introduce the metric polynomial and the first-order coefficient
\[
 Q_{rs}:=\partial_r\partial_s q(0,\xi)
 =-\tfrac13(R_{arbs}+R_{asbr})\xi_a\xi_b,
 \qquad v^b=g^{ac}\Gamma^b_{ac}-2g^{ab}\omega_a.
\]
At the center,
\[
 v^b=0,\qquad V_{rb}:=\partial_rv^b
 =\tfrac23\Ric_{rb}-\Omega_{rb}.
\]
The contractions used below are
\begin{align*}
 \sum_r Q_{rr}&=\tfrac23\Ric_{ab}\xi_a\xi_b,
 &\sum_s\xi_sQ_{rs}&=0,\\
 \frac1{\nu_{n-1}}\int_{S^{n-1}}Q_{ij}\,\dd S
 &=\frac{2}{3n}\Ric_{ij}.
\end{align*}
For the zero contraction, each of the two curvature summands contains
an antisymmetric index pair contracted against a symmetric product of
three covariables. The integral follows by replacing $\xi_a\xi_b$ by
$\delta_{ab}/n$. These are polynomial identities, not assumptions on
the curvature.

\paragraph{The inverse-power symbols quoted from the literature.}
We use the inverse-power symbols of \cite[Proposition~2.2,
equation~(3.4) and Lemma~3.3]{DSZmetric}, translated to our curvature
ordering and the convention $L=\nabla^*\nabla+E$.
We state the required symbols and their jets directly, without
rederiving the inverse-power recursion. For $a\geq1$ put
$S_{a,j}=\sigma_{-2a-j}(L^{-a})$. At $0$ they are
\begin{align*}
 S_{a,0}&=q^{-a},&S_{a,1}&=0,\\
 S_{a,2}&=-aEq^{-a-1}
       +\frac{a(a+1)}3\Ric_{bc}\xi_b\xi_cq^{-a-2},\\
 \partial_rS_{a,0}&=0,&
 \partial_rS_{a,1}&=-iaV_{rb}\xi_bq^{-a-1},\\
 \partial_r\partial_sS_{a,0}&=-aQ_{rs}q^{-a-1}.
\end{align*}
Here $q=|\xi|^2$ in these center values, and $Q,V$ are the geometric
jets defined above. The last two lines are the first and second jets
of the cited normal-coordinate expansions; $\xi_sQ_{rs}=0$ has already
been used. Taking $a=m$ gives the symbols of $L^{-m}$ needed for the
second-order insertion, and taking $a=m-1$ gives the static case.
In particular, $S_{a,1}(0)=0$ does not permit one to discard its
nonzero first derivatives in a symbol composition. The subsequent
composition and angular contractions are given in full below.

\paragraph{Static insertion.}
For $a=m-1$, the degree $-n$ term is $TS_{a,2}$. Multiplication by $T$
on the left has no covariable derivatives, even when $T$ is not constant.
Since $n=2(a+1)$, its spherical average is
\[
 \frac1{\nu_{n-1}}\int_{S^{n-1}}S_{a,2}\,\dd S
 =-aE+\frac{a(a+1)}{3n}s
 =a(s/6-E).
\]
This proves~\eqref{eq:static}. If $m=1$, the operator is $T$ and its
residue is zero, agreeing with the factor $m-1$.

\paragraph{Second-order insertion.}
The symbols of $\nabla_i\nabla_j$ at $0$ are
$p_2=-\xi_i\xi_j$, $p_1=0$ and $p_0=\partial_i\omega_j=\Omega_{ij}/2$.
At degree $-n$ all terms in the composition are
\[
 p_2S_{m,2}+p_0S_{m,0}
 +\partial_{\xi_r}p_2D_{x_r}S_{m,1}
 +\tfrac12\partial_{\xi_r}\partial_{\xi_s}p_2
                         D_{x_r}D_{x_s}S_{m,0}.
\]
After the preceding substitutions this is
\begin{align*}
 &m\xi_i\xi_jE q^{-m-1}
 -\frac{m(m+1)}3\Ric_{ab}\xi_i\xi_j\xi_a\xi_bq^{-m-2}
 +\tfrac12\Omega_{ij}q^{-m}\\
 &\qquad
 +m(\xi_jV_{ib}\xi_b+\xi_iV_{jb}\xi_b)q^{-m-1}
 -mQ_{ij}q^{-m-1}.
\end{align*}
The second and fourth spherical moments, with $n=2m$, give respectively
\begin{align*}
 \text{potential: }&\tfrac12\delta_{ij}E,\\
 \text{quartic metric term: }&-\tfrac1{12}
                         (s\delta_{ij}+2\Ric_{ij}),\\
 \text{first differentiated product: }&\tfrac12(V_{ij}+V_{ji})
                                      =\tfrac23\Ric_{ij},\\
 \text{second differentiated product: }&-\tfrac13\Ric_{ij},\\
 \text{zero-order insertion: }&\tfrac12\Omega_{ij}.
\end{align*}
Adding the Ricci coefficients gives $-1/6+2/3-1/3=1/6$.
This proves exactly the displayed expression for $H_{ij}$.
At the center $\dd x=\dd V_g$; invariance restores the density
$\dd V_g$ in arbitrary coordinates. No additional density multiplier
is inserted into the left-symbol composition.

\paragraph{First-order insertion.}
For $\nabla_iL^{-m}$ the degree $-n$ terms are
\[
 i\xi_iS_{m,1}+\omega_iS_{m,0}
 +\partial_{\xi_r}(i\xi_i)D_{x_r}S_{m,0}.
\]
All three vanish at the
center, since $S_{m,1}=\omega_i=\partial_rq=0$. This proves the
pointwise vanishing with arbitrary left multiplier $T$ and completes
the symbol-composition proof.
\end{proof}

\begin{proposition}[Universal mixed insertion]\label{prop:universal}
Let $A,B$ be smooth zero-order endomorphisms. Then
\begin{equation}\label{eq:universal}
 \Wres(ADBD^{-n+1})
 =C_m\int_M\tau_\Lambda(Ac_iBc_jH_{ij})\,\dd V_g.
\end{equation}
In particular derivatives of $B$ do not contribute to this residue.
\end{proposition}
\begin{proof}
Since $D^{-n+1}=DL^{-m}$ modulo smoothing operators, at a normal-coordinate center
\[
 ADBD=Ac_iBc_j\nabla_i\nabla_j+
 Ac_i(\nabla_iB)c_j\nabla_j.
\]
Apply Lemma~\ref{lem:contractions}. Multiplication by $A$ is on the left and introduces no symbol derivatives. The vanishing statement for the first-order insertion eliminates the entire second term, not merely its integral.
\end{proof}

\paragraph{Expansion on a test section.}
The preceding operator identity follows from
\[
 D(BDs)=c_i(\nabla_iB)c_j\nabla_js
       +c_iB(\nabla_ic_j)\nabla_js
       +c_iBc_j\nabla_i\nabla_js.
\]
In a geodesic orthonormal frame $\nabla_ic_j=0$ at the point in question.
This is a choice of frame, not a condition on $B$. Applying the
first-order vanishing with left multiplier $Ac_i(\nabla_iB)c_j$ removes
its entire local contribution. Derivatives of $A$ do not arise because
$\partial_\xi A=0$ in the left-symbol convention.

\section{The four mixed residues}\label{sec:four}
Define $F_\ell$ by the following ordered pairs in
$F_\ell=\Wres(A_\ell D B_\ell D^{-n+1})$:
\begin{equation}\label{eq:pairs}
\begin{array}{c|c|c}
\ell&A_\ell&B_\ell\\\hline
1&\hc(X)&\hc(Y)\\
2&\hc(X)\hc(Y)&c(Z)c(W)\\
3&\hc(X)\hc(Y)c(Z)&c(W)\\
4&\hc(X)\hc(Y)&\hc(Z)\hc(W).
\end{array}
\end{equation}
\begin{theorem}\label{thm:four}
For $n=2m\geq2$, $F_\ell=C_m\int_M I_\ell\dd V_g$, where
\begin{align}
 I_1&=\frac{m-1}{12}s\,g(X,Y),\label{eq:I1}\\
 I_2&=(m-1)\left(\frac{s}{12}g(X,Y)g(Z,W)+\frac12R(X,Y,Z,W)\right),\label{eq:I2}\\
 I_3&=-\frac16g(X,Y)\Ric(Z,W)-\frac{m-2}{12}s\,g(X,Y)g(Z,W)
 -\frac{m-1}{2}R(X,Y,Z,W),\label{eq:I3}\\
 I_4&=-\frac{m-1}{12}s\,K(X,Y,Z,W).\label{eq:I4}
\end{align}
These identities also hold for the corresponding local densities.
\end{theorem}

\subsection{The first and fourth insertions}
For a smooth endomorphism $B$, the covariant Leibniz rule is
\[
 D(Bs)=c_i(\nabla_iB)s+c_iB\nabla_i s.
\]
If $B=\hc(Y)$, the mixed Clifford relation gives $c_iB=-Bc_i$,
whereas compatibility of the exterior connection gives
$\nabla_iB=\hc(\nabla_iY)$. Therefore
\[
 DB=-BD+c_i\hc(\nabla_iY),\qquad
 ADBD^{1-n}=-AB L^{1-m}+Ac_i\hc(\nabla_iY)D L^{-m}.
\]
The second term has zero local residue by the first-order assertion of
Lemma~\ref{lem:contractions}. This argument retains the derivatives of
$Y$ and then proves that they do not contribute; it does not choose a
parallel extension of $Y$.

In the static term, $s/6-E=-s/12+R_{abkl}c_ac_b\hc_k\hc_l/8$.
The trace of the mixed term against $AB=\hc(X)\hc(Y)$ is zero:
the ordinary Clifford trace contributes $-\delta_{ab}$, which kills
$R_{abkl}$. The remaining trace is
\[
 \tau_\Lambda\{AB(s/6-E)\}=-\tfrac{s}{12}g(X,Y).
\]
The minus sign in $DB=-BD+\cdots$ consequently gives
$I_1=(m-1)sg(X,Y)/12$.

For the fourth insertion $B=\hc(Z)\hc(W)$ contains two hatted
generators, so $c_iB=Bc_i$. The same Leibniz calculation now gives
\[
 ADBD^{1-n}=AB L^{1-m}+Ac_i(\nabla_iB)D L^{-m}.
\]
The first-order term again has zero local residue. The mixed potential
still contains only two unhatted generators, whose trace kills its
skew curvature indices. Finally
$\tau_\Lambda(AB)=K(X,Y,Z,W)$, giving~\eqref{eq:I4}.

\subsection{The second insertion and the Riemann term}
Put $A=\hc(X)\hc(Y)$ and $B=c(Z)c(W)$. We give a pointwise reduction,
so no passage from an integrated trace identity to a local identity is
needed. At a normal-coordinate center, let
$\mathcal S(B)=\sum_{ij}c_iBc_jH_{ij}$. The contraction defining
$I_2$ is $\tau_\Lambda(A\mathcal S(B))$. Moving the first $c_i$
cyclically inside the finite-dimensional fiber trace is legitimate,
and $A$ commutes with every $c_i$. Thus
\[
 I_2=\tau_\Lambda\left(AB\sum_{ij}c_jH_{ij}c_i\right).
\]
We now compute this last endomorphism. The metric part contributes
$(n-2)s/12$. For the curvature parts, use
$E=c_ac_b\Omega_{ab}/2$ and
$[\Omega_{ij},c_k]=c(R(e_i,e_j)e_k)$. Direct contraction gives
\begin{align*}
 \sum_i c_iEc_i&=(4-n)E-s,\\
 \sum_{ij}c_j\Omega_{ij}c_i&=-2E+s.
\end{align*}
For the first identity one can instead insert~\eqref{eq:E}: the scalar
part gives $-ns/4$. For $a\ne b$ the other part uses
\[
 \sum_ic_i(c_ac_b)c_i=(4-n)c_ac_b.
\]
For the second, move $c_i$ to the left of $\Omega_{ij}$:
$\sum c_jc_i\Omega_{ij}=-2E$ and
$\sum c_j[\Omega_{ij},c_i]=s$. The latter follows from
$\sum_iR_{ijil}=-\Ric_{jl}$ and the symmetric Ricci contraction.
Consequently
\[
 \sum_{ij}c_jH_{ij}c_i
 =\frac{n-2}{12}s+\frac{2-n}{2}E
 =(m-1)(s/6-E).
\]
Thus the static reduction holds locally, including arbitrary scalar
weights, and
\[
 I_2=(m-1)\tau_\Lambda\{AB(s/6-E)\}.
\]

Write $G=g(X,Y)$ and $H=g(Z,W)$. The scalar contribution is
$sGH/12$, since $\tau_\Lambda(AB)=-GH$. For the mixed term,
\begin{align*}
 \tau_\Lambda(\hc(X)\hc(Y)\hc_k\hc_l)
 &=G\delta_{kl}-X_kY_l+X_lY_k,\\
 \tau_\Lambda(c(Z)c(W)c_ic_j)
 &=H\delta_{ij}-Z_iW_j+Z_jW_i.
\end{align*}
Terms containing either Kronecker delta vanish against $R_{ijkl}$.
The remaining contraction is
\begin{align*}
 &R_{ijkl}(-Z_iW_j+Z_jW_i)(-X_kY_l+X_lY_k)\\
 &\quad=R(Z,W,X,Y)-R(Z,W,Y,X)\\
 &\qquad\quad-R(W,Z,X,Y)+R(W,Z,Y,X)\\
 &\quad=4R(X,Y,Z,W).
\end{align*}
Multiplication by $1/8$ gives $R(X,Y,Z,W)/2$. This proves~\eqref{eq:I2}
with its sign and coefficient fixed by the two skew pairs.

\subsection{The third insertion}
Decompose the endomorphisms as
\[
 E=s/4+E_{\rm mix},\quad
 E_{\rm mix}=-\tfrac18R_{abkl}c_ac_b\hc_k\hc_l,\quad
 \Omega^c_{ij}=\tfrac14R_{ijkl}c_kc_l,\quad
 \Omega^{\widehat c}_{ij}=-\tfrac14R_{ijkl}\hc_k\hc_l.
\]
Then $H_{ij}=\Ric_{ij}/6+s\delta_{ij}/24
+\delta_{ij}E_{\rm mix}/2+(\Omega^c_{ij}+\Omega^{\widehat c}_{ij})/2$.
We compute the four resulting contributions separately.

For any symmetric matrix $S$, the Clifford relations imply
\[
 \sum_{ij}S_{ij}c_ic(W)c_j=(\tr S)c(W)-2c(SW).
\]
Indeed $c_ic(W)=-c(W)c_i-2W_i$, and
$\sum S_{ij}c_ic_j=-\tr S$. With $S=\Ric/6+s\Id/24$,
this gives $(n+2)sc(W)/24-c(\Ric^\sharp W)/3$.

The curvature identity
$\sum_jc_j\Omega^c_{ij}=-c(\Ric^\sharp e_i)/2$
follows by expanding $R_{ijkl}c_jc_kc_l/4$: its three-generator
antisymmetric part vanishes by Bianchi, and its contractions give
the displayed Ricci term. In particular $\sum c_ic_j\Omega^c_{ij}=s/2$.
Moving $c(W)$ past $c_i$ now yields
\[
 \frac12\sum_{ij}c_ic(W)c_j\Omega^c_{ij}
 =-\frac{s}{4}c(W)+\frac12c(\Ric^\sharp W).
\]
Adding these two contributions proves
\begin{equation}\label{eq:spincontraction}
 \sum_{ij}c_ic(W)c_j
 \left(\frac{\Ric_{ij}}6+\frac{s\delta_{ij}}{24}
       +\frac12\Omega^c_{ij}\right)
 =\frac16c(\Ric^\sharp W)+\frac{n-4}{24}s\,c(W).
\end{equation}
Taking the trace after $\hc(X)\hc(Y)c(Z)$ gives
\[
 -\frac16g(X,Y)\Ric(Z,W)
 -\frac{n-4}{24}s\,g(X,Y)g(Z,W).
\]

For the mixed potential, use
$\sum_ic_ic(W)c_i=(n-2)c(W)$, obtained from
$c_ic(W)c_i=c(W)-2W_ic_i$. Its contribution is
\begin{align*}
 &-\frac1{16}R_{abkl}\sum_i\tau_\Lambda
 (\hc(X)\hc(Y)c(Z)c_ic(W)c_ic_ac_b\hc_k\hc_l)\\
 &\qquad=-\frac{n-2}{16}\cdot4R(X,Y,Z,W)
 =-\frac{n-2}{4}R(X,Y,Z,W).
\end{align*}
The factor four is exactly the skew contraction calculated above.
Finally,
\[
 \tau_\Lambda(c(Z)c_ic(W)c_j)
 =Z_iW_j-g(Z,W)\delta_{ij}+Z_jW_i
\]
is symmetric in $i,j$. Its contraction with $R_{ijkl}$ is zero, so
$\Omega^{\widehat c}$ contributes nothing. Since $n=2m$, the sum is
precisely~\eqref{eq:I3}. Every step is a local contraction; no derivative
of an inserted field has been discarded without the vanishing statement
in Lemma~\ref{lem:contractions}.

\section{Curvature recovery and classical checks}\label{sec:curvature}
For $n>2$ define
\begin{equation}\label{eq:Rfunctional}
 \RR_D(u,v,w,z)=\frac1{n-2}\Wres
 \left([\hc(u),\hc(v)][c(w),c(z)]|D|^{2-n}\right).
\end{equation}
Here vectors and covectors are identified by $g$. The factor is $n-2=2(m-1)$.
\begin{corollary}\label{cor:R}
On a closed Riemannian manifold,
\begin{equation}
 \RR_D(X,Y,Z,W)=C_m\int_MR(X,Y,Z,W)\dd V_g.
\end{equation}
With an arbitrary scalar weight $f$ inserted on the left, the same formula holds with $f$ in the integrand. Consequently the weighted family determines the pointwise curvature tensor.
\end{corollary}
\begin{proof}
Antisymmetrize the static form of $F_2$ in both pairs. The scalar term cancels and the Riemann term is multiplied by four. Divide by $2(m-1)$.
\end{proof}
The recovered classical tensor is skew in each pair, symmetric under interchange of the pairs, and obeys the first Bianchi identity. These follow from the recovered tensor, not from a claim that arbitrary noncommutative commutators satisfy the same identities.

Another useful check is
\begin{equation}\label{eq:Einstein}
 F_2+F_3=-\frac{C_m}{6}\int_Mg(X,Y)
 \left(\Ric(Z,W)-\frac s2g(Z,W)\right)\dd V_g.
\end{equation}
In dimension four, $C_2=2/\pi^2$ and
\[
\begin{aligned}
 I_2&=s\,g(X,Y)g(Z,W)/12+R(X,Y,Z,W)/2,\\
 I_3&=-g(X,Y)\Ric(Z,W)/6-R(X,Y,Z,W)/2.
\end{aligned}
\]
In dimension two all four densities vanish, using $\Ric=sg/2$. Accordingly~\eqref{eq:Rfunctional} is not defined at $n=2$; this residue order does not recover surface curvature. Surface curvature can instead be recovered by a different scalar residue or heat coefficient. Removing the factor $n-2$ by fiat would not repair this degeneracy.

For constant sectional curvature $K$, substitution of
$R(X,Y,Z,W)=K(g(X,W)g(Y,Z)-g(X,Z)g(Y,W))$
checks the sign of the mixed terms. On a flat manifold all four densities vanish for arbitrary smooth inserted vector fields. This last statement is specific to these zero-order Clifford insertions and is not true for arbitrary products of direction covariant derivatives.

\subsection{Recovery by antisymmetrizing the first two inputs}
To make the four arguments explicit, write $I_j(X,Y;Z,W)$ for the densities in Theorem~\ref{thm:four}. The variables $Z,W$ are held fixed throughout this subsection. Since $g(X,Y)=g(Y,X)$ and $R(Y,X,Z,W)=-R(X,Y,Z,W)$, direct subtraction gives
\begin{align}
 I_2(X,Y;Z,W)-I_2(Y,X;Z,W)&=(m-1)R(X,Y,Z,W),\\
 I_3(X,Y;Z,W)-I_3(Y,X;Z,W)&=-(m-1)R(X,Y,Z,W).
\end{align}
The cancellation in the second line includes both the scalar-curvature term and the Ricci term: each is multiplied by the symmetric factor $g(X,Y)$. Consequently, for $m>1$,
\begin{align}
 R(X,Y,Z,W)&=\frac{I_2(X,Y;Z,W)-I_2(Y,X;Z,W)}{m-1},\\
 R(X,Y,Z,W)&=-\frac{I_3(X,Y;Z,W)-I_3(Y,X;Z,W)}{m-1}.
\end{align}
Thus either insertion independently recovers the full Riemann tensor. Antisymmetrizing both pairs is convenient for an algebraic definition but is not necessary for classical recovery.

For completeness let $F_{j,f}$ denote the same residue with an arbitrary smooth scalar function $f$ inserted on the far left. There are no derivatives of $f$ in a left multiplication. The operator-level statements are
\begin{align}
 F_{2,f}(X,Y;Z,W)-F_{2,f}(Y,X;Z,W)
 &=C_m(m-1)\int_M fR(X,Y,Z,W)\,\dd V_g,\\
 F_{3,f}(X,Y;Z,W)-F_{3,f}(Y,X;Z,W)
 &=-C_m(m-1)\int_M fR(X,Y,Z,W)\,\dd V_g.
\end{align}
If these weighted expressions vanish for all $f$ supported in a coordinate neighborhood, the fundamental lemma for distributions implies that the corresponding curvature component vanishes there. Taking local frame inputs and varying all four indices therefore determines the complete smooth curvature tensor, rather than only its integral. At $m=1$ the prefactor is zero, so division is unavailable and no two-dimensional recovery follows from these formulas.

\section{Boundary residues in general even dimension}\label{sec:boundary}
\subsection{The operator whose boundary residue is being computed}
On a manifold with boundary, the closed expression alone does not specify a boundary functional. Fix $n=2m\geq4$, extend the data smoothly across the boundary, and use the class-zero transmission operators
\begin{equation}\label{eq:boundarydef}
 \mathcal B_\ell=
 \pi^+(A_\ell D B_\ell D^{-2})\circ\pi^+(D^{-n+3}).
\end{equation}
The factors have orders $-1$ and $3-n$. Their interior product is $A_\ell D B_\ell D^{-n+1}$. The boundary residue below is the Fedosov--Golse--Leichtnam--Schrohe residue \cite{FGLS,Schrohe,WangBoundary}, scaled by $(2\pi)^{-n}$ to agree with~\eqref{eq:Wres}. No APS, absolute, or relative realization of $D$ is intended. Different factorizations of the same interior operator can have different singular Green remainders and different boundary residues.

Use an inward normal coordinate $x_n\geq0$ and assume
\begin{equation}\label{eq:collar}
 g=h(x_n)^{-1}g_{\partial M}+\dd x_n^2,\qquad h(0)=1,\qquad h'(0)=\eta.
\end{equation}
At the point of calculation choose geodesic boundary coordinates and the orthonormal frame $e_a=\sqrt h\,\partial_a$, $a<n$, $e_n=\partial_{x_n}$. Trivialize the exterior bundle in this orthonormal frame. The matrices $c_i,\hc_i$ are constant in this trivialization. Write
\begin{equation}
 u=\xi',\quad z=\xi_n,\quad |u|=1,\quad
 U=c(u),\quad N=c_n,\quad C=U+zN,\quad q=1+\xi_n^2,\quad k=m-1.
\end{equation}
The full expression before restriction is $q(x,\xi)=h(x_n)|\xi'|^2+\xi_n^2$. In particular
\begin{equation}\label{eq:jets}
 q_{x_n}={\eta}|u|^2,\quad C_{x_n}=\tfrac{\eta}2U,
 \quad q_{x_a}=C_{x_a}=0\quad(a<n)
\end{equation}
at the chosen point. The restriction $|u|=1$ is made only after these derivatives.

\subsection{Connection and vector-field jets}
The nonzero frame connection coefficients at the center are
\begin{equation}
 \nabla_{e_a}e_b=\tfrac{\eta}2\delta_{ab}e_n,
 \qquad\nabla_{e_a}e_n=-\tfrac{\eta}2e_a,
 \qquad\nabla_{e_n}e_j=0.
\end{equation}
Therefore the order-zero Dirac symbol is
\begin{equation}\label{eq:w}
 w=\sigma_0(D)=\frac{\eta}4
 \left[-(n-1)N-\sum_{a<n}c_a\hc_a\hc_n\right].
\end{equation}
The factor $1/4$ includes both skew connection entries. For a vector field $V=V^je_j$, its ordinary coefficient derivative in this trivialization is
\begin{equation}\label{eq:Vjets}
 V_{,a}=\nabla_{e_a}V-\tfrac{\eta}2(V_ae_n-V_ne_a),\qquad
 V_{,n}=\nabla_{e_n}V.
\end{equation}
Consequently $\partial_jc(V)=c(V_{,j})$ and $\partial_j\hc(V)=\hc(V_{,j})$. Products are differentiated by the ordinary Leibniz rule in their displayed order. This convention retains arbitrary first jets; it does not set the vector-field derivatives to zero.

\paragraph{Ordinary coefficients versus covariant derivatives.}
The connection formulas are
\[
 \nabla_{e_a}e_b=\tfrac\eta2\delta_{ab}N,\qquad
 \nabla_{e_a}N=-\tfrac\eta2e_a.
\]
Differentiating $V=V_be_b+V_nN$ therefore gives
\begin{align}
 (\nabla_{e_a}V)_b&=e_a(V_b)-\tfrac12\eta{}V_n\delta_{ab},\\
 (\nabla_{e_a}V)_n&=e_a(V_n)+\tfrac12\eta{}V_a.
\end{align}
Solving these identities for the coefficient derivatives gives~\eqref{eq:Vjets}. Their tangential trace gives
$\mathcal D(V)=\operatorname{div}_{\partial M}V^T-(n-1)\eta{}V_n/2$.
Taking $V=N$ gives the mean-curvature conversion~\eqref{eq:meanconversion}. Thus replacing $\eta$ by $K_{\partial}$ is the final geometric step, after all coefficient-derivative corrections have been included.

\subsection{All homogeneous symbols that can contribute}
The principal and subprincipal inverse-square symbols are
\begin{align}
 \sigma_{-2}(D^{-2})&=q^{-1},\\
 \ell_3:=\sigma_{-3}(D^{-2})
 &=-i(Cw+wC+NC_{x_n})q^{-2}-2i{\eta} zq^{-3}.\label{eq:ell3}
\end{align}
Indeed $\sigma_1(D^2)=i(Cw+wC+NC_{x_n})$ and the composition correction is $-q^{-1}\partial_{\xi_n}qD_tq^{-1}=-2i{\eta} zq^{-3}$.
For the inverse powers of $L$ one has
\begin{equation}\label{eq:Lksub}
 \sigma_{-2k}(L^{-k})=q^{-k},\qquad
 \sigma_{-2k-1}(L^{-k})
 =kq^{-k+1}\ell_3-i k(k-1){\eta} zq^{-k-2}.
\end{equation}
To see the second term, choose the one differentiated pair among $k$ leading factors. There are $k(k-1)/2$ choices and
$\partial_{\xi_n}q^{-1}D_tq^{-1}=-2i{\eta} zq^{-4}$.

For the differentiated-pair term in~\eqref{eq:Lksub}, label the $k$
leading inverse factors by $1,\ldots,k$. A first-order composition
correction differentiates factor $r$ in $\xi_n$ and factor $s$ in
$x_n$ with $r<s$. For each $s=2,\ldots,k$ there are $s-1$ possible
left factors, so the number is $1+\cdots+(k-1)=k(k-1)/2$.
The remaining $k-2$ factors supply $q^{-k+2}$; multiplication by
$-2i\eta zq^{-4}$ gives exactly
$-ik(k-1)\eta zq^{-k-2}$. The order of the two positions is dictated
by the left-symbol convention; reversing the positions does not create
another term.

For one of the pairs $(A,B)$ in~\eqref{eq:pairs}, set
$P=ADB D^{-2}$ and $Q=D^{-n+3}=DL^{-k}$. Define
\begin{align}
 p&:=\sigma_{-1}(P)=iACBq^{-1},\label{eq:p}\\
 p_-&:=\sigma_{-2}(P)
 =iACB\ell_3+A\left(wB+\sum_jc_jB_{,j}\right)q^{-1}
 -{\eta} ANBq^{-2},\label{eq:psub}\\
 p_{x_n}&=i(A_{,n}CB+AC_{x_n}B+ACB_{,n})q^{-1}
       -i{\eta} ACBq^{-2},\label{eq:pt}\\
 t_0&:=\sigma_{3-n}(Q)=iCq^{-k},\\
 t_-&:=\sigma_{2-n}(Q)
 =ikCq^{-k+1}\ell_3+k(k-1){\eta} zCq^{-k-2}
 +wq^{-k}-k{\eta} Nq^{-k-1},\label{eq:tsub}\\
 t_{x_n}:=\partial_{x_n}t_0&=iC_{x_n}q^{-k}-ik{\eta} Cq^{-k-1}.\label{eq:tt}
\end{align}
Every term is explicit in ${\eta}$, the Clifford matrices, vector-field values, and first derivatives. For example the last term of~\eqref{eq:psub} is
$\partial_{\xi_n}(iACB)D_tq^{-1}=-{\eta} ANBq^{-2}$;
its sign would be lost if the symbol product were replaced by pointwise multiplication.

\paragraph{Deriving the left factor and its normal derivative.}
For a test section $s$, first expand $ADB$ rather than reorder its
endomorphisms:
\[
 ADBs=Ac_jB\partial_js+A(c_jB_{,j}+wB)s.
\]
Thus its degree-one symbol is $iACB$ and its degree-zero symbol is
$A(c_jB_{,j}+wB)$. Multiplication with $D^{-2}$ has, at order $-2$,
exactly the three terms
\[
 (iACB)\ell_3+A(c_jB_{,j}+wB)q^{-1}
       +\partial_{\xi_n}(iACB)D_{x_n}q^{-1}.
\]
Since $D_{x_n}q^{-1}=i\eta q^{-2}$ on $|\xi'|=1$, the last term is
$-\eta ANBq^{-2}$. This proves~\eqref{eq:psub}. Ordinary differentiation
of $p=iACBq^{-1}$ gives~\eqref{eq:pt}, including both $A_{,n}$ and
$B_{,n}$. Similarly $t_{x_n}$ means the derivative of $t_0$ with the
covariable held fixed, not a new homogeneous symbol. Explicitly,
\[
 \partial_{x_n}(iCq^{-k})
 =iC_{x_n}q^{-k}-ikCq^{-k-1}q_{x_n}.
\]
These derivatives are taken before imposing $|\xi'|=1$. Here
$D_{x_n}=(1/i)\partial_{x_n}$ is the symbol derivative; it is not the
Dirac operator $D$.

\subsection{Five possible cases}
We write $\Phi_\ell$ for the boundary density before the exterior-rank and tangential-sphere prefactor, and $\phi_{\ell,j}$ for its five normalized case contributions. The notation $\pi_{\xi_n}^+$ denotes the Hardy projection in the normal covariable; $\pi^+$ on operators denotes truncation. These are the conventions of the cited boundary computations. The fixed power parameter is $k=m-1$; we use $v$ for the normal-derivative counter in the composition formula to avoid confusing these two integers.

The boundary composition formula requires
\begin{equation}
 r+\ell-j-v-|\alpha|-1=-n,\quad r\leq-1,\quad\ell\leq3-n.
\end{equation}
Writing $r=-1-a$, $\ell=3-n-b$ gives
$a+b+j+v+|\alpha|=1$. Thus there are exactly five types:
\begin{center}
\begin{tabular}{llll}\toprule
Case & Nonzero loss & Contribution before trace and integration & Status\\\midrule
(a.I) & $|\alpha|=1$ & $-\partial_{\xi'}\pi_{\xi_n}^+p\,\partial_{x'}\partial_{\xi_n} t_0$ & $0$\\
(a.II) & $j=1$ & $-\tfrac12\pi_{\xi_n}^+p_{x_n}\,\partial_{\xi_n}^2t_0$ & retained\\
(a.III) & $v=1$ & $-\tfrac12\partial_{\xi_n}\pi_{\xi_n}^+p\,\partial_{\xi_n}t_{x_n}$ & retained\\
(b) & $a=1$ & $-i\pi_{\xi_n}^+p_-\,\partial_{\xi_n}t_0$ & retained\\
(c) & $b=1$ & $-i\pi_{\xi_n}^+p\,\partial_{\xi_n}t_-$ & retained\\\bottomrule
\end{tabular}
\end{center}
The first case vanishes because the right principal symbol has zero tangential derivative at the coordinate center. It does not require tangential constancy of the inserted fields.

\subsection{Finite evaluation of the normal and angular variables}
We now remove all contour integrals from the answer. For integers $p,b,d\geq0$ with $b\geq1$, define
\begin{equation}\label{eq:Tcoef}
 T(p,b,d)=\sum_{j=0}^{\min(p,d)}
 \binom pj i^{p-j}(-1)^{d-j}
 \binom{b+d-j-1}{d-j}(2i)^{-b-d+j}.
\end{equation}
This is the coefficient of $(\xi_n-i)^d$ in $\xi_n^p(\xi_n+i)^{-b}$. For $a,b\geq1$ put
\begin{equation}\label{eq:Jcoef}
 \mathcal J(a,p;b,q)
 =2\pi i\sum_{v=1}^{a}T(p,a,a-v)T(q,b,v+b-1).
\end{equation}
Then, with the upper-half-plane Hardy projection,
\begin{equation}\label{eq:Jidentity}
 \int_{\mathbb R}\pi_{\xi_n}^+\!\left(\frac{\xi_n^p}{(1+\xi_n^2)^a}\right)
 \frac{\xi_n^q}{(1+\xi_n^2)^b}\,\dd \xi_n=\mathcal J(a,p;b,q)
\end{equation}
whenever the resulting product is integrable. In the applications here it always is. To prove the identity, the principal part of the first factor at $i$ is
$\sum_{v=1}^a T(p,a,a-v)(\xi_n-i)^{-v}$.
Multiply by the second factor and read the coefficient of $(\xi_n-i)^{-1}$. The lower pole contributes no enclosed residue. Polynomial parts, if present, are omitted by this projection.

\paragraph{Projection and the residue theorem.}
For the rational functions used here, the Hardy projection in the
boundary convention of \cite{WangGravity,WWJMP} retains the principal
part at the upper pole $z=i$. For example,
\[
 (1+z^2)^{-1}=\frac1{2i(z-i)}-\frac1{2i(z+i)},\qquad
 \pi_z^+(1+z^2)^{-1}=\frac1{2i(z-i)}.
\]
For higher poles write $z=i+t$ and expand the analytic factor
$(2i+t)^{-a}$; the negative powers of $t$ are exactly the retained
principal part. To integrate the product, close the real axis by a
counterclockwise upper semicircle. The rational products in the case
formulas decay at least as $|z|^{-2}$, so the arc integral tends to zero.
The residue theorem then gives $2\pi i$ times the residue at $i$.
A pole of order $r$ contributes
\[
 \Res_{z=i}\frac{F(z)}{(z-i)^r}
 =\frac{F^{(r-1)}(i)}{(r-1)!}.
\]
Thus the real-line integral is converted to a residue by this contour
argument, not by identifying integration with the projection itself.
In~\eqref{eq:Tcoef}, expanding
$(i+t)^p=\sum_j\binom pj i^{p-j}t^j$ and
$(2i+t)^{-b}=\sum_l(-1)^l\binom{b+l-1}{l}(2i)^{-b-l}t^l$
and setting $j+l=d$ gives every coefficient in the displayed finite sum.

For a tangential multiindex $\beta\in\mathbb N^{n-1}$, let
\begin{equation}\label{eq:mom}
 \mathcal M_{n-1}(\beta)=
 \begin{cases}
 0,&\text{some $\beta_j$ is odd},\\
 \displaystyle\frac{\prod_j(\beta_j-1)!!}
 {(n-1)(n+1)\cdots(n+2s-3)},&|\beta|=2s,
 \end{cases}
\end{equation}
with empty products equal to $1$. This is the spherical average of $u^\beta$. It follows by the Gaussian moment formula and division by the radial moment.

Every symbol in~\eqref{eq:p}--\eqref{eq:tt} is a finite sum of terms
$F_{a,p,\beta}u^\beta \xi_n^p/q^a$ with numerical Clifford coefficients $F$. Define the finite bilinear contraction
\begin{equation}\label{eq:bracket}
 \mathcal C(F,G)=\sum_{a,p,\beta;b,q,\gamma}
 \tau_\Lambda(F_{a,p,\beta}G_{b,q,\gamma})
 \mathcal M_{n-1}(\beta+\gamma)\mathcal J(a,p;b,q).
\end{equation}
There is no pseudodifferential symbol or unevaluated integral on the right. Differentiation in $\xi_n$ is performed by
\begin{equation}\label{eq:dzrule}
 \partial_{\xi_n}\left(\frac{\xi_n^p}{q^a}\right)
 =p\frac{\xi_n^{p-1}}{q^a}-2a\frac{\xi_n^{p+1}}{q^{a+1}}.
\end{equation}
These formulas use finite sums for every integer $m$; they are not a low-dimensional interpolation.

\begin{theorem}[Boundary residue in even dimension]\label{thm:boundary}
Under~\eqref{eq:collar}, for every even $n=2m\geq4$ and each pair in~\eqref{eq:pairs},
\begin{equation}\label{eq:boundaryfinal}
 \widetilde{\Wres}(\mathcal B_\ell)
 =C_m\int_M I_\ell\dd V_g
 +\frac{r\nu_{n-2}}{(2\pi)^n}\int_{\partial M}\Phi_\ell\dd V_{\partial M},
\end{equation}
where $I_\ell$ is given in Theorem~\ref{thm:four} and the fully finite boundary coefficient is
\begin{equation}\label{eq:Qfinite}
 \boxed{\Phi_\ell=
 -\tfrac12\mathcal C(p_{x_n},\partial_{\xi_n}^2t_0)
 -\tfrac12\mathcal C(\partial_{\xi_n}p,\partial_{\xi_n}t_{x_n})
 -i\mathcal C(p_-,\partial_{\xi_n}t_0)
 -i\mathcal C(p,\partial_{\xi_n}t_-).}
\end{equation}
Equivalently, the fully geometric expressions are given in equations~\eqref{eq:geomsums12}--\eqref{eq:geomsum4} below. All the inputs on the right are the explicit rational polynomials
\eqref{eq:w}--\eqref{eq:tt}; $\mathcal C$ is the binomial-and-moment sum~\eqref{eq:bracket}. Thus the formula includes all first jets of arbitrary $X,Y,Z,W$ and requires no contour integration or symbol recursion at the point of use.
\end{theorem}
\begin{proof}
The interior symbol is the same as in Theorem~\ref{thm:four}. The boundary order equation gives precisely the five cases displayed above. Their coefficients follow from
$(-i)^{|\alpha|+j+v+1}/(\alpha!(j+v+1)!)$.
Insert~\eqref{eq:p}--\eqref{eq:tt}, apply~\eqref{eq:Jidentity} termwise, and apply~\eqref{eq:mom}. The matrix trace contributes $r$ times its normalized value; the tangential integration contributes $\nu_{n-2}$. These are exactly the prefactors in~\eqref{eq:boundaryfinal}.
\end{proof}

\paragraph{From the symbol density to the geometric boundary measure.}
The boundary formula is initially expressed in coordinates as a
tangential sphere integral followed by $\dd x'$. At a geodesic boundary
center $g_{\partial,ab}=\delta_{ab}$, so
$\dd x'=\dd V_{\partial M}$ there. Invariant assembly restores
$\sqrt{\det g_\partial}\,\dd x'$ away from that center. Normal integration
is in the covariable $\xi_n$, not the base coordinate $x_n$; hence its
residue contributes a coefficient to the boundary density rather than
an extra factor $\dd x_n$.

\subsection{Geometric notation and evaluation of the dimension coefficients}
We now reduce Theorem~\ref{thm:boundary} to geometric contractions. This is the form of the boundary answer used in the five-case calculations of \cite{WWJMP}, with coefficients recomputed for the exterior operator and the present factorization. Put
\begin{equation}\label{eq:boundaryscale}
 k=m-1,\quad d=n-1=2k+1,\quad
 \beta_k=\frac{\binom{2k}{k}}{4^k},\qquad
 b_n=\frac{2^n\nu_{n-2}\pi\beta_k}{(2\pi)^n}.
\end{equation}
Write $N=e_n$, $V_n=g(V,N)$, $V^T=V-V_nN$, and $\dot V=\nabla_NV$. All derivatives below are evaluated at the boundary. Define
\begin{equation}
 G=g(X,Y),\qquad H=g(Z,W),\qquad P_1=Z_nW_n.
\end{equation}
\begin{align}
 J&=g(X,Y)g(Z,W)-g(X,Z)g(Y,W)+g(X,W)g(Y,Z),\\
 \mathcal D(V)&=\sum_{a<n}g(\nabla_{e_a}V,e_a),\\
 T_-&=g(\nabla_{Z^T}W,N)-g(\nabla_{W^T}Z,N)
       -Z_n\mathcal D(W)+W_n\mathcal D(Z),\\
 T_+&=g(\nabla_{Z^T}W,N)+Z_n\mathcal D(W),\\
 J_B&=J(X,Y,\dot Z,W)+J(X,Y,Z,\dot W),\\
 L&=W_nX_n g(Y,Z)-W_nY_n g(X,Z)
       +Z_nX_n g(Y,W)-Z_nY_n g(X,W).
\end{align}
Here $J(U,V,S,T)$ denotes the same four-linear metric polynomial as $J$. A dot on a scalar contraction is its normal derivative. For example,
\begin{align}
 \dot G&=g(\dot X,Y)+g(X,\dot Y),\qquad
 \dot H=g(\dot Z,W)+g(Z,\dot W),\\
 \dot P_1&=g(\dot Z,N)W_n+Z_ng(\dot W,N).
\end{align}
The last identity uses $\nabla_NN=0$ in the collar. The tangential divergence $\mathcal D(V)$ is the displayed ambient-connection contraction, not an unannounced boundary divergence. In particular
\[
 \mathcal D(V)=\operatorname{div}_{\partial M}V^T-\frac{d{\eta}}{2}V_n.
\]
The sign convention for extrinsic curvature used in this section is fixed by
\begin{equation}\label{eq:meanconversion}
 K_{\partial}:=\sum_{a<n}g(\nabla_{e_a}N,e_a)
 =-\frac d2{\eta},\qquad {\eta}=-\frac{2K_{\partial}}d.
\end{equation}
The equality itself fixes the sign: $N$ is the inward coordinate normal, regardless of alternative names for the extrinsic-curvature convention.

For clarity, let $\phi_{\ell,j}$ denote the contribution of case $j$, divided by $\pi\beta_k$, before multiplication by $2^n\nu_{n-2}/(2\pi)^n$. Thus
\begin{equation}
 \Phi_\ell=\pi\beta_k\sum_{j=1}^5\phi_{\ell,j}.
\end{equation}
The following six coefficients will occur repeatedly:
\begin{align}\label{eq:abcde}
 a_k&=\frac{3k(2k+1)}{8(k+1)(k+2)},&
 b_k&=\frac{k}{2(k+1)},& t_k&=\frac1{2(k+1)},\notag\\
 c_k&=\frac{k(2k+1)^2}{8(k+1)(k+2)},&
 d_k&=\frac{k(4k^2+16k+7)}{8(k+1)(k+2)},&
 e_k&=\frac{k(2k+1)(4k+5)}{8(k+1)(k+2)}.
\end{align}
They satisfy $a_k-c_k-d_k+e_k=0$.

Here is an explicit reduction of the contour coefficients, so that the rational functions in the answer do not rely on dimensional interpolation. With $a,p,q$ fixed and $b=k+o$, equation~\eqref{eq:Jcoef} becomes
\begin{align}\label{eq:binomialreduced}
 \frac{\mathcal J(a,p;b,q)}{\pi\beta_k}
 &=i^{p+q}\sum_{v=1}^{a}\sum_{j=0}^{\min(p,a-v)}\sum_{l=0}^{q}
 (-1)^{j+l}2^{2-2a-2o+j+l}\binom pj\binom ql\notag\\
 &\quad\times\binom{2a-v-j-1}{a-v-j}
 \frac{(2b+v-l-2)!(k!)^2}
 {(b+v-l-1)!(b-1)!(2k)!}.
\end{align}
Terms outside the valid binomial range are zero. In every application, cancellation of factorial ratios yields the rational functions displayed below. This is also a useful way to avoid ambiguous empty-product conventions in low dimensions.

The auxiliary integer $o$ in $b=k+o$ is a fixed shift of the denominator
exponent, not an additional dimension or a derivative index. For
example, a factor $q^{-k-2}$ has $o=2$. Substituting $b=k+o$ isolates
the central binomial factor $\beta_k$; factorial ratios with fixed
shifts then reduce to rational functions of $k$. This is the purpose
of~\eqref{eq:binomialreduced}.

\paragraph{Tangential and normal decomposition.}
At $x_0\in\partial M$, write $V=V^T+V_nN$ and $\xi=\xi'+\xi_n\dd x_n$. Rotational invariance gives
\begin{align}
 \int_{|\xi'|=1}\xi_a\xi_b\,\sigma(\xi')&=\frac{\nu_{n-2}}{n-1}\delta_{ab},\\
 \int_{|\xi'|=1}g(V,\xi')g(W,\xi')\,\sigma(\xi')
 &=\frac{\nu_{n-2}}{n-1}\{g(V,W)-V_nW_n\}.
\end{align}
The first identity follows because the off-diagonal moments change sign under a coordinate reflection, the diagonal moments are equal, and their sum is the sphere volume. The second is obtained by expanding the two tangential components. All odd tangential moments vanish by $\xi'\mapsto-\xi'$.

\paragraph{Normalization of the detailed case integrands.}
For each insertion $\ell$ and case $j$, let $\mathfrak I_{\ell,j}(\xi_n)$ denote the normalized angular average of the full traced case integrand, with the trace divided by $2^n$. It includes the factor $-1/2$ or $-i$ appearing in the boundary composition formula. Consequently
\[
 \phi_{\ell,j}=\frac1{\pi\beta_k}\int_{\mathbb R}\mathfrak I_{\ell,j}(\xi_n)\dd\xi_n.
\]
The geometric symbols $M_\rho$ and polynomials $P_\rho$ below are local to their case; their indices restart for each case. No Clifford trace or sphere integral remains in these lists.

The lists are obtained by applying the projection before multiplying the two factors. In particular, for $a=1,2,3$ the necessary projections are given by the explicit principal-part formula
\begin{equation}
 \pi_{\xi_n}^+\frac{\xi_n^p}{(1+\xi_n^2)^a}
 =\sum_{v=1}^a\sum_{j=0}^{\min(p,a-v)}
 \binom pj i^{p-j}(-1)^{a-v-j}
 \binom{2a-v-j-1}{a-v-j}
 \frac{(2i)^{-2a+v+j}}{(\xi_n-i)^v}.
\end{equation}
For example,
\begin{align}
 \pi_{\xi_n}^+q^{-1}&=\frac1{2i(\xi_n-i)},&
 \pi_{\xi_n}^+(\xi_nq^{-1})&=\frac1{2(\xi_n-i)},\\
 \pi_{\xi_n}^+q^{-2}&=-\frac1{4(\xi_n-i)^2}-\frac{i}{4(\xi_n-i)},&
 \pi_{\xi_n}^+(\xi_nq^{-2})&=-\frac{i}{4(\xi_n-i)^2}.
\end{align}
These identities follow by setting $\xi_n=i+t$, expanding $(2i+t)^{-a}$, and keeping precisely the negative powers of $t$. Projection of the product of the two symbols would give a different operation and is not used.

For a polynomial $P$, the last step in every list is Cauchy's formula
\begin{equation}\label{eq:expandedcauchy}
 \int_{\mathbb R}\frac{P(\xi_n)\dd\xi_n}
 {(\xi_n-i)^{k+a}(\xi_n+i)^{k+b}}
 =\frac{2\pi i}{(k+a-1)!}
 \left[\partial_{\xi_n}^{k+a-1}
 \frac{P(\xi_n)}{(\xi_n+i)^{k+b}}\right]_{\xi_n=i}.
\end{equation}
For $P(\xi_n)=\xi_n^p$, evaluation of the derivative gives the completely explicit normalized result
\begin{align}
 &\frac1{\pi\beta_k}\int_{\mathbb R}
 \frac{\xi_n^p\dd\xi_n}{(\xi_n-i)^{k+a}(\xi_n+i)^{k+b}}\notag\\
 &\quad=i^{p+2-a-b}\sum_{j=0}^{p}(-1)^{a-1-j}2^{2-a-b+j}\binom pj
 \frac{(2k+a+b-j-2)!(k!)^2}
 {(k+a-j-1)!(k+b-1)!(2k)!}.
\end{align}
Terms with $j>k+a-1$ are zero. The individual values below have already simplified these factorial ratios. Thus both the projected rational integrands and their evaluated residues are displayed for each nonzero case.

\subsection{The first insertion: all five cases}
For $(A,B)=(\hc(X),\hc(Y))$, the two-principal-symbol trace is
\[
 \tau_\Lambda(\hc(X)c(v)\hc(Y)c(w))=Gg(v,w).
\]
Consequently only the tangential contraction $|u|^2=1$ and the normal contraction $g(N,N)=1$ survive. The mixed contractions $g(u,N)$ vanish.

\paragraph{Case (a.I).}

\[
 \Phi_{1,1}=-\int_{|\xi'|=1}\int_{\mathbb R}
 \sum_{a<n}\tr_\Lambda\left[
 \partial_{\xi_a}\pi_{\xi_n}^+\sigma_{-1}(P)
 \partial_{x_a}\partial_{\xi_n}\sigma_{3-n}(Q)\right]
 \dd\xi_n\,\sigma(\xi').
\]
At $x_0$, $\partial_{x_a}g^{ij}=0$ for $a<n$, and the chosen orthonormal trivialization gives $\partial_{x_a}c(\xi)=0$. Therefore
\[
 \partial_{x_a}\sigma_{3-n}(Q)
 =i\,\partial_{x_a}c(\xi)|\xi|^{-2k}
 -ikc(\xi)|\xi|^{-2k-2}\partial_{x_a}|\xi|^2=0.
\]
The derivative $\partial_{\xi_n}$ of this zero expression is also zero. This proves the vanishing before trace or integration and does not differentiate the inserted fields in the first factor.
The right principal symbol has zero tangential derivative at the chosen point. Hence $\phi_{1,1}=0$, without restricting derivatives of $X$ or $Y$.

\paragraph{Case (a.II).}
Separate $p_{x_n}$ into the derivatives of the inserted fields, the derivative $C_{x_n}={\eta} U/2$, and the denominator derivative $-{\eta} q^{-2}$. The first part gives $-b_k\dot G$. For the metric part, writing $t_0=i(U+zN)q^{-k}$, direct differentiation gives
\[
 t_0''=i\{-4kzNq^{-k-1}-2k(U+zN)q^{-k-1}
             +4k(k+1)\xi_n^2(U+zN)q^{-k-2}\}.
\]
The three products are evaluated using~\eqref{eq:binomialreduced}; their sum is $a_k{\eta} G$. Therefore
\begin{equation}\phi_{1,2}=a_k{\eta} G-b_k\dot G.\end{equation}

\paragraph{Projected integrand and individual residues.}
After the Hardy projection, the Clifford trace and the normalized tangential sphere integral, the complete scalar integrand for this case is
\begin{equation}
\mathfrak I_{1,2}(\xi_n)=\sum_{\rho}M_\rho\frac{P_\rho(\xi_n)}{(\xi_n-i)^{k+3}(\xi_n+i)^{k+2}}.
\end{equation}
The prefactor of the case, including its sign and any factor $1/2$ or $i$, is included here. The following list gives every nonzero geometric monomial and its numerator; terms with the same monomial have already been collected.
\begin{gather*}
M_{1}=g(\dot X,Y),\\
P_{1}(\xi_n)=\frac{k \left(2 k - 1\right)}{2}\xi_n^{3} - \frac{i k \left(2 k + 1\right)}{2}\xi_n^{2} - \frac{3 k}{2}\xi_n + \frac{i k}{2},\\
\frac{1}{\pi\beta_k}\int_{\mathbb R}\frac{P_{1}(\xi_n)\,\dd\xi_n}{(\xi_n-i)^{k+3}(\xi_n+i)^{k+2}}=- \frac{k}{2 \left(k + 1\right)}.
\end{gather*}
\begin{gather*}
M_{2}=\eta g(X,Y),\\
P_{2}(\xi_n)=\frac{i k \left(2 k - 1\right)}{4}\xi_n^{2} + \frac{k}{2}\xi_n - \frac{i k}{4},\\
\frac{1}{\pi\beta_k}\int_{\mathbb R}\frac{P_{2}(\xi_n)\,\dd\xi_n}{(\xi_n-i)^{k+3}(\xi_n+i)^{k+2}}=\frac{3 k \left(2 k + 1\right)}{8 \left(k + 1\right) \left(k + 2\right)}.
\end{gather*}
\begin{gather*}
M_{3}=g(X,\dot Y),\\
P_{3}(\xi_n)=\frac{k \left(2 k - 1\right)}{2}\xi_n^{3} - \frac{i k \left(2 k + 1\right)}{2}\xi_n^{2} - \frac{3 k}{2}\xi_n + \frac{i k}{2},\\
\frac{1}{\pi\beta_k}\int_{\mathbb R}\frac{P_{3}(\xi_n)\,\dd\xi_n}{(\xi_n-i)^{k+3}(\xi_n+i)^{k+2}}=- \frac{k}{2 \left(k + 1\right)}.
\end{gather*}
Multiplying each last line by its $M_\rho$ and summing gives the geometric case coefficient stated above. Thus cancellations inside this case, including any zero individual integrals, are retained explicitly.
\paragraph{Case (a.III).}
Here all normal differentiation falls on the right principal symbol. Insert $t_{x_n}=i{\eta} Uq^{-k}/2-ik{\eta}(U+zN)q^{-k-1}$ and differentiate in $\xi_n$. After the same trace contractions,
\begin{equation}\phi_{1,3}=-c_k{\eta} G.\end{equation}
There are no field derivatives in this case.

\paragraph{Projected integrand and individual residues.}
After the Hardy projection, the Clifford trace and the normalized tangential sphere integral, the complete scalar integrand for this case is
\begin{equation}
\mathfrak I_{1,3}(\xi_n)=\sum_{\rho}M_\rho\frac{P_\rho(\xi_n)}{(\xi_n-i)^{k+3}(\xi_n+i)^{k+2}}.
\end{equation}
The prefactor of the case, including its sign and any factor $1/2$ or $i$, is included here. The following list gives every nonzero geometric monomial and its numerator; terms with the same monomial have already been collected.
\begin{gather*}
M_{1}=\eta g(X,Y),\\
P_{1}(\xi_n)=- \frac{i k}{4}\xi_n^{2} - \frac{k^{2}}{2}\xi_n + \frac{i k}{4},\\
\frac{1}{\pi\beta_k}\int_{\mathbb R}\frac{P_{1}(\xi_n)\,\dd\xi_n}{(\xi_n-i)^{k+3}(\xi_n+i)^{k+2}}=- \frac{k \left(2 k + 1\right)^{2}}{8 \left(k + 1\right) \left(k + 2\right)}.
\end{gather*}
Multiplying each last line by its $M_\rho$ and summing gives the geometric case coefficient stated above. Thus cancellations inside this case, including any zero individual integrals, are retained explicitly.
\paragraph{Case (b).}
\paragraph{A normal contour integral evaluated completely.}
The principal-part projection is
\[
 \pi_{\xi_n}^+\frac1{1+\xi_n^2}=\frac1{2i(\xi_n-i)}.
\]
For $k\ge1$, Cauchy's formula and the binomial theorem imply
\begin{align}
 \int_{\mathbb R}\pi_{\xi_n}^+\frac1{1+\xi_n^2}
             (1+\xi_n^2)^{-k}\dd\xi_n
 &=\pi\beta_k/2,\\
 \int_{\mathbb R}\pi_{\xi_n}^+\frac1{1+\xi_n^2}
       \xi_n^2(1+\xi_n^2)^{-k-1}\dd\xi_n
 &=\frac{\pi\beta_k}{4(k+1)}.
\end{align}
The second equality also follows by subtracting the first integral at $k+1$ from that at $k$ and using $\beta_{k+1}/\beta_k=(2k+1)/(2k+2)$. Their combination with coefficient $-2k$ is $\pi\beta_k/(2(k+1))$, which is exactly the coefficient $t_k$ of the derivative term in case (b). This calculation separates the contour factor, the exterior rank $2^n$, and the sphere volume, preventing a spinor-rank normalization from entering the de Rham result.

The three sources in $p_-$ are the inverse-square correction, the zero-order symbol of $DB$, and the composition correction $-{\eta} ANBq^{-2}$. The field-derivative part is
\[
 -i\int\tau_\Lambda\!\left[
 \pi^+(\hc(X)c_j\hc(Y_{,j})q^{-1})\,
 \partial_{\xi_n}\{i(U+zN)q^{-k}\}\right].
\]
Tangential covariable terms are odd. The normal derivative contributes
\[
 \frac{\mathcal J(1,0;k,0)-2k\mathcal J(1,0;k+1,2)}{\pi\beta_k}
 g(X,\dot Y)=t_kg(X,\dot Y).
\]
The connection terms and the conversion from $Y_{,a}$ to $\nabla_{e_a}Y$ must both be retained. Together with the two remaining sources they give $-d_k{\eta} G$. Thus
\begin{equation}\phi_{1,4}=-d_k{\eta} G+t_kg(X,\dot Y).\end{equation}

\paragraph{Projected integrand and individual residues.}
After the Hardy projection, the Clifford trace and the normalized tangential sphere integral, the complete scalar integrand for this case is
\begin{equation}
\mathfrak I_{1,4}(\xi_n)=\sum_{\rho}M_\rho\frac{P_\rho(\xi_n)}{(\xi_n-i)^{k+3}(\xi_n+i)^{k+1}}.
\end{equation}
The prefactor of the case, including its sign and any factor $1/2$ or $i$, is included here. The following list gives every nonzero geometric monomial and its numerator; terms with the same monomial have already been collected.
\begin{gather*}
M_{1}=\eta g(X,Y),\\
P_{1}(\xi_n)=- \frac{\left(2 k - 1\right) \left(2 k + 1\right)}{8}\xi_n^{2} + \frac{i \left(2 k^{2} + 3 k - 1\right)}{4}\xi_n + \frac{2 k + 3}{8},\\
\frac{1}{\pi\beta_k}\int_{\mathbb R}\frac{P_{1}(\xi_n)\,\dd\xi_n}{(\xi_n-i)^{k+3}(\xi_n+i)^{k+1}}=- \frac{k \left(2 k + 1\right) \left(2 k + 7\right)}{8 \left(k + 1\right) \left(k + 2\right)}.
\end{gather*}
\begin{gather*}
M_{2}=g(X,\dot Y),\\
P_{2}(\xi_n)=\frac{i \left(2 k - 1\right)}{2}\xi_n^{3} + \frac{2 k - 1}{2}\xi_n^{2} - \frac{i}{2}\xi_n - \frac{1}{2},\\
\frac{1}{\pi\beta_k}\int_{\mathbb R}\frac{P_{2}(\xi_n)\,\dd\xi_n}{(\xi_n-i)^{k+3}(\xi_n+i)^{k+1}}=\frac{1}{2 \left(k + 1\right)}.
\end{gather*}
Multiplying each last line by its $M_\rho$ and summing gives the geometric case coefficient stated above. Thus cancellations inside this case, including any zero individual integrals, are retained explicitly.
\paragraph{Case (c).}
Substitute the four summands of $t_-$ in~\eqref{eq:tsub}. The full exterior connection, including $c_a\hc_a\hc_n$, gives
\begin{equation}\phi_{1,5}=e_k{\eta} G.\end{equation}
The cancellation identity following~\eqref{eq:abcde} shows that the total metric contribution is zero. The surviving boundary density depends on the normal first jets of the inserted fields.

\paragraph{Projected integrand and individual residues.}
After the Hardy projection, the Clifford trace and the normalized tangential sphere integral, the complete scalar integrand for this case is
\begin{equation}
\mathfrak I_{1,5}(\xi_n)=\sum_{\rho}M_\rho\frac{P_\rho(\xi_n)}{(\xi_n-i)^{k+3}(\xi_n+i)^{k+3}}.
\end{equation}
The prefactor of the case, including its sign and any factor $1/2$ or $i$, is included here. The following list gives every nonzero geometric monomial and its numerator; terms with the same monomial have already been collected.
\begin{gather*}
M_{1}=\eta g(X,Y),\\
\begin{aligned}P_{1}(\xi_n)&=- \frac{k \left(2 k - 1\right) \left(2 k + 1\right)}{4}\xi_n^{4} + \frac{i k \left(2 k + 1\right)}{4}\xi_n^{3}\\ &\quad {} - \frac{k \left(2 k + 1\right) \left(4 k - 1\right)}{4}\xi_n^{2} + \frac{i k \left(8 k + 7\right)}{4}\xi_n + \frac{k \left(2 k + 1\right)}{2}\end{aligned},\\
\frac{1}{\pi\beta_k}\int_{\mathbb R}\frac{P_{1}(\xi_n)\,\dd\xi_n}{(\xi_n-i)^{k+3}(\xi_n+i)^{k+3}}=\frac{k \left(2 k + 1\right) \left(4 k + 5\right)}{8 \left(k + 1\right) \left(k + 2\right)}.
\end{gather*}
Multiplying each last line by its $M_\rho$ and summing gives the geometric case coefficient stated above. Thus cancellations inside this case, including any zero individual integrals, are retained explicitly.
\subsection{The second insertion: all five cases}
Now $A=\hc(X)\hc(Y)$ and $B=c(Z)c(W)$. The basic trace is
\begin{equation}
 \tau_\Lambda(Ac(v)Bc(w))
 =G\{g(v,Z)g(W,w)-g(v,W)g(Z,w)+g(v,w)H\}.
\end{equation}
When both $v,w$ are the same tangential covariable, the first two terms cancel; the same is true when both are $N$. Cross terms have odd tangential degree. This explains why the principal-symbol cases depend on $GH$.

\paragraph{Case (a.I).} As before, $\phi_{2,1}=0$.

\[
 \Phi_{2,1}=-\int_{|\xi'|=1}\int_{\mathbb R}
 \sum_{a<n}\tr_\Lambda\left[
 \partial_{\xi_a}\pi_{\xi_n}^+\sigma_{-1}(P)
 \partial_{x_a}\partial_{\xi_n}\sigma_{3-n}(Q)\right]
 \dd\xi_n\,\sigma(\xi').
\]
At $x_0$, $\partial_{x_a}g^{ij}=0$ for $a<n$, and the chosen orthonormal trivialization gives $\partial_{x_a}c(\xi)=0$. Therefore
\[
 \partial_{x_a}\sigma_{3-n}(Q)
 =i\,\partial_{x_a}c(\xi)|\xi|^{-2k}
 -ikc(\xi)|\xi|^{-2k-2}\partial_{x_a}|\xi|^2=0.
\]
The derivative $\partial_{\xi_n}$ of this zero expression is also zero. This proves the vanishing before trace or integration and does not differentiate the inserted fields in the first factor.
\paragraph{Case (a.II).} Differentiate every factor of $GH$, not merely $H$. The contour coefficients are those already evaluated above:
\begin{equation}\phi_{2,2}=a_k{\eta} GH-b_k\partial_n(GH).\end{equation}

\paragraph{Projected integrand and individual residues.}
After the Hardy projection, the Clifford trace and the normalized tangential sphere integral, the complete scalar integrand for this case is
\begin{equation}
\mathfrak I_{2,2}(\xi_n)=\sum_{\rho}M_\rho\frac{P_\rho(\xi_n)}{(\xi_n-i)^{k+3}(\xi_n+i)^{k+2}}.
\end{equation}
The prefactor of the case, including its sign and any factor $1/2$ or $i$, is included here. The following list gives every nonzero geometric monomial and its numerator; terms with the same monomial have already been collected.
\begin{gather*}
M_{1}=g(W,Z) g(\dot X,Y),\\
P_{1}(\xi_n)=\frac{k \left(2 k - 1\right)}{2}\xi_n^{3} - \frac{i k \left(2 k + 1\right)}{2}\xi_n^{2} - \frac{3 k}{2}\xi_n + \frac{i k}{2},\\
\frac{1}{\pi\beta_k}\int_{\mathbb R}\frac{P_{1}(\xi_n)\,\dd\xi_n}{(\xi_n-i)^{k+3}(\xi_n+i)^{k+2}}=- \frac{k}{2 \left(k + 1\right)}.
\end{gather*}
\begin{gather*}
M_{2}=g(W,Z) g(X,\dot Y),\\
P_{2}(\xi_n)=\frac{k \left(2 k - 1\right)}{2}\xi_n^{3} - \frac{i k \left(2 k + 1\right)}{2}\xi_n^{2} - \frac{3 k}{2}\xi_n + \frac{i k}{2},\\
\frac{1}{\pi\beta_k}\int_{\mathbb R}\frac{P_{2}(\xi_n)\,\dd\xi_n}{(\xi_n-i)^{k+3}(\xi_n+i)^{k+2}}=- \frac{k}{2 \left(k + 1\right)}.
\end{gather*}
\begin{gather*}
M_{3}=\eta g(W,Z) g(X,Y),\\
P_{3}(\xi_n)=\frac{i k \left(2 k - 1\right)}{4}\xi_n^{2} + \frac{k}{2}\xi_n - \frac{i k}{4},\\
\frac{1}{\pi\beta_k}\int_{\mathbb R}\frac{P_{3}(\xi_n)\,\dd\xi_n}{(\xi_n-i)^{k+3}(\xi_n+i)^{k+2}}=\frac{3 k \left(2 k + 1\right)}{8 \left(k + 1\right) \left(k + 2\right)}.
\end{gather*}
\begin{gather*}
M_{4}=g(W,\dot Z) g(X,Y),\\
P_{4}(\xi_n)=\frac{k \left(2 k - 1\right)}{2}\xi_n^{3} - \frac{i k \left(2 k + 1\right)}{2}\xi_n^{2} - \frac{3 k}{2}\xi_n + \frac{i k}{2},\\
\frac{1}{\pi\beta_k}\int_{\mathbb R}\frac{P_{4}(\xi_n)\,\dd\xi_n}{(\xi_n-i)^{k+3}(\xi_n+i)^{k+2}}=- \frac{k}{2 \left(k + 1\right)}.
\end{gather*}
\begin{gather*}
M_{5}=g(\dot W,Z) g(X,Y),\\
P_{5}(\xi_n)=\frac{k \left(2 k - 1\right)}{2}\xi_n^{3} - \frac{i k \left(2 k + 1\right)}{2}\xi_n^{2} - \frac{3 k}{2}\xi_n + \frac{i k}{2},\\
\frac{1}{\pi\beta_k}\int_{\mathbb R}\frac{P_{5}(\xi_n)\,\dd\xi_n}{(\xi_n-i)^{k+3}(\xi_n+i)^{k+2}}=- \frac{k}{2 \left(k + 1\right)}.
\end{gather*}
Multiplying each last line by its $M_\rho$ and summing gives the geometric case coefficient stated above. Thus cancellations inside this case, including any zero individual integrals, are retained explicitly.
\paragraph{Case (a.III).} The right-symbol derivative gives
\begin{equation}\phi_{2,3}=-c_k{\eta} GH.\end{equation}

\paragraph{Projected integrand and individual residues.}
After the Hardy projection, the Clifford trace and the normalized tangential sphere integral, the complete scalar integrand for this case is
\begin{equation}
\mathfrak I_{2,3}(\xi_n)=\sum_{\rho}M_\rho\frac{P_\rho(\xi_n)}{(\xi_n-i)^{k+3}(\xi_n+i)^{k+2}}.
\end{equation}
The prefactor of the case, including its sign and any factor $1/2$ or $i$, is included here. The following list gives every nonzero geometric monomial and its numerator; terms with the same monomial have already been collected.
\begin{gather*}
M_{1}=\eta g(W,Z) g(X,Y),\\
P_{1}(\xi_n)=- \frac{i k}{4}\xi_n^{2} - \frac{k^{2}}{2}\xi_n + \frac{i k}{4},\\
\frac{1}{\pi\beta_k}\int_{\mathbb R}\frac{P_{1}(\xi_n)\,\dd\xi_n}{(\xi_n-i)^{k+3}(\xi_n+i)^{k+2}}=- \frac{k \left(2 k + 1\right)^{2}}{8 \left(k + 1\right) \left(k + 2\right)}.
\end{gather*}
Multiplying each last line by its $M_\rho$ and summing gives the geometric case coefficient stated above. Thus cancellations inside this case, including any zero individual integrals, are retained explicitly.
\paragraph{Case (b).} The Leibniz expansion
\[
 B_{,j}=c(Z_{,j})c(W)+c(Z)c(W_{,j})
\]
produces both normal and tangential derivatives. Its normal contraction is $G\dot H$. For the tangential contraction, the four-Clifford trace gives the four terms
\[
 G\{g(\nabla_{Z^T}W,N)-g(\nabla_{W^T}Z,N)
       -Z_n\mathcal D(W)+W_n\mathcal D(Z)\}=GT_-.
\]
Using~\eqref{eq:Vjets} for the conversion of ordinary coefficients and retaining the connection part of $w$ yields
\begin{equation}\phi_{2,4}=-d_k{\eta} GH+t_kG(\dot H+T_-).\end{equation}

\paragraph{Projected integrand and individual residues.}
After the Hardy projection, the Clifford trace and the normalized tangential sphere integral, the complete scalar integrand for this case is
\begin{equation}
\mathfrak I_{2,4}(\xi_n)=\sum_{\rho}M_\rho\frac{P_\rho(\xi_n)}{(\xi_n-i)^{k+3}(\xi_n+i)^{k+1}}.
\end{equation}
The prefactor of the case, including its sign and any factor $1/2$ or $i$, is included here. The following list gives every nonzero geometric monomial and its numerator; terms with the same monomial have already been collected.
\begin{gather*}
M_{1}=\eta g(W,Z) g(X,Y),\\
P_{1}(\xi_n)=- \frac{\left(2 k - 1\right) \left(2 k + 1\right)}{8}\xi_n^{2} + \frac{i \left(2 k^{2} + 3 k - 1\right)}{4}\xi_n + \frac{2 k + 3}{8},\\
\frac{1}{\pi\beta_k}\int_{\mathbb R}\frac{P_{1}(\xi_n)\,\dd\xi_n}{(\xi_n-i)^{k+3}(\xi_n+i)^{k+1}}=- \frac{k \left(2 k + 1\right) \left(2 k + 7\right)}{8 \left(k + 1\right) \left(k + 2\right)}.
\end{gather*}
\begin{gather*}
M_{2}=\eta g(N,W) g(N,X) g(Y,Z),\\
P_{2}(\xi_n)=- \frac{i k \left(2 k - 1\right)}{4 \left(2 k + 1\right)}\xi_n^{3} - \frac{\left(2 k - 1\right)^{2}}{8 \left(2 k + 1\right)}\xi_n^{2} + \frac{2 k - 1}{8 \left(2 k + 1\right)},\\
\frac{1}{\pi\beta_k}\int_{\mathbb R}\frac{P_{2}(\xi_n)\,\dd\xi_n}{(\xi_n-i)^{k+3}(\xi_n+i)^{k+1}}=0.
\end{gather*}
\begin{gather*}
M_{3}=\eta g(N,Y) g(N,Z) g(W,X),\\
P_{3}(\xi_n)=- \frac{i k \left(2 k - 1\right)}{4 \left(2 k + 1\right)}\xi_n^{3} - \frac{\left(2 k - 1\right)^{2}}{8 \left(2 k + 1\right)}\xi_n^{2} + \frac{2 k - 1}{8 \left(2 k + 1\right)},\\
\frac{1}{\pi\beta_k}\int_{\mathbb R}\frac{P_{3}(\xi_n)\,\dd\xi_n}{(\xi_n-i)^{k+3}(\xi_n+i)^{k+1}}=0.
\end{gather*}
\begin{gather*}
M_{4}=\eta g(N,W) g(N,Y) g(X,Z),\\
P_{4}(\xi_n)=\frac{i k \left(2 k - 1\right)}{4 \left(2 k + 1\right)}\xi_n^{3} + \frac{\left(2 k - 1\right)^{2}}{8 \left(2 k + 1\right)}\xi_n^{2} - \frac{2 k - 1}{8 \left(2 k + 1\right)},\\
\frac{1}{\pi\beta_k}\int_{\mathbb R}\frac{P_{4}(\xi_n)\,\dd\xi_n}{(\xi_n-i)^{k+3}(\xi_n+i)^{k+1}}=0.
\end{gather*}
\begin{gather*}
M_{5}=\eta g(N,X) g(N,Z) g(W,Y),\\
P_{5}(\xi_n)=\frac{i k \left(2 k - 1\right)}{4 \left(2 k + 1\right)}\xi_n^{3} + \frac{\left(2 k - 1\right)^{2}}{8 \left(2 k + 1\right)}\xi_n^{2} - \frac{2 k - 1}{8 \left(2 k + 1\right)},\\
\frac{1}{\pi\beta_k}\int_{\mathbb R}\frac{P_{5}(\xi_n)\,\dd\xi_n}{(\xi_n-i)^{k+3}(\xi_n+i)^{k+1}}=0.
\end{gather*}
\begin{gather*}
M_{6}=g(W,\dot Z) g(X,Y),\\
P_{6}(\xi_n)=\frac{i \left(2 k - 1\right)}{2}\xi_n^{3} + \frac{2 k - 1}{2}\xi_n^{2} - \frac{i}{2}\xi_n - \frac{1}{2},\\
\frac{1}{\pi\beta_k}\int_{\mathbb R}\frac{P_{6}(\xi_n)\,\dd\xi_n}{(\xi_n-i)^{k+3}(\xi_n+i)^{k+1}}=\frac{1}{2 \left(k + 1\right)}.
\end{gather*}
\begin{gather*}
M_{7}=g(\dot W,Z) g(X,Y),\\
P_{7}(\xi_n)=\frac{i \left(2 k - 1\right)}{2}\xi_n^{3} + \frac{2 k - 1}{2}\xi_n^{2} - \frac{i}{2}\xi_n - \frac{1}{2},\\
\frac{1}{\pi\beta_k}\int_{\mathbb R}\frac{P_{7}(\xi_n)\,\dd\xi_n}{(\xi_n-i)^{k+3}(\xi_n+i)^{k+1}}=\frac{1}{2 \left(k + 1\right)}.
\end{gather*}
\begin{gather*}
M_{8}=g(\nabla_{W^T}Z,N) g(X,Y),\\
P_{8}(\xi_n)=- \frac{i \left(2 k - 1\right)}{2}\xi_n^{3} - \frac{2 k - 1}{2}\xi_n^{2} + \frac{i}{2}\xi_n + \frac{1}{2},\\
\frac{1}{\pi\beta_k}\int_{\mathbb R}\frac{P_{8}(\xi_n)\,\dd\xi_n}{(\xi_n-i)^{k+3}(\xi_n+i)^{k+1}}=- \frac{1}{2 \left(k + 1\right)}.
\end{gather*}
\begin{gather*}
M_{9}=\mathcal D(Z) g(N,W) g(X,Y),\\
P_{9}(\xi_n)=\frac{i \left(2 k - 1\right)}{2}\xi_n^{3} + \frac{2 k - 1}{2}\xi_n^{2} - \frac{i}{2}\xi_n - \frac{1}{2},\\
\frac{1}{\pi\beta_k}\int_{\mathbb R}\frac{P_{9}(\xi_n)\,\dd\xi_n}{(\xi_n-i)^{k+3}(\xi_n+i)^{k+1}}=\frac{1}{2 \left(k + 1\right)}.
\end{gather*}
\begin{gather*}
M_{10}=g(\nabla_{Z^T}W,N) g(X,Y),\\
P_{10}(\xi_n)=\frac{i \left(2 k - 1\right)}{2}\xi_n^{3} + \frac{2 k - 1}{2}\xi_n^{2} - \frac{i}{2}\xi_n - \frac{1}{2},\\
\frac{1}{\pi\beta_k}\int_{\mathbb R}\frac{P_{10}(\xi_n)\,\dd\xi_n}{(\xi_n-i)^{k+3}(\xi_n+i)^{k+1}}=\frac{1}{2 \left(k + 1\right)}.
\end{gather*}
\begin{gather*}
M_{11}=\mathcal D(W) g(N,Z) g(X,Y),\\
P_{11}(\xi_n)=- \frac{i \left(2 k - 1\right)}{2}\xi_n^{3} - \frac{2 k - 1}{2}\xi_n^{2} + \frac{i}{2}\xi_n + \frac{1}{2},\\
\frac{1}{\pi\beta_k}\int_{\mathbb R}\frac{P_{11}(\xi_n)\,\dd\xi_n}{(\xi_n-i)^{k+3}(\xi_n+i)^{k+1}}=- \frac{1}{2 \left(k + 1\right)}.
\end{gather*}
Multiplying each last line by its $M_\rho$ and summing gives the geometric case coefficient stated above. Thus cancellations inside this case, including any zero individual integrals, are retained explicitly.
\paragraph{Case (c).} The subprincipal right factor contributes
\begin{equation}\phi_{2,5}=e_k{\eta} GH.\end{equation}
Thus the explicit mean-curvature terms cancel in the sum, but the unrestricted derivative terms remain. In particular this is not a claim that the full boundary density is zero.

\paragraph{Projected integrand and individual residues.}
After the Hardy projection, the Clifford trace and the normalized tangential sphere integral, the complete scalar integrand for this case is
\begin{equation}
\mathfrak I_{2,5}(\xi_n)=\sum_{\rho}M_\rho\frac{P_\rho(\xi_n)}{(\xi_n-i)^{k+3}(\xi_n+i)^{k+3}}.
\end{equation}
The prefactor of the case, including its sign and any factor $1/2$ or $i$, is included here. The following list gives every nonzero geometric monomial and its numerator; terms with the same monomial have already been collected.
\begin{gather*}
M_{1}=\eta g(W,Z) g(X,Y),\\
\begin{aligned}P_{1}(\xi_n)&=- \frac{k \left(2 k - 1\right) \left(2 k + 1\right)}{4}\xi_n^{4} + \frac{i k \left(2 k + 1\right)}{4}\xi_n^{3}\\ &\quad {} - \frac{k \left(2 k + 1\right) \left(4 k - 1\right)}{4}\xi_n^{2} + \frac{i k \left(8 k + 7\right)}{4}\xi_n + \frac{k \left(2 k + 1\right)}{2}\end{aligned},\\
\frac{1}{\pi\beta_k}\int_{\mathbb R}\frac{P_{1}(\xi_n)\,\dd\xi_n}{(\xi_n-i)^{k+3}(\xi_n+i)^{k+3}}=\frac{k \left(2 k + 1\right) \left(4 k + 5\right)}{8 \left(k + 1\right) \left(k + 2\right)}.
\end{gather*}
\begin{gather*}
M_{2}=\eta g(N,W) g(N,Y) g(X,Z),\\
P_{2}(\xi_n)=\frac{k}{4}\xi_n^{4} + \frac{k^{2}}{2 \left(2 k + 1\right)}\xi_n^{2} - \frac{k}{4 \left(2 k + 1\right)},\\
\frac{1}{\pi\beta_k}\int_{\mathbb R}\frac{P_{2}(\xi_n)\,\dd\xi_n}{(\xi_n-i)^{k+3}(\xi_n+i)^{k+3}}=0.
\end{gather*}
\begin{gather*}
M_{3}=\eta g(N,X) g(N,Z) g(W,Y),\\
P_{3}(\xi_n)=\frac{k}{4}\xi_n^{4} + \frac{k^{2}}{2 \left(2 k + 1\right)}\xi_n^{2} - \frac{k}{4 \left(2 k + 1\right)},\\
\frac{1}{\pi\beta_k}\int_{\mathbb R}\frac{P_{3}(\xi_n)\,\dd\xi_n}{(\xi_n-i)^{k+3}(\xi_n+i)^{k+3}}=0.
\end{gather*}
\begin{gather*}
M_{4}=\eta g(N,W) g(N,X) g(Y,Z),\\
P_{4}(\xi_n)=- \frac{k}{4}\xi_n^{4} - \frac{k^{2}}{2 \left(2 k + 1\right)}\xi_n^{2} + \frac{k}{4 \left(2 k + 1\right)},\\
\frac{1}{\pi\beta_k}\int_{\mathbb R}\frac{P_{4}(\xi_n)\,\dd\xi_n}{(\xi_n-i)^{k+3}(\xi_n+i)^{k+3}}=0.
\end{gather*}
\begin{gather*}
M_{5}=\eta g(N,Y) g(N,Z) g(W,X),\\
P_{5}(\xi_n)=- \frac{k}{4}\xi_n^{4} - \frac{k^{2}}{2 \left(2 k + 1\right)}\xi_n^{2} + \frac{k}{4 \left(2 k + 1\right)},\\
\frac{1}{\pi\beta_k}\int_{\mathbb R}\frac{P_{5}(\xi_n)\,\dd\xi_n}{(\xi_n-i)^{k+3}(\xi_n+i)^{k+3}}=0.
\end{gather*}
Multiplying each last line by its $M_\rho$ and summing gives the geometric case coefficient stated above. Thus cancellations inside this case, including any zero individual integrals, are retained explicitly.
\subsection{The third insertion: all five cases}
For $A=\hc(X)\hc(Y)c(Z)$ and $B=c(W)$, the principal trace is
\begin{equation}
 \tau_\Lambda(Ac(v)Bc(w))
 =G\{g(Z,v)g(W,w)-Hg(v,w)+g(Z,w)g(v,W)\}.
\end{equation}
Tangential averaging and normal evaluation give, respectively,
\begin{equation}\label{eq:thirdtraces}
 G\left(\frac{2(H-P_1)}d-H\right),\qquad G(2P_1-H).
\end{equation}
These are different and must not be combined before the contour coefficients are evaluated. Set $\Delta_k=8(k+1)(k+2)(k+3)$ and introduce
\begin{align}
 A_H&=-\frac{k(6k^2+15k+3)}{\Delta_k},&
 A_{P_1}&=\frac{k(18k+6)}{\Delta_k},\\
 B_H&=\frac{k(4k^3+12k^2+5k+3)}{\Delta_k},&
 B_{P_1}&=-\frac{k(8k^2+10k+6)}{\Delta_k},\\
 C_H&=\frac{4k^4+24k^3+31k^2-3k}{\Delta_k},&
 C_{P_1}&=\frac{-14k^2+6k}{\Delta_k},\\
 E_H&=-\frac{k(8k^3+30k^2+25k+9)}{\Delta_k},&
 E_{P_1}&=\frac{k(16k^2+38k+18)}{\Delta_k},\\
 \lambda_k&=\frac{2k(k+3)}{\Delta_k},&
 u_k&=\frac{k}{2(k+2)},\qquad v_k=\frac{k}{(k+1)(k+2)}.
\end{align}

\paragraph{The third-insertion trace, term by term.}
From $c(V)c(W)+c(W)c(V)=-2g(V,W)$ and cyclicity,
\[
 \tau_\Lambda(c(Z)c(v)c(W)c(w))
 =g(Z,v)g(W,w)-g(Z,W)g(v,w)+g(Z,w)g(v,W).
\]
The hatted pair contributes $G=g(X,Y)$, because it commutes with the even unhatted word and its nonscalar component has zero trace against that word. With $v=w=\xi'$, the first and third terms each give $(H-P_1)/(n-1)$ after angular averaging; the middle term gives $-H$. With $v=w=N$, the three terms give $P_1,-H,P_1$. This derives the two different contractions in~\eqref{eq:thirdtraces} and explains why the coefficient of $Z_nW_n$ cannot be inferred from the coefficient of $g(Z^T,W^T)$.

\paragraph{Case (a.I).} Again $\phi_{3,1}=0$.

\[
 \Phi_{3,1}=-\int_{|\xi'|=1}\int_{\mathbb R}
 \sum_{a<n}\tr_\Lambda\left[
 \partial_{\xi_a}\pi_{\xi_n}^+\sigma_{-1}(P)
 \partial_{x_a}\partial_{\xi_n}\sigma_{3-n}(Q)\right]
 \dd\xi_n\,\sigma(\xi').
\]
At $x_0$, $\partial_{x_a}g^{ij}=0$ for $a<n$, and the chosen orthonormal trivialization gives $\partial_{x_a}c(\xi)=0$. Therefore
\[
 \partial_{x_a}\sigma_{3-n}(Q)
 =i\,\partial_{x_a}c(\xi)|\xi|^{-2k}
 -ikc(\xi)|\xi|^{-2k-2}\partial_{x_a}|\xi|^2=0.
\]
The derivative $\partial_{\xi_n}$ of this zero expression is also zero. This proves the vanishing before trace or integration and does not differentiate the inserted fields in the first factor.

\paragraph{Case (a.II).} The two contractions in~\eqref{eq:thirdtraces}, together with $C_{x_n}$ and $q_{x_n}$, give
\begin{equation}
 \phi_{3,2}={\eta} G(A_HH+A_{P_1}P)+u_k\partial_n(GH)-v_k\partial_n(GP).
\end{equation}
The derivative terms include derivatives of $X,Y,Z,W$. For example, the coefficient of $\dot G$ is $u_kH-v_kP_1$; the coefficients of $G\dot H$ and $G\dot P_1$ are $u_k$ and $-v_k$.

\paragraph{Projected integrand and individual residues.}
After the Hardy projection, the Clifford trace and the normalized tangential sphere integral, the complete scalar integrand for this case is
\begin{equation}
\mathfrak I_{3,2}(\xi_n)=\sum_{\rho}M_\rho\frac{P_\rho(\xi_n)}{(\xi_n-i)^{k+4}(\xi_n+i)^{k+2}}.
\end{equation}
The prefactor of the case, including its sign and any factor $1/2$ or $i$, is included here. The following list gives every nonzero geometric monomial and its numerator; terms with the same monomial have already been collected.
\begin{gather*}
M_{1}=g(W,Z) g(\dot X,Y),\\
P_{1}(\xi_n)=- \frac{k \left(2 k - 1\right)}{2}\xi_n^{4} + i k \left(2 k - 1\right)\xi_n^{3} + k \left(k + 1\right)\xi_n^{2} - \frac{i k \left(4 k + 1\right)}{2 k + 1}\xi_n - \frac{k \left(2 k - 1\right)}{2 \left(2 k + 1\right)},\\
\frac{1}{\pi\beta_k}\int_{\mathbb R}\frac{P_{1}(\xi_n)\,\dd\xi_n}{(\xi_n-i)^{k+4}(\xi_n+i)^{k+2}}=\frac{k}{2 \left(k + 2\right)}.
\end{gather*}
\begin{gather*}
M_{2}=g(N,W) g(N,Z) g(\dot X,Y),\\
P_{2}(\xi_n)=k \left(2 k - 1\right)\xi_n^{4} - 2 i k \left(k - 1\right)\xi_n^{3} - 2 k\xi_n^{2} + \frac{2 i k \left(3 k + 1\right)}{2 k + 1}\xi_n - \frac{k}{2 k + 1},\\
\frac{1}{\pi\beta_k}\int_{\mathbb R}\frac{P_{2}(\xi_n)\,\dd\xi_n}{(\xi_n-i)^{k+4}(\xi_n+i)^{k+2}}=- \frac{k}{\left(k + 1\right) \left(k + 2\right)}.
\end{gather*}
\begin{gather*}
M_{3}=g(W,Z) g(X,\dot Y),\\
P_{3}(\xi_n)=- \frac{k \left(2 k - 1\right)}{2}\xi_n^{4} + i k \left(2 k - 1\right)\xi_n^{3} + k \left(k + 1\right)\xi_n^{2} - \frac{i k \left(4 k + 1\right)}{2 k + 1}\xi_n - \frac{k \left(2 k - 1\right)}{2 \left(2 k + 1\right)},\\
\frac{1}{\pi\beta_k}\int_{\mathbb R}\frac{P_{3}(\xi_n)\,\dd\xi_n}{(\xi_n-i)^{k+4}(\xi_n+i)^{k+2}}=\frac{k}{2 \left(k + 2\right)}.
\end{gather*}
\begin{gather*}
M_{4}=g(N,W) g(N,Z) g(X,\dot Y),\\
P_{4}(\xi_n)=k \left(2 k - 1\right)\xi_n^{4} - 2 i k \left(k - 1\right)\xi_n^{3} - 2 k\xi_n^{2} + \frac{2 i k \left(3 k + 1\right)}{2 k + 1}\xi_n - \frac{k}{2 k + 1},\\
\frac{1}{\pi\beta_k}\int_{\mathbb R}\frac{P_{4}(\xi_n)\,\dd\xi_n}{(\xi_n-i)^{k+4}(\xi_n+i)^{k+2}}=- \frac{k}{\left(k + 1\right) \left(k + 2\right)}.
\end{gather*}
\begin{gather*}
M_{5}=g(W,\dot Z) g(X,Y),\\
P_{5}(\xi_n)=- \frac{k \left(2 k - 1\right)}{2}\xi_n^{4} + i k \left(2 k - 1\right)\xi_n^{3} + k \left(k + 1\right)\xi_n^{2} - \frac{i k \left(4 k + 1\right)}{2 k + 1}\xi_n - \frac{k \left(2 k - 1\right)}{2 \left(2 k + 1\right)},\\
\frac{1}{\pi\beta_k}\int_{\mathbb R}\frac{P_{5}(\xi_n)\,\dd\xi_n}{(\xi_n-i)^{k+4}(\xi_n+i)^{k+2}}=\frac{k}{2 \left(k + 2\right)}.
\end{gather*}
\begin{gather*}
M_{6}=g(N,W) g(N,\dot Z) g(X,Y),\\
P_{6}(\xi_n)=k \left(2 k - 1\right)\xi_n^{4} - 2 i k \left(k - 1\right)\xi_n^{3} - 2 k\xi_n^{2} + \frac{2 i k \left(3 k + 1\right)}{2 k + 1}\xi_n - \frac{k}{2 k + 1},\\
\frac{1}{\pi\beta_k}\int_{\mathbb R}\frac{P_{6}(\xi_n)\,\dd\xi_n}{(\xi_n-i)^{k+4}(\xi_n+i)^{k+2}}=- \frac{k}{\left(k + 1\right) \left(k + 2\right)}.
\end{gather*}
\begin{gather*}
M_{7}=\eta g(W,Z) g(X,Y),\\
P_{7}(\xi_n)=- \frac{i k \left(2 k - 1\right)}{4}\xi_n^{3} - \frac{k \left(2 k - 1\right)}{4}\xi_n^{2} + \frac{3 i k}{4}\xi_n + \frac{k \left(2 k - 1\right)}{4 \left(2 k + 1\right)},\\
\frac{1}{\pi\beta_k}\int_{\mathbb R}\frac{P_{7}(\xi_n)\,\dd\xi_n}{(\xi_n-i)^{k+4}(\xi_n+i)^{k+2}}=- \frac{3 k \left(2 k^{2} + 5 k + 1\right)}{8 \left(k + 1\right) \left(k + 2\right) \left(k + 3\right)}.
\end{gather*}
\begin{gather*}
M_{8}=\eta g(N,W) g(N,Z) g(X,Y),\\
P_{8}(\xi_n)=\frac{i k \left(2 k - 1\right)}{2}\xi_n^{3} - \frac{k}{2}\xi_n^{2} - \frac{3 i k}{2}\xi_n + \frac{k}{2 \left(2 k + 1\right)},\\
\frac{1}{\pi\beta_k}\int_{\mathbb R}\frac{P_{8}(\xi_n)\,\dd\xi_n}{(\xi_n-i)^{k+4}(\xi_n+i)^{k+2}}=\frac{3 k \left(3 k + 1\right)}{4 \left(k + 1\right) \left(k + 2\right) \left(k + 3\right)}.
\end{gather*}
\begin{gather*}
M_{9}=g(\dot W,Z) g(X,Y),\\
P_{9}(\xi_n)=- \frac{k \left(2 k - 1\right)}{2}\xi_n^{4} + i k \left(2 k - 1\right)\xi_n^{3} + k \left(k + 1\right)\xi_n^{2} - \frac{i k \left(4 k + 1\right)}{2 k + 1}\xi_n - \frac{k \left(2 k - 1\right)}{2 \left(2 k + 1\right)},\\
\frac{1}{\pi\beta_k}\int_{\mathbb R}\frac{P_{9}(\xi_n)\,\dd\xi_n}{(\xi_n-i)^{k+4}(\xi_n+i)^{k+2}}=\frac{k}{2 \left(k + 2\right)}.
\end{gather*}
\begin{gather*}
M_{10}=g(N,\dot W) g(N,Z) g(X,Y),\\
P_{10}(\xi_n)=k \left(2 k - 1\right)\xi_n^{4} - 2 i k \left(k - 1\right)\xi_n^{3} - 2 k\xi_n^{2} + \frac{2 i k \left(3 k + 1\right)}{2 k + 1}\xi_n - \frac{k}{2 k + 1},\\
\frac{1}{\pi\beta_k}\int_{\mathbb R}\frac{P_{10}(\xi_n)\,\dd\xi_n}{(\xi_n-i)^{k+4}(\xi_n+i)^{k+2}}=- \frac{k}{\left(k + 1\right) \left(k + 2\right)}.
\end{gather*}
Multiplying each last line by its $M_\rho$ and summing gives the geometric case coefficient stated above. Thus cancellations inside this case, including any zero individual integrals, are retained explicitly.
\paragraph{Case (a.III).} Keeping the tangential and normal contractions separate until after~\eqref{eq:binomialreduced} gives
\begin{equation}\phi_{3,3}={\eta} G(B_HH+B_{P_1}P).\end{equation}

\paragraph{Projected integrand and individual residues.}
After the Hardy projection, the Clifford trace and the normalized tangential sphere integral, the complete scalar integrand for this case is
\begin{equation}
\mathfrak I_{3,3}(\xi_n)=\sum_{\rho}M_\rho\frac{P_\rho(\xi_n)}{(\xi_n-i)^{k+4}(\xi_n+i)^{k+2}}.
\end{equation}
The prefactor of the case, including its sign and any factor $1/2$ or $i$, is included here. The following list gives every nonzero geometric monomial and its numerator; terms with the same monomial have already been collected.
\begin{gather*}
M_{1}=\eta g(W,Z) g(X,Y),\\
P_{1}(\xi_n)=\frac{i k \left(2 k - 1\right)}{4 \left(2 k + 1\right)}\xi_n^{3} + \frac{k \left(2 k + 1\right)}{4}\xi_n^{2} - \frac{i k \left(2 k - 1\right)}{4}\xi_n - \frac{k}{4},\\
\frac{1}{\pi\beta_k}\int_{\mathbb R}\frac{P_{1}(\xi_n)\,\dd\xi_n}{(\xi_n-i)^{k+4}(\xi_n+i)^{k+2}}=\frac{k \left(4 k^{3} + 12 k^{2} + 5 k + 3\right)}{8 \left(k + 1\right) \left(k + 2\right) \left(k + 3\right)}.
\end{gather*}
\begin{gather*}
M_{2}=\eta g(N,W) g(N,Z) g(X,Y),\\
P_{2}(\xi_n)=\frac{i k}{2 \left(2 k + 1\right)}\xi_n^{3} - \frac{k \left(2 k + 1\right)}{2}\xi_n^{2} - \frac{i k}{2}\xi_n + \frac{k}{2},\\
\frac{1}{\pi\beta_k}\int_{\mathbb R}\frac{P_{2}(\xi_n)\,\dd\xi_n}{(\xi_n-i)^{k+4}(\xi_n+i)^{k+2}}=- \frac{k \left(4 k^{2} + 5 k + 3\right)}{4 \left(k + 1\right) \left(k + 2\right) \left(k + 3\right)}.
\end{gather*}
Multiplying each last line by its $M_\rho$ and summing gives the geometric case coefficient stated above. Thus cancellations inside this case, including any zero individual integrals, are retained explicitly.
\paragraph{Case (b).} The derivative contribution has a particularly useful geometric form:
\paragraph{The tangential covariant-derivative contribution.}
For the third insertion, take $v=e_a$, $W$ replaced by $\nabla_{e_a}W$, and $w=N$ in the same trace. Summing $a<n$ gives
\begin{align}
 &\sum_{a<n}\tau_\Lambda\big(\hc(X)\hc(Y)c(Z)c(e_a)
              c(\nabla_{e_a}W)c(N)\big)\notag\\
 &\quad=G\sum_{a<n}\{Z_a g(\nabla_{e_a}W,N)
              +Z_ng(e_a,\nabla_{e_a}W)\}\notag\\
 &\quad=G\{g(\nabla_{Z^T}W,N)+Z_n\mathcal D(W)\}=GT_+.
\end{align}
The omitted middle contraction is zero because $g(e_a,N)=0$. The normal summand instead equals
\[
 G\{2Z_ng(\nabla_NW,N)-g(Z,\nabla_NW)\}.
\]
The second insertion has two Leibniz terms. Performing this calculation on $c(\nabla Z)c(W)$ and $c(Z)c(\nabla W)$, with their different orders, gives respectively the negative and positive directional terms in $T_-$. This establishes the relative signs of the two divergences as well.

\begin{equation}
 t_kG\{T_++2Z_ng(\dot W,N)-g(Z,\dot W)\}.
\end{equation}
The terms $g(\nabla_{Z^T}W,N)$ and $Z_n\mathcal D(W)$ arise from the two nonzero tangential contractions; the last two terms arise from $j=n$. The part of $w$ containing two hats also contributes a first-pair skew expression, namely $-{\eta}\lambda_kL$. The complete case is
\begin{equation}
 \phi_{3,4}={\eta} G(C_HH+C_{P_1}P)-{\eta}\lambda_kL
 +t_kG\{T_++2Z_ng(\dot W,N)-g(Z,\dot W)\}.
\end{equation}

\paragraph{Projected integrand and individual residues.}
After the Hardy projection, the Clifford trace and the normalized tangential sphere integral, the complete scalar integrand for this case is
\begin{equation}
\mathfrak I_{3,4}(\xi_n)=\sum_{\rho}M_\rho\frac{P_\rho(\xi_n)}{(\xi_n-i)^{k+4}(\xi_n+i)^{k+1}}.
\end{equation}
The prefactor of the case, including its sign and any factor $1/2$ or $i$, is included here. The following list gives every nonzero geometric monomial and its numerator; terms with the same monomial have already been collected.
\begin{gather*}
M_{1}=\eta g(W,Z) g(X,Y),\\
\begin{aligned}P_{1}(\xi_n)&=- \frac{i k \left(2 k - 1\right)}{2 \left(2 k + 1\right)}\xi_n^{4} + \frac{\left(2 k - 1\right) \left(4 k^{2} - 4 k + 3\right)}{8 \left(2 k + 1\right)}\xi_n^{3} - \frac{i \left(16 k^{3} + 4 k^{2} - 4 k - 5\right)}{8 \left(2 k + 1\right)}\xi_n^{2}\\ &\quad {} - \frac{\left(2 k - 1\right) \left(4 k^{2} + 8 k - 3\right)}{8 \left(2 k + 1\right)}\xi_n + \frac{i \left(4 k^{2} + 4 k + 5\right)}{8 \left(2 k + 1\right)}\end{aligned},\\
\frac{1}{\pi\beta_k}\int_{\mathbb R}\frac{P_{1}(\xi_n)\,\dd\xi_n}{(\xi_n-i)^{k+4}(\xi_n+i)^{k+1}}=\frac{k \left(4 k^{3} + 24 k^{2} + 31 k - 3\right)}{8 \left(k + 1\right) \left(k + 2\right) \left(k + 3\right)}.
\end{gather*}
\begin{gather*}
M_{2}=\eta g(N,W) g(N,Z) g(X,Y),\\
\begin{aligned}P_{2}(\xi_n)&=\frac{i k \left(k + 1\right) \left(2 k - 1\right)}{2 k + 1}\xi_n^{4} + \frac{\left(2 k - 1\right) \left(4 k^{2} + 2 k - 3\right)}{4 \left(2 k + 1\right)}\xi_n^{3}\\ &\quad {} - \frac{i \left(6 k + 5\right)}{4 \left(2 k + 1\right)}\xi_n^{2} - \frac{3 \left(4 k^{2} + 4 k - 1\right)}{4 \left(2 k + 1\right)}\xi_n - \frac{i \left(6 k + 5\right)}{4 \left(2 k + 1\right)}\end{aligned},\\
\frac{1}{\pi\beta_k}\int_{\mathbb R}\frac{P_{2}(\xi_n)\,\dd\xi_n}{(\xi_n-i)^{k+4}(\xi_n+i)^{k+1}}=- \frac{k \left(7 k - 3\right)}{4 \left(k + 1\right) \left(k + 2\right) \left(k + 3\right)}.
\end{gather*}
\begin{gather*}
M_{3}=\eta g(N,W) g(N,Y) g(X,Z),\\
P_{3}(\xi_n)=\frac{i k \left(2 k - 1\right)}{4 \left(2 k + 1\right)}\xi_n^{4} + \frac{\left(2 k - 1\right) \left(4 k - 1\right)}{8 \left(2 k + 1\right)}\xi_n^{3} - \frac{i \left(4 k^{2} + 1\right)}{8 \left(2 k + 1\right)}\xi_n^{2} - \frac{6 k - 1}{8 \left(2 k + 1\right)}\xi_n + \frac{i \left(2 k - 1\right)}{8 \left(2 k + 1\right)},\\
\frac{1}{\pi\beta_k}\int_{\mathbb R}\frac{P_{3}(\xi_n)\,\dd\xi_n}{(\xi_n-i)^{k+4}(\xi_n+i)^{k+1}}=\frac{k}{4 \left(k + 1\right) \left(k + 2\right)}.
\end{gather*}
\begin{gather*}
M_{4}=\eta g(N,X) g(N,Z) g(W,Y),\\
P_{4}(\xi_n)=- \frac{i k \left(2 k - 1\right)}{4 \left(2 k + 1\right)}\xi_n^{4} - \frac{\left(2 k - 1\right) \left(4 k - 1\right)}{8 \left(2 k + 1\right)}\xi_n^{3} + \frac{i \left(4 k^{2} + 1\right)}{8 \left(2 k + 1\right)}\xi_n^{2} + \frac{6 k - 1}{8 \left(2 k + 1\right)}\xi_n - \frac{i \left(2 k - 1\right)}{8 \left(2 k + 1\right)},\\
\frac{1}{\pi\beta_k}\int_{\mathbb R}\frac{P_{4}(\xi_n)\,\dd\xi_n}{(\xi_n-i)^{k+4}(\xi_n+i)^{k+1}}=- \frac{k}{4 \left(k + 1\right) \left(k + 2\right)}.
\end{gather*}
\begin{gather*}
M_{5}=\eta g(N,W) g(N,X) g(Y,Z),\\
P_{5}(\xi_n)=- \frac{i k \left(2 k - 1\right)}{4 \left(2 k + 1\right)}\xi_n^{4} - \frac{\left(2 k - 1\right) \left(4 k - 1\right)}{8 \left(2 k + 1\right)}\xi_n^{3} + \frac{i \left(4 k^{2} + 1\right)}{8 \left(2 k + 1\right)}\xi_n^{2} + \frac{6 k - 1}{8 \left(2 k + 1\right)}\xi_n - \frac{i \left(2 k - 1\right)}{8 \left(2 k + 1\right)},\\
\frac{1}{\pi\beta_k}\int_{\mathbb R}\frac{P_{5}(\xi_n)\,\dd\xi_n}{(\xi_n-i)^{k+4}(\xi_n+i)^{k+1}}=- \frac{k}{4 \left(k + 1\right) \left(k + 2\right)}.
\end{gather*}
\begin{gather*}
M_{6}=\eta g(N,Y) g(N,Z) g(W,X),\\
P_{6}(\xi_n)=\frac{i k \left(2 k - 1\right)}{4 \left(2 k + 1\right)}\xi_n^{4} + \frac{\left(2 k - 1\right) \left(4 k - 1\right)}{8 \left(2 k + 1\right)}\xi_n^{3} - \frac{i \left(4 k^{2} + 1\right)}{8 \left(2 k + 1\right)}\xi_n^{2} - \frac{6 k - 1}{8 \left(2 k + 1\right)}\xi_n + \frac{i \left(2 k - 1\right)}{8 \left(2 k + 1\right)},\\
\frac{1}{\pi\beta_k}\int_{\mathbb R}\frac{P_{6}(\xi_n)\,\dd\xi_n}{(\xi_n-i)^{k+4}(\xi_n+i)^{k+1}}=\frac{k}{4 \left(k + 1\right) \left(k + 2\right)}.
\end{gather*}
\begin{gather*}
M_{7}=g(\dot W,Z) g(X,Y),\\
P_{7}(\xi_n)=- \frac{i \left(2 k - 1\right)}{2}\xi_n^{4} + 1 - 2 k\xi_n^{3} + i k\xi_n^{2} + 1\xi_n - \frac{i}{2},\\
\frac{1}{\pi\beta_k}\int_{\mathbb R}\frac{P_{7}(\xi_n)\,\dd\xi_n}{(\xi_n-i)^{k+4}(\xi_n+i)^{k+1}}=- \frac{1}{2 \left(k + 1\right)}.
\end{gather*}
\begin{gather*}
M_{8}=g(N,\dot W) g(N,Z) g(X,Y),\\
P_{8}(\xi_n)=i \left(2 k - 1\right)\xi_n^{4} + 2 \left(2 k - 1\right)\xi_n^{3} - 2 i k\xi_n^{2} - 2\xi_n + i,\\
\frac{1}{\pi\beta_k}\int_{\mathbb R}\frac{P_{8}(\xi_n)\,\dd\xi_n}{(\xi_n-i)^{k+4}(\xi_n+i)^{k+1}}=\frac{1}{k + 1}.
\end{gather*}
\begin{gather*}
M_{9}=g(\nabla_{Z^T}W,N) g(X,Y),\\
P_{9}(\xi_n)=\frac{i \left(2 k - 1\right)}{2}\xi_n^{4} + 2 k - 1\xi_n^{3} - i k\xi_n^{2} - 1\xi_n + \frac{i}{2},\\
\frac{1}{\pi\beta_k}\int_{\mathbb R}\frac{P_{9}(\xi_n)\,\dd\xi_n}{(\xi_n-i)^{k+4}(\xi_n+i)^{k+1}}=\frac{1}{2 \left(k + 1\right)}.
\end{gather*}
\begin{gather*}
M_{10}=\mathcal D(W) g(N,Z) g(X,Y),\\
P_{10}(\xi_n)=\frac{i \left(2 k - 1\right)}{2}\xi_n^{4} + 2 k - 1\xi_n^{3} - i k\xi_n^{2} - 1\xi_n + \frac{i}{2},\\
\frac{1}{\pi\beta_k}\int_{\mathbb R}\frac{P_{10}(\xi_n)\,\dd\xi_n}{(\xi_n-i)^{k+4}(\xi_n+i)^{k+1}}=\frac{1}{2 \left(k + 1\right)}.
\end{gather*}
Multiplying each last line by its $M_\rho$ and summing gives the geometric case coefficient stated above. Thus cancellations inside this case, including any zero individual integrals, are retained explicitly.
\paragraph{Case (c).} Inserting $t_-$ yields
\begin{equation}\phi_{3,5}={\eta} G(E_HH+E_{P_1}P)+{\eta}\lambda_kL.\end{equation}
This case cancels the entire first-pair skew metric term of case (b). Dropping either subprincipal case would therefore give a false boundary obstruction to curvature recovery.

The coefficient identities needed for the sum are
\begin{equation}
 A_H+B_H+C_H+E_H=-\frac{k}{2(k+1)(k+2)},\qquad
 A_{P_1}+B_{P_1}+C_{P_1}+E_{P_1}=\frac{k}{k+2}.
\end{equation}
Both follow by putting the four explicit rational functions over $\Delta_k$. This completes the calculation of the mean-curvature part in arbitrary even dimension.

\paragraph{Projected integrand and individual residues.}
After the Hardy projection, the Clifford trace and the normalized tangential sphere integral, the complete scalar integrand for this case is
\begin{equation}
\mathfrak I_{3,5}(\xi_n)=\sum_{\rho}M_\rho\frac{P_\rho(\xi_n)}{(\xi_n-i)^{k+4}(\xi_n+i)^{k+3}}.
\end{equation}
The prefactor of the case, including its sign and any factor $1/2$ or $i$, is included here. The following list gives every nonzero geometric monomial and its numerator; terms with the same monomial have already been collected.
\begin{gather*}
M_{1}=\eta g(W,Z) g(X,Y),\\
\begin{aligned}P_{1}(\xi_n)&=\frac{k \left(2 k - 1\right) \left(2 k + 1\right)}{4}\xi_n^{5} - \frac{i k^{2} \left(2 k - 1\right)}{2}\xi_n^{4} + \frac{k \left(2 k - 1\right) \left(2 k + 1\right)}{2}\xi_n^{3}\\ &\quad {} - \frac{i k \left(2 k - 1\right) \left(4 k^{2} + 5 k + 3\right)}{2 \left(2 k + 1\right)}\xi_n^{2} - \frac{3 k \left(4 k + 3\right)}{4}\xi_n + \frac{i k \left(2 k - 1\right)}{2}\end{aligned},\\
\frac{1}{\pi\beta_k}\int_{\mathbb R}\frac{P_{1}(\xi_n)\,\dd\xi_n}{(\xi_n-i)^{k+4}(\xi_n+i)^{k+3}}=- \frac{k \left(8 k^{3} + 30 k^{2} + 25 k + 9\right)}{8 \left(k + 1\right) \left(k + 2\right) \left(k + 3\right)}.
\end{gather*}
\begin{gather*}
M_{2}=\eta g(N,W) g(N,Z) g(X,Y),\\
\begin{aligned}P_{2}(\xi_n)&=- \frac{k \left(2 k - 1\right) \left(2 k + 1\right)}{2}\xi_n^{5} - i k^{2}\xi_n^{4} - k \left(2 k - 1\right) \left(2 k + 1\right)\xi_n^{3}\\ &\quad {} - \frac{i k \left(4 k^{2} + 5 k + 3\right)}{2 k + 1}\xi_n^{2} + \frac{3 k \left(4 k + 3\right)}{2}\xi_n + i k\end{aligned},\\
\frac{1}{\pi\beta_k}\int_{\mathbb R}\frac{P_{2}(\xi_n)\,\dd\xi_n}{(\xi_n-i)^{k+4}(\xi_n+i)^{k+3}}=\frac{k \left(8 k^{2} + 19 k + 9\right)}{4 \left(k + 1\right) \left(k + 2\right) \left(k + 3\right)}.
\end{gather*}
\begin{gather*}
M_{3}=\eta g(N,W) g(N,X) g(Y,Z),\\
P_{3}(\xi_n)=\frac{k}{4}\xi_n^{5} + \frac{i k}{4}\xi_n^{4} + \frac{k^{2}}{2 \left(2 k + 1\right)}\xi_n^{3} + \frac{i k^{2}}{2 \left(2 k + 1\right)}\xi_n^{2} - \frac{k}{4 \left(2 k + 1\right)}\xi_n - \frac{i k}{4 \left(2 k + 1\right)},\\
\frac{1}{\pi\beta_k}\int_{\mathbb R}\frac{P_{3}(\xi_n)\,\dd\xi_n}{(\xi_n-i)^{k+4}(\xi_n+i)^{k+3}}=\frac{k}{4 \left(k + 1\right) \left(k + 2\right)}.
\end{gather*}
\begin{gather*}
M_{4}=\eta g(N,Y) g(N,Z) g(W,X),\\
P_{4}(\xi_n)=- \frac{k}{4}\xi_n^{5} - \frac{i k}{4}\xi_n^{4} - \frac{k^{2}}{2 \left(2 k + 1\right)}\xi_n^{3} - \frac{i k^{2}}{2 \left(2 k + 1\right)}\xi_n^{2} + \frac{k}{4 \left(2 k + 1\right)}\xi_n + \frac{i k}{4 \left(2 k + 1\right)},\\
\frac{1}{\pi\beta_k}\int_{\mathbb R}\frac{P_{4}(\xi_n)\,\dd\xi_n}{(\xi_n-i)^{k+4}(\xi_n+i)^{k+3}}=- \frac{k}{4 \left(k + 1\right) \left(k + 2\right)}.
\end{gather*}
\begin{gather*}
M_{5}=\eta g(N,W) g(N,Y) g(X,Z),\\
P_{5}(\xi_n)=- \frac{k}{4}\xi_n^{5} - \frac{i k}{4}\xi_n^{4} - \frac{k^{2}}{2 \left(2 k + 1\right)}\xi_n^{3} - \frac{i k^{2}}{2 \left(2 k + 1\right)}\xi_n^{2} + \frac{k}{4 \left(2 k + 1\right)}\xi_n + \frac{i k}{4 \left(2 k + 1\right)},\\
\frac{1}{\pi\beta_k}\int_{\mathbb R}\frac{P_{5}(\xi_n)\,\dd\xi_n}{(\xi_n-i)^{k+4}(\xi_n+i)^{k+3}}=- \frac{k}{4 \left(k + 1\right) \left(k + 2\right)}.
\end{gather*}
\begin{gather*}
M_{6}=\eta g(N,X) g(N,Z) g(W,Y),\\
P_{6}(\xi_n)=\frac{k}{4}\xi_n^{5} + \frac{i k}{4}\xi_n^{4} + \frac{k^{2}}{2 \left(2 k + 1\right)}\xi_n^{3} + \frac{i k^{2}}{2 \left(2 k + 1\right)}\xi_n^{2} - \frac{k}{4 \left(2 k + 1\right)}\xi_n - \frac{i k}{4 \left(2 k + 1\right)},\\
\frac{1}{\pi\beta_k}\int_{\mathbb R}\frac{P_{6}(\xi_n)\,\dd\xi_n}{(\xi_n-i)^{k+4}(\xi_n+i)^{k+3}}=\frac{k}{4 \left(k + 1\right) \left(k + 2\right)}.
\end{gather*}
Multiplying each last line by its $M_\rho$ and summing gives the geometric case coefficient stated above. Thus cancellations inside this case, including any zero individual integrals, are retained explicitly.
\subsection{The fourth insertion: all five cases}
For $A=\hc(X)\hc(Y)$ and $B=\hc(Z)\hc(W)$, the even hatted product commutes with $c(v)$, and
\[\tau_\Lambda(Ac(v)Bc(w))=-Jg(v,w).\]
We retain separately every derivative of $J$ and of its right-insertion part $J_B$.
\paragraph{Case (a.I).}

\[
 \Phi_{4,1}=-\int_{|\xi'|=1}\int_{\mathbb R}
 \sum_{a<n}\tr_\Lambda\left[
 \partial_{\xi_a}\pi_{\xi_n}^+\sigma_{-1}(P)
 \partial_{x_a}\partial_{\xi_n}\sigma_{3-n}(Q)\right]
 \dd\xi_n\,\sigma(\xi').
\]
At $x_0$, $\partial_{x_a}g^{ij}=0$ for $a<n$, and the chosen orthonormal trivialization gives $\partial_{x_a}c(\xi)=0$. Therefore
\[
 \partial_{x_a}\sigma_{3-n}(Q)
 =i\,\partial_{x_a}c(\xi)|\xi|^{-2k}
 -ikc(\xi)|\xi|^{-2k-2}\partial_{x_a}|\xi|^2=0.
\]
The derivative $\partial_{\xi_n}$ of this zero expression is also zero. This proves the vanishing before trace or integration and does not differentiate the inserted fields in the first factor.
The tangential derivative of the right principal symbol vanishes, so $\phi_{4,1}=0$.
\paragraph{Case (a.II).}
The derivative of the left principal symbol differentiates all four input fields as well as the metric symbol. It gives
\[\phi_{4,2}=-a_k\eta{}J+b_k\dot J.\]

\paragraph{Projected integrand and individual residues.}
After the Hardy projection, the Clifford trace and the normalized tangential sphere integral, the complete scalar integrand for this case is
\begin{equation}
\mathfrak I_{4,2}(\xi_n)=\sum_{\rho}M_\rho\frac{P_\rho(\xi_n)}{(\xi_n-i)^{k+3}(\xi_n+i)^{k+2}}.
\end{equation}
The prefactor of the case, including its sign and any factor $1/2$ or $i$, is included here. The following list gives every nonzero geometric monomial and its numerator; terms with the same monomial have already been collected.
\begin{gather*}
M_{1}=g(W,\dot X) g(Y,Z),\\
P_{1}(\xi_n)=- \frac{k \left(2 k - 1\right)}{2}\xi_n^{3} + \frac{i k \left(2 k + 1\right)}{2}\xi_n^{2} + \frac{3 k}{2}\xi_n - \frac{i k}{2},\\
\frac{1}{\pi\beta_k}\int_{\mathbb R}\frac{P_{1}(\xi_n)\,\dd\xi_n}{(\xi_n-i)^{k+3}(\xi_n+i)^{k+2}}=\frac{k}{2 \left(k + 1\right)}.
\end{gather*}
\begin{gather*}
M_{2}=g(W,Z) g(\dot X,Y),\\
P_{2}(\xi_n)=- \frac{k \left(2 k - 1\right)}{2}\xi_n^{3} + \frac{i k \left(2 k + 1\right)}{2}\xi_n^{2} + \frac{3 k}{2}\xi_n - \frac{i k}{2},\\
\frac{1}{\pi\beta_k}\int_{\mathbb R}\frac{P_{2}(\xi_n)\,\dd\xi_n}{(\xi_n-i)^{k+3}(\xi_n+i)^{k+2}}=\frac{k}{2 \left(k + 1\right)}.
\end{gather*}
\begin{gather*}
M_{3}=g(W,Y) g(\dot X,Z),\\
P_{3}(\xi_n)=\frac{k \left(2 k - 1\right)}{2}\xi_n^{3} - \frac{i k \left(2 k + 1\right)}{2}\xi_n^{2} - \frac{3 k}{2}\xi_n + \frac{i k}{2},\\
\frac{1}{\pi\beta_k}\int_{\mathbb R}\frac{P_{3}(\xi_n)\,\dd\xi_n}{(\xi_n-i)^{k+3}(\xi_n+i)^{k+2}}=- \frac{k}{2 \left(k + 1\right)}.
\end{gather*}
\begin{gather*}
M_{4}=g(W,X) g(\dot Y,Z),\\
P_{4}(\xi_n)=- \frac{k \left(2 k - 1\right)}{2}\xi_n^{3} + \frac{i k \left(2 k + 1\right)}{2}\xi_n^{2} + \frac{3 k}{2}\xi_n - \frac{i k}{2},\\
\frac{1}{\pi\beta_k}\int_{\mathbb R}\frac{P_{4}(\xi_n)\,\dd\xi_n}{(\xi_n-i)^{k+3}(\xi_n+i)^{k+2}}=\frac{k}{2 \left(k + 1\right)}.
\end{gather*}
\begin{gather*}
M_{5}=g(W,Z) g(X,\dot Y),\\
P_{5}(\xi_n)=- \frac{k \left(2 k - 1\right)}{2}\xi_n^{3} + \frac{i k \left(2 k + 1\right)}{2}\xi_n^{2} + \frac{3 k}{2}\xi_n - \frac{i k}{2},\\
\frac{1}{\pi\beta_k}\int_{\mathbb R}\frac{P_{5}(\xi_n)\,\dd\xi_n}{(\xi_n-i)^{k+3}(\xi_n+i)^{k+2}}=\frac{k}{2 \left(k + 1\right)}.
\end{gather*}
\begin{gather*}
M_{6}=g(W,\dot Y) g(X,Z),\\
P_{6}(\xi_n)=\frac{k \left(2 k - 1\right)}{2}\xi_n^{3} - \frac{i k \left(2 k + 1\right)}{2}\xi_n^{2} - \frac{3 k}{2}\xi_n + \frac{i k}{2},\\
\frac{1}{\pi\beta_k}\int_{\mathbb R}\frac{P_{6}(\xi_n)\,\dd\xi_n}{(\xi_n-i)^{k+3}(\xi_n+i)^{k+2}}=- \frac{k}{2 \left(k + 1\right)}.
\end{gather*}
\begin{gather*}
M_{7}=\eta g(W,Y) g(X,Z),\\
P_{7}(\xi_n)=\frac{i k \left(2 k - 1\right)}{4}\xi_n^{2} + \frac{k}{2}\xi_n - \frac{i k}{4},\\
\frac{1}{\pi\beta_k}\int_{\mathbb R}\frac{P_{7}(\xi_n)\,\dd\xi_n}{(\xi_n-i)^{k+3}(\xi_n+i)^{k+2}}=\frac{3 k \left(2 k + 1\right)}{8 \left(k + 1\right) \left(k + 2\right)}.
\end{gather*}
\begin{gather*}
M_{8}=\eta g(W,X) g(Y,Z),\\
P_{8}(\xi_n)=- \frac{i k \left(2 k - 1\right)}{4}\xi_n^{2} - \frac{k}{2}\xi_n + \frac{i k}{4},\\
\frac{1}{\pi\beta_k}\int_{\mathbb R}\frac{P_{8}(\xi_n)\,\dd\xi_n}{(\xi_n-i)^{k+3}(\xi_n+i)^{k+2}}=- \frac{3 k \left(2 k + 1\right)}{8 \left(k + 1\right) \left(k + 2\right)}.
\end{gather*}
\begin{gather*}
M_{9}=\eta g(W,Z) g(X,Y),\\
P_{9}(\xi_n)=- \frac{i k \left(2 k - 1\right)}{4}\xi_n^{2} - \frac{k}{2}\xi_n + \frac{i k}{4},\\
\frac{1}{\pi\beta_k}\int_{\mathbb R}\frac{P_{9}(\xi_n)\,\dd\xi_n}{(\xi_n-i)^{k+3}(\xi_n+i)^{k+2}}=- \frac{3 k \left(2 k + 1\right)}{8 \left(k + 1\right) \left(k + 2\right)}.
\end{gather*}
\begin{gather*}
M_{10}=g(W,X) g(Y,\dot Z),\\
P_{10}(\xi_n)=- \frac{k \left(2 k - 1\right)}{2}\xi_n^{3} + \frac{i k \left(2 k + 1\right)}{2}\xi_n^{2} + \frac{3 k}{2}\xi_n - \frac{i k}{2},\\
\frac{1}{\pi\beta_k}\int_{\mathbb R}\frac{P_{10}(\xi_n)\,\dd\xi_n}{(\xi_n-i)^{k+3}(\xi_n+i)^{k+2}}=\frac{k}{2 \left(k + 1\right)}.
\end{gather*}
\begin{gather*}
M_{11}=g(W,\dot Z) g(X,Y),\\
P_{11}(\xi_n)=- \frac{k \left(2 k - 1\right)}{2}\xi_n^{3} + \frac{i k \left(2 k + 1\right)}{2}\xi_n^{2} + \frac{3 k}{2}\xi_n - \frac{i k}{2},\\
\frac{1}{\pi\beta_k}\int_{\mathbb R}\frac{P_{11}(\xi_n)\,\dd\xi_n}{(\xi_n-i)^{k+3}(\xi_n+i)^{k+2}}=\frac{k}{2 \left(k + 1\right)}.
\end{gather*}
\begin{gather*}
M_{12}=g(W,Y) g(X,\dot Z),\\
P_{12}(\xi_n)=\frac{k \left(2 k - 1\right)}{2}\xi_n^{3} - \frac{i k \left(2 k + 1\right)}{2}\xi_n^{2} - \frac{3 k}{2}\xi_n + \frac{i k}{2},\\
\frac{1}{\pi\beta_k}\int_{\mathbb R}\frac{P_{12}(\xi_n)\,\dd\xi_n}{(\xi_n-i)^{k+3}(\xi_n+i)^{k+2}}=- \frac{k}{2 \left(k + 1\right)}.
\end{gather*}
\begin{gather*}
M_{13}=g(\dot W,X) g(Y,Z),\\
P_{13}(\xi_n)=- \frac{k \left(2 k - 1\right)}{2}\xi_n^{3} + \frac{i k \left(2 k + 1\right)}{2}\xi_n^{2} + \frac{3 k}{2}\xi_n - \frac{i k}{2},\\
\frac{1}{\pi\beta_k}\int_{\mathbb R}\frac{P_{13}(\xi_n)\,\dd\xi_n}{(\xi_n-i)^{k+3}(\xi_n+i)^{k+2}}=\frac{k}{2 \left(k + 1\right)}.
\end{gather*}
\begin{gather*}
M_{14}=g(\dot W,Z) g(X,Y),\\
P_{14}(\xi_n)=- \frac{k \left(2 k - 1\right)}{2}\xi_n^{3} + \frac{i k \left(2 k + 1\right)}{2}\xi_n^{2} + \frac{3 k}{2}\xi_n - \frac{i k}{2},\\
\frac{1}{\pi\beta_k}\int_{\mathbb R}\frac{P_{14}(\xi_n)\,\dd\xi_n}{(\xi_n-i)^{k+3}(\xi_n+i)^{k+2}}=\frac{k}{2 \left(k + 1\right)}.
\end{gather*}
\begin{gather*}
M_{15}=g(\dot W,Y) g(X,Z),\\
P_{15}(\xi_n)=\frac{k \left(2 k - 1\right)}{2}\xi_n^{3} - \frac{i k \left(2 k + 1\right)}{2}\xi_n^{2} - \frac{3 k}{2}\xi_n + \frac{i k}{2},\\
\frac{1}{\pi\beta_k}\int_{\mathbb R}\frac{P_{15}(\xi_n)\,\dd\xi_n}{(\xi_n-i)^{k+3}(\xi_n+i)^{k+2}}=- \frac{k}{2 \left(k + 1\right)}.
\end{gather*}
Multiplying each last line by its $M_\rho$ and summing gives the geometric case coefficient stated above. Thus cancellations inside this case, including any zero individual integrals, are retained explicitly.
\paragraph{Case (a.III).}
Only the right principal symbol is differentiated in the normal coordinate. Its contribution is
\[\phi_{4,3}=c_k\eta{}J.\]

\paragraph{Projected integrand and individual residues.}
After the Hardy projection, the Clifford trace and the normalized tangential sphere integral, the complete scalar integrand for this case is
\begin{equation}
\mathfrak I_{4,3}(\xi_n)=\sum_{\rho}M_\rho\frac{P_\rho(\xi_n)}{(\xi_n-i)^{k+3}(\xi_n+i)^{k+2}}.
\end{equation}
The prefactor of the case, including its sign and any factor $1/2$ or $i$, is included here. The following list gives every nonzero geometric monomial and its numerator; terms with the same monomial have already been collected.
\begin{gather*}
M_{1}=\eta g(W,Y) g(X,Z),\\
P_{1}(\xi_n)=- \frac{i k}{4}\xi_n^{2} - \frac{k^{2}}{2}\xi_n + \frac{i k}{4},\\
\frac{1}{\pi\beta_k}\int_{\mathbb R}\frac{P_{1}(\xi_n)\,\dd\xi_n}{(\xi_n-i)^{k+3}(\xi_n+i)^{k+2}}=- \frac{k \left(2 k + 1\right)^{2}}{8 \left(k + 1\right) \left(k + 2\right)}.
\end{gather*}
\begin{gather*}
M_{2}=\eta g(W,X) g(Y,Z),\\
P_{2}(\xi_n)=\frac{i k}{4}\xi_n^{2} + \frac{k^{2}}{2}\xi_n - \frac{i k}{4},\\
\frac{1}{\pi\beta_k}\int_{\mathbb R}\frac{P_{2}(\xi_n)\,\dd\xi_n}{(\xi_n-i)^{k+3}(\xi_n+i)^{k+2}}=\frac{k \left(2 k + 1\right)^{2}}{8 \left(k + 1\right) \left(k + 2\right)}.
\end{gather*}
\begin{gather*}
M_{3}=\eta g(W,Z) g(X,Y),\\
P_{3}(\xi_n)=\frac{i k}{4}\xi_n^{2} + \frac{k^{2}}{2}\xi_n - \frac{i k}{4},\\
\frac{1}{\pi\beta_k}\int_{\mathbb R}\frac{P_{3}(\xi_n)\,\dd\xi_n}{(\xi_n-i)^{k+3}(\xi_n+i)^{k+2}}=\frac{k \left(2 k + 1\right)^{2}}{8 \left(k + 1\right) \left(k + 2\right)}.
\end{gather*}
Multiplying each last line by its $M_\rho$ and summing gives the geometric case coefficient stated above. Thus cancellations inside this case, including any zero individual integrals, are retained explicitly.
\paragraph{Case (b).}
The subprincipal left factor includes $B_{,j}=\hc(Z_{,j})\hc(W)+\hc(Z)\hc(W_{,j})$. The normal derivatives give $-t_kJ_B$ and the connection and metric terms give
\[\phi_{4,4}=d_k\eta{}J-t_kJ_B.\]
The tangential derivative terms have a zero unhatted normal--tangential contraction.

\paragraph{Projected integrand and individual residues.}
After the Hardy projection, the Clifford trace and the normalized tangential sphere integral, the complete scalar integrand for this case is
\begin{equation}
\mathfrak I_{4,4}(\xi_n)=\sum_{\rho}M_\rho\frac{P_\rho(\xi_n)}{(\xi_n-i)^{k+3}(\xi_n+i)^{k+1}}.
\end{equation}
The prefactor of the case, including its sign and any factor $1/2$ or $i$, is included here. The following list gives every nonzero geometric monomial and its numerator; terms with the same monomial have already been collected.
\begin{gather*}
M_{1}=\eta g(W,X) g(Y,Z),\\
P_{1}(\xi_n)=\frac{\left(2 k - 1\right) \left(2 k + 1\right)}{8}\xi_n^{2} - \frac{i \left(2 k^{2} + 3 k - 1\right)}{4}\xi_n - \frac{2 k + 3}{8},\\
\frac{1}{\pi\beta_k}\int_{\mathbb R}\frac{P_{1}(\xi_n)\,\dd\xi_n}{(\xi_n-i)^{k+3}(\xi_n+i)^{k+1}}=\frac{k \left(2 k + 1\right) \left(2 k + 7\right)}{8 \left(k + 1\right) \left(k + 2\right)}.
\end{gather*}
\begin{gather*}
M_{2}=\eta g(W,Z) g(X,Y),\\
P_{2}(\xi_n)=\frac{\left(2 k - 1\right) \left(2 k + 1\right)}{8}\xi_n^{2} - \frac{i \left(2 k^{2} + 3 k - 1\right)}{4}\xi_n - \frac{2 k + 3}{8},\\
\frac{1}{\pi\beta_k}\int_{\mathbb R}\frac{P_{2}(\xi_n)\,\dd\xi_n}{(\xi_n-i)^{k+3}(\xi_n+i)^{k+1}}=\frac{k \left(2 k + 1\right) \left(2 k + 7\right)}{8 \left(k + 1\right) \left(k + 2\right)}.
\end{gather*}
\begin{gather*}
M_{3}=\eta g(W,Y) g(X,Z),\\
P_{3}(\xi_n)=- \frac{\left(2 k - 1\right) \left(2 k + 1\right)}{8}\xi_n^{2} + \frac{i \left(2 k^{2} + 3 k - 1\right)}{4}\xi_n + \frac{2 k + 3}{8},\\
\frac{1}{\pi\beta_k}\int_{\mathbb R}\frac{P_{3}(\xi_n)\,\dd\xi_n}{(\xi_n-i)^{k+3}(\xi_n+i)^{k+1}}=- \frac{k \left(2 k + 1\right) \left(2 k + 7\right)}{8 \left(k + 1\right) \left(k + 2\right)}.
\end{gather*}
\begin{gather*}
M_{4}=g(W,Y) g(X,\dot Z),\\
P_{4}(\xi_n)=\frac{i \left(2 k - 1\right)}{2}\xi_n^{3} + \frac{2 k - 1}{2}\xi_n^{2} - \frac{i}{2}\xi_n - \frac{1}{2},\\
\frac{1}{\pi\beta_k}\int_{\mathbb R}\frac{P_{4}(\xi_n)\,\dd\xi_n}{(\xi_n-i)^{k+3}(\xi_n+i)^{k+1}}=\frac{1}{2 \left(k + 1\right)}.
\end{gather*}
\begin{gather*}
M_{5}=g(W,X) g(Y,\dot Z),\\
P_{5}(\xi_n)=- \frac{i \left(2 k - 1\right)}{2}\xi_n^{3} - \frac{2 k - 1}{2}\xi_n^{2} + \frac{i}{2}\xi_n + \frac{1}{2},\\
\frac{1}{\pi\beta_k}\int_{\mathbb R}\frac{P_{5}(\xi_n)\,\dd\xi_n}{(\xi_n-i)^{k+3}(\xi_n+i)^{k+1}}=- \frac{1}{2 \left(k + 1\right)}.
\end{gather*}
\begin{gather*}
M_{6}=g(W,\dot Z) g(X,Y),\\
P_{6}(\xi_n)=- \frac{i \left(2 k - 1\right)}{2}\xi_n^{3} - \frac{2 k - 1}{2}\xi_n^{2} + \frac{i}{2}\xi_n + \frac{1}{2},\\
\frac{1}{\pi\beta_k}\int_{\mathbb R}\frac{P_{6}(\xi_n)\,\dd\xi_n}{(\xi_n-i)^{k+3}(\xi_n+i)^{k+1}}=- \frac{1}{2 \left(k + 1\right)}.
\end{gather*}
\begin{gather*}
M_{7}=g(\dot W,Y) g(X,Z),\\
P_{7}(\xi_n)=\frac{i \left(2 k - 1\right)}{2}\xi_n^{3} + \frac{2 k - 1}{2}\xi_n^{2} - \frac{i}{2}\xi_n - \frac{1}{2},\\
\frac{1}{\pi\beta_k}\int_{\mathbb R}\frac{P_{7}(\xi_n)\,\dd\xi_n}{(\xi_n-i)^{k+3}(\xi_n+i)^{k+1}}=\frac{1}{2 \left(k + 1\right)}.
\end{gather*}
\begin{gather*}
M_{8}=g(\dot W,X) g(Y,Z),\\
P_{8}(\xi_n)=- \frac{i \left(2 k - 1\right)}{2}\xi_n^{3} - \frac{2 k - 1}{2}\xi_n^{2} + \frac{i}{2}\xi_n + \frac{1}{2},\\
\frac{1}{\pi\beta_k}\int_{\mathbb R}\frac{P_{8}(\xi_n)\,\dd\xi_n}{(\xi_n-i)^{k+3}(\xi_n+i)^{k+1}}=- \frac{1}{2 \left(k + 1\right)}.
\end{gather*}
\begin{gather*}
M_{9}=g(\dot W,Z) g(X,Y),\\
P_{9}(\xi_n)=- \frac{i \left(2 k - 1\right)}{2}\xi_n^{3} - \frac{2 k - 1}{2}\xi_n^{2} + \frac{i}{2}\xi_n + \frac{1}{2},\\
\frac{1}{\pi\beta_k}\int_{\mathbb R}\frac{P_{9}(\xi_n)\,\dd\xi_n}{(\xi_n-i)^{k+3}(\xi_n+i)^{k+1}}=- \frac{1}{2 \left(k + 1\right)}.
\end{gather*}
Multiplying each last line by its $M_\rho$ and summing gives the geometric case coefficient stated above. Thus cancellations inside this case, including any zero individual integrals, are retained explicitly.
\paragraph{Case (c).}
The subprincipal right factor gives
\[\phi_{4,5}=-e_k\eta{}J.\]
The identity $-a_k+c_k+d_k-e_k=0$ cancels the total explicit metric contribution; the final derivative density is $b_k\dot J-t_kJ_B$.

\paragraph{Projected integrand and individual residues.}
After the Hardy projection, the Clifford trace and the normalized tangential sphere integral, the complete scalar integrand for this case is
\begin{equation}
\mathfrak I_{4,5}(\xi_n)=\sum_{\rho}M_\rho\frac{P_\rho(\xi_n)}{(\xi_n-i)^{k+3}(\xi_n+i)^{k+3}}.
\end{equation}
The prefactor of the case, including its sign and any factor $1/2$ or $i$, is included here. The following list gives every nonzero geometric monomial and its numerator; terms with the same monomial have already been collected.
\begin{gather*}
M_{1}=\eta g(W,Y) g(X,Z),\\
\begin{aligned}P_{1}(\xi_n)&=- \frac{k \left(2 k - 1\right) \left(2 k + 1\right)}{4}\xi_n^{4} + \frac{i k \left(2 k + 1\right)}{4}\xi_n^{3}\\ &\quad {} - \frac{k \left(2 k + 1\right) \left(4 k - 1\right)}{4}\xi_n^{2} + \frac{i k \left(8 k + 7\right)}{4}\xi_n + \frac{k \left(2 k + 1\right)}{2}\end{aligned},\\
\frac{1}{\pi\beta_k}\int_{\mathbb R}\frac{P_{1}(\xi_n)\,\dd\xi_n}{(\xi_n-i)^{k+3}(\xi_n+i)^{k+3}}=\frac{k \left(2 k + 1\right) \left(4 k + 5\right)}{8 \left(k + 1\right) \left(k + 2\right)}.
\end{gather*}
\begin{gather*}
M_{2}=\eta g(W,X) g(Y,Z),\\
\begin{aligned}P_{2}(\xi_n)&=\frac{k \left(2 k - 1\right) \left(2 k + 1\right)}{4}\xi_n^{4} - \frac{i k \left(2 k + 1\right)}{4}\xi_n^{3}\\ &\quad {} + \frac{k \left(2 k + 1\right) \left(4 k - 1\right)}{4}\xi_n^{2} - \frac{i k \left(8 k + 7\right)}{4}\xi_n - \frac{k \left(2 k + 1\right)}{2}\end{aligned},\\
\frac{1}{\pi\beta_k}\int_{\mathbb R}\frac{P_{2}(\xi_n)\,\dd\xi_n}{(\xi_n-i)^{k+3}(\xi_n+i)^{k+3}}=- \frac{k \left(2 k + 1\right) \left(4 k + 5\right)}{8 \left(k + 1\right) \left(k + 2\right)}.
\end{gather*}
\begin{gather*}
M_{3}=\eta g(W,Z) g(X,Y),\\
\begin{aligned}P_{3}(\xi_n)&=\frac{k \left(2 k - 1\right) \left(2 k + 1\right)}{4}\xi_n^{4} - \frac{i k \left(2 k + 1\right)}{4}\xi_n^{3}\\ &\quad {} + \frac{k \left(2 k + 1\right) \left(4 k - 1\right)}{4}\xi_n^{2} - \frac{i k \left(8 k + 7\right)}{4}\xi_n - \frac{k \left(2 k + 1\right)}{2}\end{aligned},\\
\frac{1}{\pi\beta_k}\int_{\mathbb R}\frac{P_{3}(\xi_n)\,\dd\xi_n}{(\xi_n-i)^{k+3}(\xi_n+i)^{k+3}}=- \frac{k \left(2 k + 1\right) \left(4 k + 5\right)}{8 \left(k + 1\right) \left(k + 2\right)}.
\end{gather*}
Multiplying each last line by its $M_\rho$ and summing gives the geometric case coefficient stated above. Thus cancellations inside this case, including any zero individual integrals, are retained explicitly.
\subsection{The final geometric boundary terms}
Equations~\eqref{eq:boundaryfinal}--\eqref{eq:Qfinite} can now be written without symbols, Clifford traces, angular averages or contour coefficients:
\begin{equation}\label{eq:geomtheorem}
 \boxed{\widetilde{\Wres}(\mathcal B_\ell)
 =C_m\int_M I_\ell\dd V_g+b_n\int_{\partial M}S_\ell\dd V_{\partial M}.}
\end{equation}
The four geometric densities are
\begin{align}\label{eq:geomsums12}
 S_1&=-b_k\dot G+t_kg(X,\dot Y),\\
 S_2&=-b_k\partial_n(GH)+t_kG(\dot H+T_-),\notag
\end{align}
\begin{align}\label{eq:geomsum3}
 S_3={}&u_k\partial_n(GH)-v_k\partial_n(GP)\notag\\
 &+t_kG\{T_++2Z_ng(\dot W,N)-g(Z,\dot W)\}\notag\\
 &+\frac{kK_{\partial}}{d(k+1)(k+2)}G\{H-2(k+1)P_1\},
\end{align}
and
\begin{equation}\label{eq:geomsum4}
 S_4=b_k\dot J-t_kJ_B.
\end{equation}
All quantities in these formulas were defined geometrically above. If desired, $H=g(Z^T,W^T)+P_1$ separates tangential and normal components in the mean-curvature term. No condition of parallel extension, tangential constancy or vanishing normal derivative has been imposed.

For $n=6$, $k=2$, $\beta_2=3/8$, $b_6=1/\pi^3$, and the formulas specialize to
\begin{align}
 S_1&=-\tfrac13\dot G+\tfrac16g(X,\dot Y),\\
 S_2&=-\tfrac13\partial_n(GH)+\tfrac16G(\dot H+T_-),\\
 S_3&=\tfrac14\partial_n(GH)-\tfrac16\partial_n(GP)
       +\tfrac16G\{T_++2Z_ng(\dot W,N)-g(Z,\dot W)\}\notag\\
 &\qquad+\tfrac1{30}K_{\partial}G(H-6P_1),\\
 S_4&=\tfrac13\dot J-\tfrac16J_B.
\end{align}
These are the geometric form of the six-dimensional result, while~\eqref{eq:geomtheorem} retains arbitrary even dimension. At $k=1$ all denominators remain nonzero, so dimension four requires no limiting convention.

As a non-product-collar check, prescribe vanishing covariant first jets at a point. Then $S_1=S_2=S_4=0$, whereas
\[
 S_3=\frac{kK_{\partial}}{d(k+1)(k+2)}G\{H-2(k+1)P_1\}.
\]
For $X=Y=e_1$, $Z=W=e_2$ tangential and ${\eta}=1$, the coefficients $\Phi_3/\pi$ in dimensions $4,6,8$ are $-1/24,-1/32,-3/128$. These agree with the direct five-case calculation. An independent implementation using exterior Clifford bit masks was compared with the geometric formulas for 300 individual case values in dimensions $4,6,8$, including nonzero normal and tangential field jets. These tests supplement the symbolic derivation; they are not the source of its dimension dependence.

\subsection{The exceptional two-dimensional order}
For $n=2$ the right factor of~\eqref{eq:boundarydef} would be $D$, a differential operator of positive order. It is not the negative-order class-zero factor used in Theorem~\ref{thm:boundary}. We do not claim that the negative-order boundary theorem applies unchanged. If one instead defines the two-dimensional extension by the single truncation $(A_\ell D B_\ell D^{-1})_+$ with zero prescribed Green part, its residue is zero, because its interior residue density is zero. This is an explicit different convention. Thus the general boundary theorem covers every even dimension in the nondegenerate curvature range $n\geq4$; the surface case and its factorization issue are recorded rather than concealed.

\section{Boundary extensions of the Riemann curvature functional}\label{sec:commboundary}
\subsection{Two factorizations and the precise geometric statement}
We now extend Corollary~\ref{cor:R} using the two inserted products suggested by the computations of Wu and Wang~\cite{WuWangInserted}. The notation and the five-case organization follow that reference, but its unhatted Clifford insertion is replaced here by a mixed double commutator. The ordinary exterior trace, of rank $r=2^n$, is essential. Throughout this section the collar metric is exactly~\eqref{eq:collar}, with inward normal $N$; no constancy or tangentiality of the four input fields is assumed.

Put $n=2m$, $k=m-1$, $d=n-1=2k+1$, and
\begin{equation}\label{eq:Tdouble}
 T=[\hc(U),\hc(V)][c(W),c(Z)].
\end{equation}
The order of $W,Z$ in this formula is the order of the last pair in the resulting Riemann tensor. For $m\geq2$ define
\begin{align}
 \mathcal B_1(T)&=\widetilde{\Wres}\big[\pi^+(TD^{-1})\circ\pi^+(D^{-2m+3})\big],\label{eq:Bcomm1}\\
 \mathcal B_2(T)&=\widetilde{\Wres}\big[\pi^+(TD^{-2})\circ\pi^+(D^{-2m+4})\big].\label{eq:Bcomm2}
\end{align}
Here $\widetilde{\Wres}$ has the same normalization as Section~\ref{sec:boundary}. Both interior products equal $TD^{2-n}$ modulo smoothing operators. This observation identifies their interior densities; it does not identify their boundary densities, which are calculated separately below.

\begin{theorem}[Two boundary extensions of the curvature residue]\label{thm:doubleboundary}
For the warped collar above, in every even dimension $n\geq4$,
\begin{equation}\label{eq:doubleanswer}
 \boxed{\mathcal B_1(T)=\mathcal B_2(T)
 =(n-2)C_m\int_M R(U,V,W,Z)\,\dd V_g.}
\end{equation}
For each factorization the full boundary density vanishes pointwise. The same identities hold with an arbitrary smooth scalar weight $f$ multiplying $T$ on the far left, with $f$ in the interior integrand. Consequently either functional divided by $n-2$ extends Corollary~\ref{cor:R} to this class of manifolds with boundary and recovers the pointwise Riemann tensor through scalar localization.
\end{theorem}
The conclusion is a zero boundary correction, not an omitted boundary calculation. For the odd-power factorization, two case integrands below are generally nonzero; their normal integrals vanish separately. For the even-power factorization, all five traces vanish before the normal integration. These are different mechanisms.

\subsection{Interior density and the trace identities needed at the boundary}
Let $H=[\hc(U),\hc(V)]$ and $K=[c(W),c(Z)]$. The scalar components of both commutators vanish. Since mixed hatted--unhatted contractions are zero, Clifford trace reduction gives
\begin{equation}\label{eq:selection}
 \tau_\Lambda(TF(c))=0,\qquad
 \tau_\Lambda(TF(\hc))=0
\end{equation}
for any word involving only the indicated Clifford family. Derivatives of $T$ satisfy the corresponding first identity when the accompanying word has only unhatted generators. Indeed differentiating $H$ gives a sum of hatted commutators, each still with zero scalar component. This includes derivatives of the field components and of the metric-dependent actions in the collar frame.

For arbitrary $A,B$, direct use of the four-generator trace gives
\begin{align}
 \tau_\Lambda(H\hc(A)\hc(B))
 &=2\{g(U,B)g(V,A)-g(U,A)g(V,B)\},\label{eq:hatbivtrace}\\
 \tau_\Lambda(Kc(A)c(B))
 &=2\{g(W,B)g(Z,A)-g(W,A)g(Z,B)\}.\label{eq:cbivtrace}
\end{align}
The mixed trace of the product of these two even words is the product of these expressions. Thus scalar-curvature terms in the static contraction vanish. Substituting the mixed part of $E$ from~\eqref{eq:E}, and using skew symmetry in both pairs, gives
\begin{equation}
 \tau_\Lambda(TE)=-2R(U,V,W,Z),\qquad
 \tau_\Lambda\big[T(s/6-E)\big]=2R(U,V,W,Z).
\end{equation}
Formula~\eqref{eq:static} therefore gives $2kC_mR(U,V,W,Z)\dd V_g$, proving the common interior term in~\eqref{eq:doubleanswer}.

At the boundary put $U_n=g(U,N)$ and similarly for the other inputs. Introduce tangential one-forms
\[
 A_a=g(U,e_a)V_n-g(V,e_a)U_n,\qquad
 B_a=g(W,e_a)Z_n-g(Z,e_a)W_n\quad(a<n).
\]
Their contraction is the geometric quantity
\begin{align}\label{eq:Gmixed}
 \mathcal G:=\sum_{a<n}A_aB_a
 &=V_nZ_n g(U^T,W^T)-V_nW_n g(U^T,Z^T)\notag\\
 &\quad-U_nZ_n g(V^T,W^T)+U_nW_n g(V^T,Z^T).
\end{align}
Replacing the tangential scalar products by full scalar products leaves the same expression: the four pure-normal terms cancel. Equations~\eqref{eq:hatbivtrace}--\eqref{eq:cbivtrace} give, in particular,
\[
 \tau_\Lambda(H\hc(e_a)\hc(N))=-2A_a,\quad
 \tau_\Lambda(Kc(e_b)c(N))=-2B_b,
\]
and hence $\tau_\Lambda(Tc(e_b)c(N)\hc(e_a)\hc(N))=4A_aB_b$. These identities fix both the sign and the factor four in the calculation.

\subsection{Symbols for the odd-power factorization}
Recall $\eta=h'(0)$, $z=\xi_n$, $q=1+z^2$, $C=c(\xi')+zc(N)$ and $N_c=c(N)$. Thus $C_{x_n}=\eta c(\xi')/2$ at $|\xi'|=1$. Distinguish $N_c$, a Clifford matrix, from the normal vector $N$. The zero-order symbol of $D$ is
\begin{equation}\label{eq:wcomm}
 w=-\frac{d\eta}{4}N_c-\frac\eta4\sum_{a<n}c(e_a)\hc(e_a)\hc(N).
\end{equation}
Let $\ell_3$ be the explicit inverse-square subprincipal symbol in~\eqref{eq:ell3}. For $P=TD^{-1}$ and $Q=D^{3-n}=DL^{-k}$, the required symbols are
\begin{align}
 p&=iTCq^{-1},\qquad t=iCq^{-k},\label{eq:oddlead}\\
 p_-&=iTC\ell_3+Twq^{-1}-\eta TN_cq^{-2},\label{eq:oddleftsub}\\
 t_-&=ikCq^{-k+1}\ell_3+k(k-1)\eta zCq^{-k-2}
       +wq^{-k}-k\eta N_cq^{-k-1},\label{eq:oddrightsub}\\
 p_{x_n}&=iT_{,n}Cq^{-1}+iTC_{x_n}q^{-1}-i\eta TCq^{-2},\\
 t_{x_n}&=iC_{x_n}q^{-k}-ik\eta Cq^{-k-1}.
\end{align}
There is no derivative of $T$ in $p_-$: multiplication by $T$ is on the far left. There is a derivative of $T$ in $p_{x_n}$, and it must be retained until the trace is taken. The term $-\eta TN_cq^{-2}$ is the symbol-composition correction, not a connection term.

For clarity define $\langle F\rangle$ as normalized tangential sphere integration followed by $\tau_\Lambda$:
\[
 \langle F\rangle=\frac1{\vol(S^{n-2})}\int_{|\xi'|=1}\tau_\Lambda(F)\,\sigma(\xi').
\]
Thus the geometric case densities below must be multiplied by $r\vol(S^{n-2})/(2\pi)^n$. We do not include $\pi\beta_k$ in this bracket. The moments used are $\langle\xi_a\rangle=0$ and $\langle\xi_a\xi_b\rangle=\delta_{ab}/d$ for scalar functions; only these moments are needed here.

\subsection{Odd powers: case (a.I), tangential differentiation}
The order constraint is $r+\ell-j-v-|\alpha|-1=-n$, with leading orders $r=-1$, $\ell=3-n$. As in Section~\ref{sec:boundary}, exactly one loss remains. The first case is
\[
 \Phi^{(1)}_{a.I}=-\int_{\mathbb R}\sum_{a<n}
 \left\langle \partial_{\xi_a}\pi_z^+p\,
 \partial_{x_a}\partial_z t\right\rangle\,\dd z.
\]
In geodesic boundary coordinates, with the specified frame, $\partial_{x_a}q=0$ and $\partial_{x_a}C=0$ at the center. Therefore
\[
 \partial_{x_a}t=i(\partial_{x_a}C)q^{-k}
 -ikCq^{-k-1}\partial_{x_a}q=0,
 \qquad \Phi^{(1)}_{a.I}=0.
\]
The right factor contains no input fields. This vanishing imposes no condition on their tangential derivatives.

\subsection{Odd powers: case (a.II), normal differentiation of the left factor}
Here
\[
 \Phi^{(1)}_{a.II}=-\frac12\int_{\mathbb R}
 \left\langle\pi_z^+(p_{x_n})\partial_z^2t\right\rangle\,\dd z.
\]
Differentiate the right principal symbol explicitly:
\begin{align}
 \partial_z t&=iN_cq^{-k}-2ikzCq^{-k-1},\\
 \partial_z^2t&=-4ikzN_cq^{-k-1}-2ikCq^{-k-1}
                  +4ik(k+1)z^2Cq^{-k-2}.
\end{align}
Each product with the three terms of $p_{x_n}$ is a scalar rational function times either $T$ or $T_{,n}$ followed by unhatted Clifford matrices only. Its trace is zero by~\eqref{eq:selection}, before angular or normal integration. The projection acts on the scalar rational factors and does not alter the Clifford selection rule. Thus $\Phi^{(1)}_{a.II}=0$, including all first derivatives of the inserted vector fields.

\subsection{Odd powers: case (a.III), normal differentiation of the right factor}
The contribution is
\[
 \Phi^{(1)}_{a.III}=-\frac12\int_{\mathbb R}
 \left\langle\partial_z\pi_z^+p\,\partial_z t_{x_n}\right\rangle\,\dd z.
\]
The required derivatives are
\begin{align}
 \partial_z p&=iTN_cq^{-1}-2izTCq^{-2},\\
 \partial_z t_{x_n}&=-2ikzC_{x_n}q^{-k-1}
 -ik\eta N_cq^{-k-1}+2ik(k+1)\eta zCq^{-k-2}.
\end{align}
Again every matrix word after $T$ is unhatted. Its hatted commutator has zero scalar component, so its trace vanishes identically. Therefore $\Phi^{(1)}_{a.III}=0$. This proof is independent of the numerical values of the scalar coefficients.

\subsection{Odd powers: case (b), the left subprincipal symbol}
This is the first case in which the hatted part of the connection symbol can pair with $H$. We must compute rather than apply~\eqref{eq:selection} to the entire expression:
\begin{equation}
 \Phi^{(1)}_b=-i\int_{\mathbb R}\left\langle\pi_z^+p_-\,\partial_z t\right\rangle\,\dd z.
\end{equation}
To keep projection separate from multiplication, let $v$ denote the normal variable of the first factor and retain $z$ in the second. Put $q_z=1+z^2$. Substitution of~\eqref{eq:oddleftsub},~\eqref{eq:ell3} and~\eqref{eq:wcomm}, followed by the two bivector traces and the second angular moment, gives
\begin{equation}\label{eq:rawbcomm}
 \left\langle p_-(v)\partial_z t(z)\right\rangle
 =\frac{i\eta\mathcal G}{d(1+v^2)^2q_z^{k+1}}\mathcal P_b(v,z),
\end{equation}
where the numerator, written by powers of $v$, is
\begin{align}
 \mathcal P_b(v,z)
 &=v^2\{(4k^2-1)z^2-2k-1\}+4kvz\notag\\
 &\quad+(4k^2-4k+1)z^2+1-2k.
\end{align}

The individual summands can also be tracked before collection. Denote the three terms of $p_-$ in~\eqref{eq:oddleftsub} by $p_-^{[1]},p_-^{[2]},p_-^{[3]}$ in their displayed order. Their reduced traces are
\begin{align}
 \langle p_-^{[1]}(v)\partial_z t(z)\rangle
 &=\frac{-2i\eta\mathcal G\{-2kvz+(2k-1)z^2-1\}}
 {d(1+v^2)^2q_z^{k+1}},\\
 \langle p_-^{[2]}(v)\partial_z t(z)\rangle
 &=\frac{i\eta\mathcal G\{(2k-1)z^2-1\}}
 {(1+v^2)q_z^{k+1}},\\
 \langle p_-^{[3]}(v)\partial_z t(z)\rangle&=0.
\end{align}
Putting the first two fractions over the denominator $d(1+v^2)^2q_z^{k+1}$ gives exactly $i\eta\mathcal G\mathcal P_b$. This displays the connection and inverse-symbol contributions separately and fixes the numerator without relying on its final integral.

Here terms from $-d\eta N_c/4$ in $w$, from $N_cC_{x_n}$, and from the scalar $zq^{-3}$ term have zero hatted trace. Only $-\eta\sum c_a\hc_a\hc_n/4$ and its appearances in $Cw+wC$ remain. Applying the trace rules produces $A_aB_b$; angular averaging replaces $\xi_a\xi_b$ by $\delta_{ab}/d$, and summing gives $\mathcal G$. This accounts for every geometric contraction in~\eqref{eq:rawbcomm}.

For a pole of order two the projection is evaluated by
\[
 \pi_v^+F(v)=\frac{[(v-i)^2F(v)]_{v=i}}{(v-i)^2}
 +\frac{[\partial_v((v-i)^2F(v))]_{v=i}}{v-i}.
\]
Apply this formula to~\eqref{eq:rawbcomm}, put $v=z$, and include the prefactor $-i$. The complete projected scalar integrand simplifies to
\begin{align}\label{eq:kernelbcomm}
 \mathcal K_b(z)
 &=-\frac{i\eta\mathcal G(2k-1)(z+i)(2kz^2+iz-1)}{2d\,q_z^{k+2}}\notag\\
 &=\frac{\eta\mathcal G(2k-1)}{2d}
 \frac{dz^2-1+2iz-2ikz^3}{q_z^{k+2}}.
\end{align}
The odd part integrates to zero. For the even part, write
\begin{align}
 J_0&=\int_{\mathbb R}q_z^{-k-2}\,\dd z
 =\frac{\sqrt\pi\,\Gamma(k+3/2)}{\Gamma(k+2)},\\
 J_2&=\int_{\mathbb R}z^2q_z^{-k-2}\,\dd z
 =\frac{\Gamma(3/2)\Gamma(k+1/2)}{\Gamma(k+2)}.
\end{align}
Since $\Gamma(k+3/2)=(k+1/2)\Gamma(k+1/2)$, $J_0=dJ_2$. Consequently
\begin{equation}
 \Phi^{(1)}_b=\frac{\eta\mathcal G(2k-1)}{2d}(dJ_2-J_0)=0.
\end{equation}
Equivalently, the even numerator is a total derivative:
\[
 \frac{dz^2-1}{q_z^{k+2}}
 =-\partial_z\big(zq_z^{-k-1}\big).
\]
This supplies a second check of the cancellation. The trace in~\eqref{eq:rawbcomm} is generally nonzero; only its projected normal integral vanishes.

\subsection{Odd powers: case (c), the right subprincipal symbol}
The last case is
\[
 \Phi^{(1)}_c=-i\int_{\mathbb R}\left\langle\pi_z^+p\,\partial_z t_-\right\rangle\,\dd z.
\]
Differentiation of the four summands in~\eqref{eq:oddrightsub}, before taking the trace, gives
\begin{align}
 \partial_z t_-&=ik\{N_cq^{1-k}\ell_3
 +2(1-k)zCq^{-k}\ell_3+Cq^{1-k}\partial_z\ell_3\}\notag\\
 &\quad+k(k-1)\eta\{(C+zN_c)q^{-k-2}
              -2(k+2)z^2Cq^{-k-3}\}\notag\\
 &\quad-2kzwq^{-k-1}+2k(k+1)\eta zN_cq^{-k-2}.
\end{align}
As in case (b), all terms with no hatted connection factor vanish. The surviving traces, reduced with~\eqref{eq:hatbivtrace}--\eqref{eq:cbivtrace}, give the particularly simple two-variable expression
\begin{equation}\label{eq:rawccomm}
 \left\langle p(v)\partial_z t_-(z)\right\rangle
 =-\frac{2i\eta k\mathcal G(vz+1)(dz^2-1)}{d(1+v^2)q_z^{k+2}}.
\end{equation}

To see the collection explicitly, write the four summands of $t_-$ in~\eqref{eq:oddrightsub} as $t_-^{[1]},\ldots,t_-^{[4]}$. Their differentiated traces are
\begin{align}
 \langle p(v)\partial_z t_-^{[1]}(z)\rangle
 &=\frac{-2i\eta k\mathcal G\{-(2k+2)vz+dz^2-1\}}
 {d(1+v^2)q_z^{k+2}},\\
 \langle p(v)\partial_z t_-^{[2]}(z)\rangle&=0,\\
 \langle p(v)\partial_z t_-^{[3]}(z)\rangle
 &=\frac{-2i\eta k\mathcal G vz}{(1+v^2)q_z^{k+1}},\\
 \langle p(v)\partial_z t_-^{[4]}(z)\rangle&=0.
\end{align}
The common numerator is simplified using
\[
 -(2k+2)vz+dz^2-1+dvz(1+z^2)
 =(vz+1)(dz^2-1).
\]
This proves~\eqref{eq:rawccomm} by summing the four source terms before projection.

The first factor has a simple pole at $v=i$. Its projected scalar part uses
\[
 \pi_v^+\frac{vz+1}{1+v^2}
 =\frac{iz+1}{2i(v-i)}.
\]
After setting $v=z$, one has $(iz+1)/(z-i)=i$. Including $-i$ from the case formula therefore yields
\begin{equation}\label{eq:kernelccomm}
 \mathcal K_c(z)=-\frac{\eta k\mathcal G}{d}
 \frac{dz^2-1}{q_z^{k+2}},\qquad
 \Phi^{(1)}_c=-\frac{\eta k\mathcal G}{d}(dJ_2-J_0)=0.
\end{equation}
This case vanishes separately, rather than cancelling a nonzero integrated case (b). The projected kernels provide an explicit normal-integral verification in arbitrary dimension.

\subsection{Even-power factorization: symbols and all five cases}
Now put $a=m-2=k-1$, so $Q=D^{4-n}=L^{-a}$, and $P=TL^{-1}$. For $n\geq6$, $a\geq1$ and both factors have negative order. Their leading and next symbols are
\begin{align}
 p&=Tq^{-1},&p_-&=T\ell_3,\\
 t&=q^{-a},&t_-&=a q^{-a+1}\ell_3-ia(a-1)\eta zq^{-a-2}.\label{eq:evensymbols}
\end{align}
The order equation again leaves one loss, so the five derivative patterns are precisely those of Section~\ref{sec:boundary}. To prove their vanishing explicitly, first expand the Clifford content of $\ell_3$. The identity
\[
 Cc(e_b)+c(e_b)C=-2\xi_b,\qquad
 C\hc(e_b)\hc(N)=\hc(e_b)\hc(N)C
\]
gives
\begin{equation}\label{eq:ellselection}
 \ell_3=-i q^{-2}\left\{\frac{d\eta z}{2}\Id
 +\frac\eta2\sum_{b<n}\xi_b\hc(e_b)\hc(N)
 +\frac\eta2N_c c(\xi')\right\}-2i\eta zq^{-3}\Id.
\end{equation}
Thus every summand is scalar, purely hatted, or purely unhatted. No summand has both hatted degree two and unhatted degree two. It follows that
\begin{equation}\label{eq:elltracezero}
 \tau_\Lambda(T)=\tau_\Lambda(T_{,n})
 =\tau_\Lambda(T\ell_3)=\tau_\Lambda(T\partial_z\ell_3)=0.
\end{equation}
These equalities are preserved by the scalar projection $\pi_z^+$.

\paragraph{Case (a.I).}
The complete integrand is $-\sum_{b<n}\tau_\Lambda[(\partial_{\xi_b}\pi_z^+p)(\partial_{x_b}\partial_z t)]$. Since $\partial_{x_b}q=0$ at the center, $\partial_{x_b}t=-a q^{-a-1}\partial_{x_b}q=0$. This case is zero before taking a trace.

\paragraph{Case (a.II).}
The symbol derivatives are
\[
 p_{x_n}=T_{,n}q^{-1}-\eta Tq^{-2},\qquad
 \partial_z^2t=-2a q^{-a-1}+4a(a+1)z^2q^{-a-2}.
\]
Their product in $-\tfrac12\tau_\Lambda[\pi_z^+(p_{x_n})\partial_z^2t]$ is a sum of scalar factors multiplying $T_{,n}$ or $T$. Both traces in~\eqref{eq:elltracezero} vanish. No derivatives of $U,V,W,Z$ are suppressed by an assumption.

\paragraph{Case (a.III).}
Here $\partial_z p=-2zTq^{-2}$ and $t_{x_n}=-a\eta q^{-a-1}$, so $\partial_z t_{x_n}=2a(a+1)\eta zq^{-a-2}$. The integrand $-\tfrac12\tau_\Lambda[(\partial_z\pi_z^+p)(\partial_z t_{x_n})]$ is a scalar multiple of $\tau_\Lambda(T)=0$.

\paragraph{Case (b): the left subprincipal symbol.}
One has $p_-=T\ell_3$ and $\partial_z t=-2azq^{-a-1}$. Thus
\[
 -i\tau_\Lambda[(\pi_z^+p_-)\partial_z t]
 =2iazq^{-a-1}\pi_z^+\tau_\Lambda(T\ell_3)=0.
\]
Equation~\eqref{eq:ellselection} proves this before either sphere integration or normal integration.

\paragraph{Case (c): the right subprincipal symbol.}
Differentiating~\eqref{eq:evensymbols} gives
\begin{align}
 \partial_z t_-&=a q^{1-a}\partial_z\ell_3
       +2a(1-a)zq^{-a}\ell_3\notag\\
 &\quad-ia(a-1)\eta\{q^{-a-2}-2(a+2)z^2q^{-a-3}\}.
\end{align}
Since $\pi_z^+p=T\pi_z^+(q^{-1})$, the trace in $-i\tau_\Lambda[(\pi_z^+p)\partial_z t_-]$ is a combination of precisely the zero traces in~\eqref{eq:elltracezero}. This last case also vanishes identically. Some references interchange the labels (b) and (c) for the even factorization; here they consistently denote loss in the left and right symbol, respectively.

\subsection{Geometric conclusion, scalar weights and exceptional dimensions}
Restoring the common boundary prefactor, the two geometric boundary densities are
\begin{align}
 \Phi^{(1)}&=\frac{2^n\vol(S^{n-2})}{(2\pi)^n}
       (0+0+0+0+0)\,\dd V_\partial=0,\\
 \Phi^{(2)}&=\frac{2^n\vol(S^{n-2})}{(2\pi)^n}
       (0+0+0+0+0)\,\dd V_\partial=0.
\end{align}
The remaining geometric invariant for the odd powers was $\eta\mathcal G$, equivalently $-2K_\partial\mathcal G/(n-1)$ with the mean-curvature convention~\eqref{eq:meanconversion}. Its coefficient is zero by $dJ_2-J_0=0$. There are no remaining normal-derivative or tangential-divergence terms. Together with the interior computation, this proves Theorem~\ref{thm:doubleboundary}.

For a scalar weight $f$, replace $T$ by $fT$. The subprincipal symbols are simply multiplied by $f$. The additional term in the normal derivative is $f_{,n}T$, whose trace vanishes in case (a.II). In case (a.I) the right tangential derivative still vanishes. The two nonzero projected kernels are multiplied by $f|_\partial$ and still integrate to zero. This proves the weighted assertion without an integration-by-parts argument on $M$.

In dimension four the even factorization has $a=0$ and $Q=\Id$. Its leading normal derivative and all lower symbols are zero, so the boundary composition contributes nothing directly. The interior formula remains valid. This checks the endpoint that would be missed if one insisted that both factors have strictly negative order.

The normalized recovery statement is not asserted in dimension two: $n-2=0$, and the interior operator is the multiplication operator $T$, whose residue is zero. In that dimension the displayed products involve positive differential powers on the right rather than the two-negative-order setup used in the theorem. They do not supply a nondegenerate extension of Corollary~\ref{cor:R}; surface curvature requires a different residue order. In particular no division by $n-2$ is made at this endpoint.

On a product collar, $\eta=0$ makes the odd kernels zero before integration; on a non-product warped collar, they may be nonzero and their integrals still vanish. If all four fields are tangential, $\mathcal G=0$ provides another direct check. Neither check replaces the general calculation above, which permits arbitrary normal components and arbitrary first jets.

\subsection{Comparison with the ordered insertions of Section~\ref{sec:boundary}}
Let $\mathcal F_{\ell,f}$ be the residue of the factorization~\eqref{eq:boundarydef} with $fA_\ell$ in the first factor. Define $\operatorname{Alt}_{XY}F=F(X,Y;Z,W)-F(Y,X;Z,W)$, without a factor $1/2$.
\begin{corollary}[Boundary version of Corollary~\ref{cor:R}]\label{cor:boundaryR}
For the warped collar and factorization of Theorem~\ref{thm:boundary}, every smooth scalar weight $f$ and $n=2m\ge4$ satisfy
\begin{align}\label{eq:boundaryrecovery}
 \frac{\operatorname{Alt}_{XY}\mathcal F_{2,f}}{m-1}
 &=C_m\int_M fR(X,Y,Z,W)\dd V_g,\\
 -\frac{\operatorname{Alt}_{XY}\mathcal F_{3,f}}{m-1}
 &=C_m\int_M fR(X,Y,Z,W)\dd V_g.\notag
\end{align}
In particular there is no boundary correction in these two antisymmetrized functionals.
\end{corollary}
\begin{proof}
The interior differences were computed in Section~\ref{sec:curvature}. Every term in $S_2$ and $S_3$ contains $G$ or its normal derivative and is symmetric under $X\leftrightarrow Y$. In the third insertion the intermediate skew polynomial $L$ cancels between cases (b) and (c), as explicitly shown above. Thus both unweighted boundary differences vanish.

A left scalar weight contributes $f$ times each boundary term and, in case (a.II) only, an additional $\partial_nf$ term. Its coefficient is $-b_kGH$ for the second insertion and $u_kGH-v_kGP_1$ for the third. Both coefficients are again symmetric in $X,Y$. Tangential derivatives of $f$ do not occur because the left composition formula differentiates the first factor only in the covariable, except for the stated normal derivative. Hence the weighted boundary differences also vanish.
\end{proof}
One may alternatively antisymmetrize in both pairs and divide by $2(m-1)$; this gives the same curvature integral. This is a boundary extension of the closed commutator recovery through the specified $F_2$ factorization. It does not identify all possible truncations of the closed commutator expression.

For local recovery, weights supported in the interior determine each component of $R$ there. Smoothness of the metric and fields up to the boundary determines its boundary values by continuity. This recovers the ambient Riemann tensor restricted to the boundary; it should not be confused with the intrinsic curvature tensor of $\partial M$, which is related to ambient curvature by the Gauss equation.

\section{General de Rham-type spectral data}\label{sec:data}
The second Clifford action is extra information. An ordinary spectral triple $(\A,\HH,D)$ gives the represented one-forms $a[D,b]$ but does not, by itself, prescribe a map $c\mapsto\hc$. The distinction matters in the examples below.

\begin{definition}[Differential de Rham spectral data]\label{def:DR}
A differential de Rham datum consists of a faithfully represented unital smooth $*$-algebra $\A$, a Hilbert space $\HH$ with finite degree decomposition $\HH=\bigoplus_{p=0}^N\HH^p$, and a densely defined degree-one operator $\dd$ on a common invariant smooth core, such that:
\begin{enumerate}[label=(\roman*)]
\item $\dd^2=0$ on the core; $D=\dd+\dd^*$ is essentially self-adjoint with compact resolvent, and $[D,a]$ is bounded for $a\in\A$;
\item $[\dd,a]$ and $[\dd^*,a]$ extend to bounded operators, and the representation of $\A$ preserves degree;
\item $D$ defines a regular finitely summable spectral triple with an integer spectral dimension $n>2$, together with a specified pseudodifferential algebra and a residue trace $\fres$;
\item a bimodule $\mathcal E^1$ of one-forms and two represented maps $c,\hc$ are specified, with
$c(\dd a)=[\dd+\dd^*,a]$ and $\hc(\dd a)=[\dd-\dd^*,a]$.
\end{enumerate}
The maps in (iv) must be well defined on the chosen quotient of universal one-forms. Equivalently, both representations must vanish on its defining relations. This condition is part of the data, not a consequence of (i).
\end{definition}

\paragraph{Meaning of the common invariant smooth core.}
A core is a dense linear subspace $\mathcal C\subset\HH$ on which the
unbounded operators are initially defined and whose graph closure
recovers the specified closed operators. ``Common'' means that the
same $\mathcal C$ is used for $\dd,\dd^*,D$ and the products under
discussion; ``invariant'' means, in particular,
$a\mathcal C,\dd\mathcal C,\dd^*\mathcal C\subset\mathcal C$.
This makes identities such as $\dd^2=0$ and the commutator formulas
well-defined before taking closures. The word core does not mean the
kernel of an operator. In the classical example one may use smooth
forms, and on the smooth torus one uses
$\A_\theta\otimes\Lambda^*\mathbb C^n$.

On a closed manifold take the ordinary de Rham complex, its degree operator, the usual one-form bimodule, and the two exterior Clifford actions. On a flat noncommutative torus the construction in Section~\ref{sec:torus} provides another example. A convenient general quotient for (iv) is the universal one-form module modulo the intersection of the two representation kernels, provided the resulting quotient has the required module properties. It need not coincide with $\Omega_D^1$.

Definition~\ref{def:DR} is a Hilbert-complex definition. A stronger differential-calculus version also asks that $\dd$ extend a graded derivation on an algebra of forms and annihilate its unit. We explicitly distinguish these requirements: a nilpotent connection differential need not annihilate a distinguished constant section.

\begin{definition}[Two-action spectral data]
A two-action datum retains a regular spectral triple, a residue trace, a specified one-form bimodule, and two bounded representations $c,\hc$, but does not assume the existence of a nilpotent degree-one $\dd$ with $D=\dd+\dd^*$. All residue expressions below are defined on this larger class when their indicated operators belong to the pseudodifferential algebra.
\end{definition}
In a noncommutative bimodule we do not impose the classical anticommutation relations for arbitrary algebra-valued coefficients. For example, with central frame $e_i$ and left coefficients,
\begin{equation}\label{eq:NC-CAR}
 \{a c_i,b\hc_j\}=[a,b]c_i\hc_j.
\end{equation}
The right side generally does not vanish. Imposing the classical relation without qualifications would exclude the natural noncommutative torus module.

\subsection{The curvature residue and the tensor question}
For any two-action datum of dimension $n>2$, define
\begin{equation}\label{eq:NC-R}
 \RR_D(u,v,w,z)=\frac1{n-2}\fres
 \left([\hc(u),\hc(v)][c(w),c(z)]|D|^{2-n}\right).
\end{equation}
It is complex-multilinear and skew in each displayed pair. These properties follow directly from the definition. It is a \emph{curvature residue functional}; it is called a Riemann tensor functional only if the additional identities under discussion have been proved for the particular datum.

For a central-frame bimodule, a necessary balancing identity is
\begin{equation}\label{eq:balance}
 \RR_D(ua,v,w,z)=\RR_D(u,av,w,z),\qquad a\in\A,
\end{equation}
with corresponding identities at the other adjacent tensor positions. Pair symmetry and the first Bianchi identity would require
\begin{align}
 \RR_D(u,v,w,z)&=\RR_D(w,z,u,v),\label{eq:pairsym}\\
 \RR_D(u,v,w,z)+\RR_D(v,w,u,z)+\RR_D(w,u,v,z)&=0.\label{eq:Bianchi}
\end{align}
In a braided calculus different balancing or symmetry notions may be appropriate, but then the braiding must be specified. We use~\eqref{eq:balance}--\eqref{eq:Bianchi} as explicit classical-type tests. A trace alone does not prove them: it moves entire products cyclically, not individual noncommuting coefficients through $|D|^{2-n}$.

\subsubsection{Two ordered curvature functionals from Section 5.1}
Assume now that the spectral dimension is even, $n=2m>2$, and put $k=m-1$. In addition to~\eqref{eq:NC-R}, define
\begin{align}\label{eq:newR2}
 \RR_D^{(2)}(u,v,w,z)
 &=\frac1k\fres\big([\hc(u),\hc(v)]D c(w)c(z)D^{1-n}\big),\\
 \RR_D^{(3)}(u,v,w,z)
 &=-\frac1k\fres\big([\hc(u),\hc(v)]c(w)D c(z)D^{1-n}\big).
 \label{eq:newR3}
\end{align}
These definitions require two-action data, including the specified hatted representation and a pseudodifferential calculus containing the indicated products. An ordinary triple without this additional representation does not determine them. Signed inverse powers mean $D^{1-n}=D(D^2)^{-m}$ modulo the kernel projection. In particular one must not replace this odd inverse power by $|D|^{1-n}$.

Both expressions are complex-multilinear and skew in their first pair. Unlike~\eqref{eq:NC-R}, they need not be skew in their last pair. Balancing, pair symmetry and Bianchi remain separate tests. The normalizing minus sign in~\eqref{eq:newR3} is essential: Section 5.1 gives, on classical Hodge data,
\begin{equation}\label{eq:threeclassical}
 \RR_D^{(2)}(X,Y,Z,W)=\RR_D^{(3)}(X,Y,Z,W)
 =\RR_D(X,Y,Z,W)=C_m\int_MR(X,Y,Z,W)\dd V_g.
\end{equation}
The equality also holds with a scalar weight on the left. Corollary~\ref{cor:boundaryR} supplies the corresponding boundary extensions for the chosen truncation products. This classical equality is a theorem, not an identification imposed on general spectral data.

To see exactly what is being changed, abbreviate $H=[\hc(u),\hc(v)]$, $B=c(w)c(z)$ and $Q=|D|^{2-n}$. The following algebraic decompositions give
\begin{equation}\label{eq:ordereddefect}
 \RR_D^{(2)}-\RR_D
 =\frac1{2k}\fres\big(H\{c(w),c(z)\}Q\big)
   +\frac1k\fres\big(H[D,B]D^{1-n}\big).
\end{equation}
Indeed $DBD^{1-n}=BQ+[D,B]D^{1-n}$ and
$B=([c(w),c(z)]+\{c(w),c(z)\})/2$. Neither term on the right is automatically zero on noncommutative data. This identity also explains why imposing last-pair skewness afterwards would change the functional. We retain the two definitions as written and calculate their values rather than presume tensoriality.

\section{Exact noncommutative torus calculations}\label{sec:torus}
Let $\A_\theta=C^\infty(\mathbb T^n_\theta)$ be generated by unitaries with
$U_iU_j=e^{2\pi i\theta_{ij}}U_jU_i$. Use the real derivations
$\partial_j(U^\alpha)=i\alpha_jU^\alpha$ and the canonical trace $\tau(U^\alpha)=\delta_{\alpha,0}$. On
\begin{equation}
 \HH=L^2(\A_\theta,\tau)\otimes\Lambda^*\mathbb C^n,
 \qquad \mathcal E=\A_\theta\otimes\Lambda^*\mathbb C^n,
\end{equation}
put
\begin{equation}
 \dd_0=\varepsilon_j\partial_j,\quad
 \dd_0^*=-\iota_j\partial_j,\quad
 D_0=c_j\partial_j,\quad D_0^2=-\sum_j\partial_j^2.
\end{equation}
The finite-dimensional vector space $V=\mathbb C^n$ is the exterior factor requested in this construction. The module of one-forms has central frame $e^j$, with
$c(a_je^j)=a_jc_j$ and $\hc(a_je^j)=a_j\hc_j$. The operators are left multipliers in their algebra coefficients. The nilpotence of $\dd_0$ follows from commuting derivations and anticommuting exterior multiplications.

The toroidal calculus and residue are taken in the standard sense of \cite{HaLeePonge,Ponge}. The toroidal residue convention is
\begin{equation}\label{eq:torres}
 \fres_\theta(P)=\int_{S^{n-1}}(\tau\otimes\tr_\Lambda)
 (\sigma_{-n}(P)(\xi))\dd S(\xi).
\end{equation}
For a torus of period $2\pi$ this agrees with~\eqref{eq:Wres}, since $\tau$ is the spatial average. Set $\mathsf V_n=r\nu_{n-1}$.

\begin{proposition}[The flat torus]
For $D_0$ and arbitrary algebra-valued one-forms $u,v,w,z$,
\begin{equation}
 \RR_{D_0}(u,v,w,z)=0.
\end{equation}
\end{proposition}
\begin{proof}
The symbol of $|D_0|^{2-n}$ is exactly $|\xi|^{2-n}$ modulo the irrelevant zero-mode smoothing term. Left multiplication by the four-form insertion has no $\xi$ derivatives in the composition formula. There is no homogeneous term of degree $-n$. The conclusion holds even when the coefficients of the four one-forms do not commute.
\end{proof}

\paragraph{The differential and flatness calculation in detail.}
The relations $\varepsilon_i\varepsilon_j=-\varepsilon_j\varepsilon_i$
and $[\partial_i,\partial_j]=0$ give
$\dd_0^2=\frac12\varepsilon_i\varepsilon_j[\partial_i,\partial_j]=0$.
Integration by parts for the canonical trace gives
$\partial_j^*=-\partial_j$, hence
$\dd_0^*=-\iota_j\partial_j$ and
$D_0^2=-\sum_j\partial_j^2$. The one-form insertion is an order-zero
left multiplier. Its star product with the constant symbol
$|\xi|^{2-n}$ has no derivative terms because its covariable
derivatives vanish. Thus absence of a degree $-n$ term proves flatness
even for nonconstant, noncommuting input coefficients.

\subsection{An ordered formula for zero-order perturbations}
Let $\Psi=\Psi^*$ be a smooth $\A_\theta$-valued exterior endomorphism and set $D_\Psi=D_0+\Psi$. It is a self-adjoint regular elliptic spectral triple with the same principal symbol; it is even when $\Psi$ is odd. It need not be a differential de Rham datum. Write
\begin{align}
 D_\Psi^2&=-\partial_j\partial_j+b_j\partial_j+V_0,\\
 b_j&=c_j\Psi+\Psi c_j,\qquad V_0=c_j\partial_j\Psi+\Psi^2,\\
 \mathcal E_\Psi&=V_0-\tfrac12\partial_jb_j+\tfrac14\sum_jb_jb_j\notag\\
 &=\Psi^2+\tfrac12[c_j,\partial_j\Psi]
       +\tfrac14\sum_j\{c_j,\Psi\}^2.\label{eq:Epsi}
\end{align}
The order in these products is part of the formula.

\paragraph{Squaring the perturbation and completing the connection.}
On a smooth vector $s$,
\[
 (D_0+\Psi)^2s=-\partial_j\partial_js
 +(c_j\Psi+\Psi c_j)\partial_js+(c_j\partial_j\Psi+\Psi^2)s.
\]
Set $\nabla_j^\Psi=\partial_j-b_j/2$. Expanding without reordering
coefficients gives
\[
 -\sum_j(\nabla_j^\Psi)^2
 =-\partial_j\partial_j+b_j\partial_j
   +\tfrac12\partial_jb_j-\tfrac14b_jb_j.
\]
Subtracting this expression from $D_\Psi^2$ proves~\eqref{eq:Epsi}.
In particular the derivative term in the potential is the commutator
$[c_j,\partial_j\Psi]/2$, not an anticommutator.

\begin{lemma}[Flat-principal-symbol recursion]\label{lem:flatrec}
For $q=|\xi|^2$, $b_\xi=b_j\xi_j$, and $a\geq1$,
\begin{align}
 \sigma_{-2a}(L^{-a})&=q^{-a},\\
 \sigma_{-2a-1}(L^{-a})&=-ai b_\xi q^{-a-1},\\
 \sigma_{-2a-2}(L^{-a})&=-aV_0q^{-a-1}
 -\tfrac{a(a+1)}2 b_\xi^2q^{-a-2}
 +a(a+1)(\partial_i b_j)\xi_i\xi_jq^{-a-2}.
\end{align}
\end{lemma}
\begin{proof}
For $a=1$ this is~\eqref{eq:inverse-rec} with constant $l_2=q$. Multiplying the expansion by the symbol of $L^{-1}$ gives the induction. At the second loss the only products are those with losses $(0,2),(1,1),(2,0)$ and the differentiated leading--subleading product. No coefficient is reordered. In particular $b_\xi^2$ remains the displayed ordered square.
\end{proof}

\paragraph{The coefficients in the flat recursion.}
For $a=1$ the three symbols are
\[
 q^{-1},\quad -ib_\xi q^{-2},\quad
 -V_0q^{-2}-b_\xi^2q^{-3}
       +2(\partial_i b_j)\xi_i\xi_jq^{-3}.
\]
The last summand comes from
$-q^{-1}(2\xi_i)D_{x_i}(-ib_\xi q^{-2})$.
In the induction from $a$ to $a+1$, the coefficients of $b_\xi^2$
add as $-a(a+1)/2-1-a=-(a+1)(a+2)/2$.
The derivative coefficients add as $a(a+1)+2+2a=(a+1)(a+2)$;
the last $2a$ is the differentiated principal--subprincipal product.
These counts retain the order $b_ib_j$ and require no commutativity.
At $a=m-1$, angular averaging replaces $\xi_i\xi_j$ by
$\delta_{ij}/n$, yielding
\[
 -aV_0-\frac a4\sum_jb_jb_j+\frac a2\partial_jb_j
 =-a\mathcal E_\Psi.
\]

\begin{theorem}[Perturbed torus curvature residue]\label{thm:torus}
For $n=2m>2$, define $F(u,v,w,z)=[\hc(u),\hc(v)][c(w),c(z)]$. Then
\begin{equation}\label{eq:torusR}
 \boxed{\RR_{D_\Psi}(u,v,w,z)
 =-\frac{\mathsf V_n}{2}\,
 \tau\!\left(\tau_\Lambda\big(F(u,v,w,z)\mathcal E_\Psi\big)\right).}
\end{equation}
This formula holds for arbitrary smooth noncommuting coefficients of the one-forms and of $\Psi$.
\end{theorem}
\begin{proof}
Set $a=m-1$ in Lemma~\ref{lem:flatrec} and integrate its last line over the unit sphere. Since $n=2(a+1)$ the result before taking any trace is
$-a\nu_{n-1}(V_0-\partial_jb_j/2+b_jb_j/4)$.
Left multiplication by $F$ introduces no differentiated symbol terms. Divide by $n-2=2a$ and use~\eqref{eq:torres}. The derivative commutator in~\eqref{eq:Epsi} cannot be discarded inside a trace weighted by the Clifford-valued $F$.
\end{proof}

\subsection{A nilpotent example with failure of balancing}
Let $\mu\in\mathbb R\setminus\{0\}$ and
\begin{equation}\label{eq:mu}
 \dd_\mu=\dd_0+\mu\varepsilon_1,\qquad
 D_\mu=\dd_\mu+\dd_\mu^*=D_0+\mu\hc_1.
\end{equation}
Because the additional exterior multiplication is constant,
$\dd_\mu^2=0$. The commutators with $\A_\theta$ are unchanged, the perturbation is odd and self-adjoint, and
\begin{equation}
 D_\mu^2=D_0^2+\mu^2,\qquad \mathcal E_{\mu\hc_1}=\mu^2\Id.
\end{equation}
This is a differential de Rham datum in the Hilbert-complex sense. It is not a unital differential-calculus extension: $\dd_\mu1=\mu e^1\neq0$.

For any fixed frame indices $i,j$ and $a,b,c,d\in\A_\theta$,
\begin{equation}\label{eq:commexample}
 \RR_{D_\mu}(ae^i,be^i,ce^j,de^j)
 =\frac{\mathsf V_n\mu^2}{2}\tau([a,b][c,d]).
\end{equation}
Indeed $[a\hc_i,b\hc_i]=[a,b]$ while $[cc_j,dc_j]=-[c,d]$.
Choose $U,V$ with $UV=e^{i\alpha}VU$ and $e^{i\alpha}\neq1$. Then
\begin{align}
 \RR_{D_\mu}(Ue^i,Ve^i,V^*e^j,U^*e^j)
 &=\mathsf V_n\mu^2(1-\cos\alpha)>0,\label{eq:positive}\\
 \RR_{D_\mu}(e^i,UV e^i,V^*e^j,U^*e^j)&=0.\label{eq:zero}
\end{align}
The first equality uses $[V^*,U^*]=[U,V]^*$ and
$\tau([U,V][U,V]^*)=2-2\cos\alpha$.
Since the frame is central,~\eqref{eq:positive} and~\eqref{eq:zero} contradict exactly the balancing identity~\eqref{eq:balance}. Hence even the existence of a nilpotent Hilbert differential does not make~\eqref{eq:NC-R} a tensor over $\A_\theta$.

\paragraph{Checking the nonzero and balanced values.}
The relevant anticommutation and square identities are
\[
 \{\dd_0,\varepsilon_1\}=0,\quad\varepsilon_1^2=0,
 \qquad\{D_0,\hc_1\}=0,\quad\hc_1^2=1.
\]
They give $\dd_\mu^2=0$ and
$D_\mu^2=D_0^2+\mu^2$. The relation $UV=e^{i\alpha}VU$ gives
\[
 [U,V]=(1-e^{-i\alpha})UV.
\]
Consequently
\[
 [U,V][U,V]^*=|1-e^{-i\alpha}|^2I
                  =2(1-\cos\alpha)I.
\]
The comparison inputs $(e^i,UV e^i)$ give $[1,UV]=0$ in the first
commutator. Moving the coefficient $U$ from the first slot to the
second slot is exactly the balancing operation, so the displayed
positive and zero values contradict that identity. Nilpotence of a
Hilbert differential imposes no such balancing identity on the residue.

For numerical frame inputs the same curvature residue is zero. This is not a contradiction: algebra-valued one-forms test information invisible to the numerical frame, and the failure of balancing prevents that information from being interpreted as ordinary frame components of a tensor. The usual three-$c$ torsion residue of this example also vanishes, since its possible term contains three $c$'s and one $\hc$ and has zero exterior trace. Thus torsion-freeness does not repair balancing.

\subsection{Pair symmetry and Bianchi can fail as well}
In every even $n\geq4$, take the constant odd self-adjoint perturbation
\begin{equation}\label{eq:badpsi}
 \Psi=c_1\hc_3\hc_1+c_2\hc_1\hc_4.
\end{equation}
Both $c_i$ and the indicated hatted bivectors are skew-adjoint and commute with one another, so their products are self-adjoint. Let $H_1=\hc_3\hc_1$, $H_2=\hc_1\hc_4$. Then
\begin{equation}
 \mathcal E_\Psi=c_1c_2[H_1,H_2]
 =2c_1c_2\hc_3\hc_4.
\end{equation}
The first equality follows by cancellation of $-\sum H_i^2$ in $\Psi^2$ against $\sum H_i^2$ from $\sum\{c_i,\Psi\}^2/4$. Formula~\eqref{eq:torusR} gives
\begin{align}
 \RR_{D_\Psi}(e^3,e^4,e^1,e^2)&=-4\mathsf V_n,\\
 \RR_{D_\Psi}(e^1,e^2,e^3,e^4)&=0,\\
 \RR_{D_\Psi}(e^4,e^1,e^3,e^2)
 =\RR_{D_\Psi}(e^1,e^3,e^4,e^2)&=0.
\end{align}
Therefore pair symmetry and the first Bianchi identity both fail. For example the Bianchi sum with $(u,v,w,z)=(e^3,e^4,e^1,e^2)$ equals $-4\mathsf V_n$. These are exact Clifford traces, not numerical approximations. This example belongs to the two-action class; no nilpotent $\dd$ is asserted for~\eqref{eq:badpsi}.

\paragraph{The Clifford calculation in the symmetry counterexample.}
Here $H_1^2=H_2^2=-1$, $H_1H_2=\hc_3\hc_4$ and
$H_2H_1=-\hc_3\hc_4$. Since the hatted bivectors commute with all
$c_i$, one has $b_j=-2H_j$ and
\[
 \Psi^2=-H_1^2-H_2^2+c_1c_2[H_1,H_2],\qquad
 \tfrac14\sum_jb_j^2=H_1^2+H_2^2.
\]
This proves the displayed potential. For the first input tuple the
insertion is $4\hc_3\hc_4c_1c_2$, and its product with the potential
is $8I$, since $(c_1c_2)^2=(\hc_3\hc_4)^2=-1$.
The prefactor $-\mathsf V_n/2$ therefore gives $-4\mathsf V_n$.
For each of the other three tuples a nonempty hatted Clifford monomial
remains after multiplication; its trace is zero. This proves both the
pair-interchange failure and the nonzero Bianchi sum explicitly.

\begin{remark}[What can still be called curvature]
The functional~\eqref{eq:NC-R} always has its explicitly proved two pair-skew symmetries. The three preceding computations show that neither flat-principal-symbol perturbations nor spectral torsion tests establish the remaining tensor axioms. One can project a multilinear functional onto a chosen algebraic symmetry class, but such a projection does not restore module balancing and changes the definition. A connection-based or braided curvature theory requires additional structures and a comparison theorem.
\end{remark}

\subsection{The two ordered functionals on perturbed tori}
Here and below $H=[\hc(u),\hc(v)]$. We first give a formula for arbitrary algebra-valued inputs, retaining their derivatives. This is needed because commuting a one-form through $D$ is not justified in general.

For a smooth zero-order endomorphism $B$, define the ordered expression
\begin{align}\label{eq:Kpsi}
 \mathcal K_\Psi(B)={}&c_i\partial_i(B\Psi)+\Psi B\Psi\notag\\
 &+\tfrac12\big(c_i(\partial_iB)c_j+c_jB\Psi+\Psi Bc_j\big)b_j
 +\tfrac12c_iBc_iV_0\notag\\
 &+\tfrac18c_iBc_j\left(\delta_{ij}\sum_l b_lb_l+b_ib_j+b_jb_i\right)\notag\\
 &+\tfrac14c_iBc_j\left(\partial_ib_j+\partial_jb_i
                         -\delta_{ij}\partial_lb_l\right).
\end{align}
All products are in the displayed order; $b_j,V_0$ are defined above. Then
\begin{align}\label{eq:orderedtorus}
 \RR_{D_\Psi}^{(2)}&=\frac{\mathsf V_n}{k}
 \tau\tau_\Lambda\big(H\mathcal K_\Psi(c(w)c(z))\big),\\
 \RR_{D_\Psi}^{(3)}&=-\frac{\mathsf V_n}{k}
 \tau\tau_\Lambda\big(Hc(w)\mathcal K_\Psi(c(z))\big).\notag
\end{align}
This formula involves no unspecified symbol coefficient or normal integral.

\paragraph{Derivation.}
Write $ADBD^{1-n}=(ADBD)L^{-m}$. The degree-two, degree-one and degree-zero insertion symbols are
\begin{align}
 p_2&=-Ac_iBc_j\xi_i\xi_j,\\
 p_1&=iA\big(c_i(\partial_iB)c_j+c_jB\Psi+\Psi Bc_j\big)\xi_j,\\
 p_0&=A\{c_i\partial_i(B\Psi)+\Psi B\Psi\}.
\end{align}
Since the principal symbol of $L$ is constant in $x$, the degree $-n$ part of the product is precisely
\[
 p_2\sigma_{-n-2}(L^{-m})+p_1\sigma_{-n-1}(L^{-m})
 +p_0|\xi|^{-n}+\partial_{\xi_l}p_2D_{x_l}\sigma_{-n-1}(L^{-m}).
\]
Use Lemma~\ref{lem:flatrec}. The second angular moment contributes $m/n=1/2$; the fourth contributes $m(m+1)/(n(n+2))=1/4$. The $b_ib_j$ term consequently has coefficient $1/8$. The differentiated product changes the symmetric derivative coefficient from $-1/4$ to $+1/4$, leaving the trace term $-\delta_{ij}\partial_lb_l/4$. These are exactly the five lines of~\eqref{eq:Kpsi}. Left multiplication by $A$ introduces no derivatives of $A$.

\paragraph{The ordered insertion and the differentiated product.}
The expansion on a test vector is
\begin{align*}
 D_\Psi B D_\Psi s
 ={}&c_iBc_j\partial_i\partial_js
 +\bigl(c_i(\partial_iB)c_j+c_jB\Psi+\Psi Bc_j\bigr)\partial_js\\
 &+\{c_i\partial_i(B\Psi)+\Psi B\Psi\}s.
\end{align*}
For instance, $c_i\partial_i(B\Psi s)$ supplies both
$c_i\partial_i(B\Psi)s$ and $c_iB\Psi\partial_i s$; the latter index
is renamed $j$ in the displayed first-order coefficient. Multiplying
by $A$ on the left gives the three symbols above, using
$\sigma(\partial_i\partial_j)=-\xi_i\xi_j$.

The $p_2S_{m,2}$ product contributes the derivative expression
\[
 -\tfrac14Ac_iBc_j
 (\delta_{ij}\partial_lb_l+\partial_ib_j+\partial_jb_i).
\]
On the other hand,
$D_{x_l}S_{m,1}=-m(\partial_lb_j)\xi_jq^{-m-1}$ and
$\partial_{\xi_l}p_2=-A(c_lBc_j\xi_j+c_iBc_l\xi_i)$.
Their spherical average is
$Ac_iBc_j(\partial_ib_j+\partial_jb_i)/2$.
Adding these two terms gives exactly the last line of~\eqref{eq:Kpsi}.
The other products give respectively $p_0$, the coefficient $1/2$
times the first-order coefficient followed by $b_j$, the potential
term $c_iBc_iV_0/2$, and the fourth-moment coefficient $1/8$ multiplying
$\delta_{ij}\sum_l b_lb_l+b_ib_j+b_jb_i$. Thus all terms of the
ordered formula arise from explicit symbol products.

\paragraph{Flat torus.}
When $\Psi=0$, all terms in $\mathcal K_\Psi(B)$ vanish, even for nonconstant algebra-valued $B$. Thus
\begin{equation}
 \RR_{D_0}^{(2)}(u,v,w,z)=\RR_{D_0}^{(3)}(u,v,w,z)=0.
\end{equation}
This proves flatness of all three residue functionals on the full represented one-form module, not merely on the numerical frame.

\paragraph{The nilpotent perturbation $\Psi=\mu\hc_1$.}
In this example $b_j=0$ and $V_0=\mu^2$. Therefore
\[
 \mathcal K_\Psi(B)=\mu c_i(\partial_iB)\hc_1
                  +\mu^2\hc_1B\hc_1+\frac{\mu^2}{2}\sum_ic_iBc_i.
\]
The derivative term has an odd number of hatted generators after either insertion and hence zero exterior trace. Write $u=u_ie^i$, $v=v_ie^i$, $w=w_je^j$, $z=z_je^j$ and put
\[
 h_0=\sum_i[u_i,v_i],\qquad B_0=\sum_jw_jz_j.
\]
Only the scalar hatted part of $H$ contributes. Decomposing $c(w)c(z)$ into scalar and bivector parts and using
$\sum_ic_i c(z)c_i=(n-2)c(z)$ gives
\begin{align}\label{eq:muordered}
 \RR_{D_\mu}^{(2)}(u,v,w,z)&=\mathsf V_n\mu^2\tau(h_0B_0),\\
 \RR_{D_\mu}^{(3)}(u,v,w,z)&=\frac{k-1}{k}\mathsf V_n\mu^2\tau(h_0B_0).\notag
\end{align}
These formulas hold for arbitrary smooth algebra coefficients. By comparison,
\[
 \RR_{D_\mu}(u,v,w,z)=\frac{\mathsf V_n\mu^2}{2}
 \tau\left(h_0\sum_j[w_j,z_j]\right).
\]
For $UV=e^{i\alpha}VU$ and the inputs $(Ue^i,Ve^i,V^*e^j,U^*e^j)$, one obtains
\begin{align}
 \RR^{(2)}_{D_\mu}&=\mathsf V_n\mu^2(1-e^{-i\alpha}),\\
 \RR^{(3)}_{D_\mu}&=\frac{k-1}{k}\mathsf V_n\mu^2(1-e^{-i\alpha}),\\
 \RR_{D_\mu}&=\mathsf V_n\mu^2(1-\cos\alpha).
\end{align}
The balanced comparison with first two inputs $(e^i,UV e^i)$ is zero. Hence the second ordered functional fails balancing; the third does so for $k>1$. At $k=1$ the third functional vanishes identically in this particular example, and no balancing obstruction for it follows from this test. All three vanish on numerical frame inputs.

\paragraph{The constant perturbation of~\eqref{eq:badpsi}.}
Set $H_1=\hc_3\hc_1$, $H_2=\hc_1\hc_4$ and $H_j=0$ otherwise. Then $\Psi=c_jH_j$, $b_j=-2H_j$, and all coefficient derivatives vanish. Substitution into~\eqref{eq:Kpsi} yields
\begin{align}
 \tau_\Lambda([\hc_3,\hc_4]\mathcal K_\Psi(c_1c_2))&=-4k,\\
 \tau_\Lambda([\hc_3,\hc_4]c_1\mathcal K_\Psi(c_2))&=4k.
\end{align}
For each $a=2,3$, the resulting values are
\begin{align}
 \RR^{(a)}_{D_\Psi}(e^3,e^4,e^1,e^2)&=-4\mathsf V_n,\\
 \RR^{(a)}_{D_\Psi}(e^1,e^2,e^3,e^4)&=0,\\
 \RR^{(a)}_{D_\Psi}(e^4,e^1,e^3,e^2)
 =\RR^{(a)}_{D_\Psi}(e^1,e^3,e^4,e^2)&=0.
\end{align}
Thus both new functionals fail pair interchange and Bianchi on these data. In this test their values agree with the original functional, while~\eqref{eq:muordered} proves that agreement is not universal.

\section{Some other examples}\label{sec:examples}
We distinguish an example as originally specified from an exterior lift supplied with additional data. The original spinorial examples in \cite{DSZtorsion} do not automatically include the second action required by~\eqref{eq:NC-R}. Assigning a value to an undefined expression would obscure the main obstruction. The computations below state the added structure each time.

\subsection{Einstein--Yang--Mills: an explicitly chosen exterior lift}
The original example uses spinors and the matrix algebra acting by left multiplication on $M_N(\mathbb C)$. Replace the spinor bundle by the exterior bundle and choose the induced connection on
$\Lambda^*T^*M\otimes\operatorname{End}(\mathbb C^N)$. Let $\nabla^A$ be the adjoint gauge connection, with curvature $F^A_{ij}$. Its Dirac-type operator is
\begin{equation}
 D_A=c_i(\nabla_i^\Lambda\otimes1+1\otimes\nabla_i^A),
 \qquad E_A=E_\Lambda\otimes1+\tfrac12c_ic_j\otimes F^A_{ij}.
\end{equation}
For numerical geometric one-forms, take $c,\hc$ tensored with the identity on the gauge fiber. The gauge term has zero trace against $[\hc(X),\hc(Y)][c(Z),c(W)]$: the hatted bivector has zero normalized trace and the gauge term contains no hats. Hence
\begin{equation}\label{eq:EYM}
 \RR_{D_A}(X,Y,Z,W)=N^2 C_m\int_MR(X,Y,Z,W)\dd V_g.
\end{equation}
The gauge curvature is invisible to this restricted test. This does not assert that the gauge connection is flat. With matrix-valued inserted forms the static formula~\eqref{eq:static}, keeping the ordered matrix products and $E_A$, gives the exact answer and can have additional terms. Moreover the covariant exterior differential satisfies $(\dd_A)^2=F^A\wedge$; when the induced adjoint curvature is nonzero, this lift is two-action data but not nilpotent differential de Rham data.

\paragraph{The adjoint connection and the meaning of its curvature.}
Locally, on $T\in\operatorname{End}(\mathbb C^N)$,
$\nabla_i^AT=\partial_iT+[A_i,T]$. Define the fundamental curvature
$\mathcal F_{ij}=\partial_iA_j-\partial_jA_i+[A_i,A_j]$.
The product rule gives
\begin{align*}
 \nabla_i^A\nabla_j^AT={}&\partial_i\partial_jT
 +[\partial_iA_j,T]+[A_j,\partial_iT]+[A_i,\partial_jT]
 +[A_i,[A_j,T]].
\end{align*}
Subtract the same expression with $i,j$ interchanged. The derivatives
of $T$ cancel, and the Jacobi identity
$[A_i,[A_j,T]]-[A_j,[A_i,T]]=[[A_i,A_j],T]$ gives
\[
 [\nabla_i^A,\nabla_j^A]T=[\mathcal F_{ij},T],\qquad
 F^A_{ij}=\operatorname{ad}(\mathcal F_{ij}).
\]
In a noncoordinate frame the left side is
$([\nabla_X^A,\nabla_Y^A]-\nabla_{[X,Y]}^A)T$.
The tensor product connection has curvature
$\Omega^\Lambda_{ij}\otimes I+I\otimes F^A_{ij}$, since the two
factor connections commute. Squaring its Dirac operator at a geodesic
center gives $D_A^2=\widetilde\nabla^*\widetilde\nabla+E_A$, with the
potential displayed above.

For $T_\Lambda=[\hc(X),\hc(Y)][c(Z),c(W)]$, the trace identities give
\[
 \tau_\Lambda(T_\Lambda)=0,\qquad
 \tau_\Lambda(T_\Lambda c_ic_j)=0,\qquad
 \tau_\Lambda(T_\Lambda E_\Lambda)=-2R(X,Y,Z,W).
\]
The middle equality follows because the hatted commutator has zero
scalar component and $[c(Z),c(W)]c_ic_j$ has no hatted generators.
Thus the gauge term disappears before the internal trace is taken.
Since $\tr_{\operatorname{End}(\mathbb C^N)}I=N^2$, formula~\eqref{eq:static}
divided by $n-2$ proves~\eqref{eq:EYM}. The factor would be $N$ for the
fundamental rank-$N$ representation, so the representation matters.
Moreover $\dd_A^2T=F^A\wedge T$ involves the adjoint curvature. A
nonzero central fundamental curvature can have zero adjoint action;
the obstruction to nilpotence is specifically $F^A\ne0$.
For matrix-valued one-forms, a hatted commutator can acquire a scalar
Clifford component proportional to a coefficient commutator. In that
case the preceding zero-trace argument is not available, and the
ordered static formula must be used.

\subsection{The two-sheet model and a finite-direction computation}

\paragraph{Definition of the two-sheet model.}
Let $M$ be the closed even-dimensional manifold used above. The algebra
and Hilbert space of the chosen exterior model are
\[
 \mathcal A_2=C^\infty(M)\otimes(\mathbb C\oplus\mathbb C),\qquad
 \mathcal H_2=L^2(M,\Lambda^*T^*M)\otimes\mathbb C^2.
\]
An element $(a_1,a_2)$ acts by $\operatorname{diag}(a_1,a_2)$.
Thus the commutative algebra describes two copies of $M$, while the
off-diagonal finite Dirac operator couples the copies. A geometric
one-form is lifted identically on both copies; a finite one-form is
represented by an off-diagonal matrix and describes the discrete
direction between them. The latter is not an additional tangent vector
of the base manifold. This definition specifies the exterior lift used
here, rather than identifying it with the original spinorial model.
For a concrete exterior lift choose $\mu>0$ and
\begin{equation}
 D_\mu^{(2)}=D_\Lambda\otimes1+\gamma\otimes D_F,
 \quad \gamma=(-1)^{\deg},\quad
 D_F=\mu(E_{12}+E_{21}).
\end{equation}
Then $(D_\mu^{(2)})^2=L\otimes1+\mu^2$. For geometric inputs lifted diagonally, the added scalar potential has zero trace against the mixed commutator insertion, so
\begin{equation}
 \RR_{D_\mu^{(2)}}(X,Y,Z,W)=2C_m\int_MR(X,Y,Z,W)\dd V_g.
\end{equation}
This is a statement about the geometric submodule only.

\paragraph{The square and the geometric residue.}
The degree parity satisfies $\gamma^2=I$ and
$D_\Lambda\gamma+\gamma D_\Lambda=0$. Also
$D_F^2=\mu^2(E_{11}+E_{22})$. Expanding the square therefore gives
\[
 (D_\mu^{(2)})^2
 =L\otimes I+(D_\Lambda\gamma+\gamma D_\Lambda)\otimes D_F
       +I\otimes D_F^2=L\otimes I+\mu^2I.
\]
For diagonal geometric inputs the insertion is $T_\Lambda\otimes I_2$.
Its trace against $\mu^2I$ is zero because
$\tr_\Lambda T_\Lambda=0$. The remaining trace is twice the ordinary
exterior trace, proving the geometric formula above.

To specify the missing finite hatted action, choose the two-term complex with $\dd_F=\mu E_{21}$. Set $Q_F=\dd_F-\dd_F^*$ and represent finite one-forms by
\begin{equation}
 c_F(a\dd b)=a[D_F,b],\qquad
 \hc_F(a\dd b)=a[Q_F,b],\qquad a,b\in\mathbb C\oplus\mathbb C.
\end{equation}
The total differential $\dd_\Lambda\otimes1+\gamma\otimes\dd_F$ is nilpotent. For a finite one-form whose $c_F$ matrix is
$C=\begin{psmallmatrix}0&p\\q&0\end{psmallmatrix}$, its hatted matrix is
$H=\begin{psmallmatrix}0&-p\\q&0\end{psmallmatrix}$.
These maps descend on the finite two-dimensional off-diagonal module.

\paragraph{The finite differential and the two actions.}
Take the first finite basis vector in degree zero and the second in
degree one. Then $\dd_F=\mu E_{21}$ has degree one,
$\dd_F^2=0$, $\dd_F^*=\mu E_{12}$ and
$Q_F=\mu(E_{21}-E_{12})$. For $b=\operatorname{diag}(b_1,b_2)$,
\[
 [D_F,b]=\mu\begin{pmatrix}0&b_2-b_1\\b_1-b_2&0\end{pmatrix},\qquad
 [Q_F,b]=\mu\begin{pmatrix}0&-(b_2-b_1)\\b_1-b_2&0\end{pmatrix}.
\]
Left multiplication by a diagonal $a$ preserves the rule that the
upper-right entry changes sign and the lower-left entry does not.
The rule holds for sums of represented one-forms. In particular, zero
ordinary matrix implies zero hatted matrix, proving descent to the
represented off-diagonal module.

Writing $\dd_{\rm tot}=\dd_\Lambda\otimes I+\gamma\otimes\dd_F$, one has
\[
 \dd_{\rm tot}^2=\dd_\Lambda^2\otimes I
 +(\dd_\Lambda\gamma+\gamma\dd_\Lambda)\otimes\dd_F
 +I\otimes\dd_F^2=0,
 \qquad \dd_{\rm tot}+\dd_{\rm tot}^*=D_\mu^{(2)}.
\]
This is a nilpotent Hilbert-complex differential. It does not by itself
establish module balancing for any residue functional.

For four constant finite one-forms, let their matrices be $C_j,H_j$; the total actions include a factor $\gamma$, which cancels in each commutator. The curvature residue is exactly
\begin{equation}\label{eq:finiteR}
 \RR_{D_\mu^{(2)}}(\omega_1,\omega_2,\omega_3,\omega_4)
 =-\frac{C_m}{2}\int_M\left(\frac{s}{12}+\mu^2\right)
 \tr_2([H_1,H_2][C_3,C_4])\dd V_g.
\end{equation}
To prove it, apply~\eqref{eq:static} to $L+\mu^2$. The mixed exterior potential has zero trace with a finite-only insertion, and the normalized exterior trace of $s/6-E_\Lambda$ is $-s/12$. For the explicit choice
$C_1=C_3=\mu E_{12}$, $C_2=C_4=\mu E_{21}$, one has
$\tr_2([H_1,H_2][C_3,C_4])=-2\mu^4$, and therefore
\begin{equation}
 \RR_{D_\mu^{(2)}}(\omega_1,\omega_2,\omega_1,\omega_2)
 =C_m\mu^4\int_M\left(\frac{s}{12}+\mu^2\right)\dd V_g.
\end{equation}
On a flat base this is $C_m\mu^6\vol(M)$. The value depends on the explicitly chosen finite differential and hatted action; it is not inferred from the original spinorial triple alone. Neither balancing nor all Riemann symmetries are claimed for the full finite-direction functional.

\paragraph{Complete evaluation for finite inputs.}
For finite forms, $[D,I\otimes b]=\gamma\otimes[D_F,b]$, so the total
actions are $\gamma\otimes C_j$ and $\gamma\otimes H_j$.
Each commutator contains two factors $\gamma$, leaving the insertion
$I_\Lambda\otimes K_F$ with $K_F=[H_1,H_2][C_3,C_4]$.
The mixed term of $E_\Lambda$ has trace zero because
\[
 R_{ijkl}\tau_\Lambda(c_ic_j\hc_k\hc_l)
 =-R_{ijkl}\delta_{ij}\delta_{kl}=0.
\]
Consequently
\[
 \tr_\Lambda(s/6-E_\Lambda-\mu^2)
 =2^n(s/6-s/4-\mu^2)=-2^n(s/12+\mu^2).
\]
Combining this identity with the factor
$\nu_{n-1}/(2(2\pi)^n)$ from the static residue proves~\eqref{eq:finiteR}
without a heat-kernel argument.

More explicitly write $C_j=\left(\begin{smallmatrix}0&p_j\\q_j&0\end{smallmatrix}\right)$
and $\Delta_{ab}=p_aq_b-p_bq_a$. Matrix multiplication gives
\[
 [C_a,C_b]=\Delta_{ab}\operatorname{diag}(1,-1),\qquad
 [H_a,H_b]=-\Delta_{ab}\operatorname{diag}(1,-1),
 \qquad \tr_2K_F=-2\Delta_{12}\Delta_{34}.
\]
For the displayed example, the forms $E_{11}\dd E_{22}$ and
$E_{22}\dd E_{11}$ have matrices $\mu E_{12}$ and $\mu E_{21}$.
Thus $[H_1,H_2]=-\mu^2\operatorname{diag}(1,-1)$ and
$[C_1,C_2]=\mu^2\operatorname{diag}(1,-1)$, giving $K_F=-\mu^4I_2$.
The factor $\mu^6$ on a flat base consists of $\mu^4$ from these four
represented inputs and $\mu^2$ from the potential. These inputs
themselves depend on $\mu$; the expression is not a statement about
the curvature of a fixed unit tangent plane of $M$.

\subsection{Conformal noncommutative tori: a vanishing but inadequate lift}
In the conformal torus example of \cite{DSZtorsion}, the positive smooth conformal factor belongs to the opposite algebra in the commutant, and $D_\kappa=\kappa D_0\kappa$. Use the exterior amplification of $D_0$ and make the explicit choice
\begin{equation}\label{eq:conformalactions}
 c_\kappa(e^j)=\kappa^2c_j,\qquad
 \hc_\kappa(e^j)=\kappa^2\hc_j.
\end{equation}
The first action is compatible with the commutators with the left algebra. The second action is an added choice. It is not the assertion that $\kappa D_0\kappa$ is the geometric Hodge operator for a conformally changed metric.

\paragraph{Left and right actions and ordered derivatives.}
To distinguish a coefficient from its representation, write
$\kappa=R_h$ with $h$ a positive invertible smooth element and
$L_a\xi=a\xi$, $R_h\xi=\xi h$. Associativity gives
$L_aR_h\xi=a(\xi h)=(a\xi)h=R_hL_a\xi$.
Thus $[\kappa,L_a]=0$ does not assert $ha=ah$ in the original algebra.
The right operators belong to the commutant of the left representation,
but that commutant need not be commutative. Indeed
\[
 R_hR_b=R_{bh},\qquad
 [\partial_j,R_h]=R_{\partial_jh},\qquad
 [R_h,R_{\partial_jh}]=R_{(\partial_jh)h-h(\partial_jh)}.
\]
Here $\partial_j\kappa$ means the coefficient derivative
$[\partial_j,R_h]$, not the full composition $\partial_j\circ R_h$.
For example, if $UV=e^{it}VU$ with $e^{it}\ne1$, take
$h=5+U+U^*+V+V^*\geq1$. This is smooth, positive and invertible.
For $\partial_1U=iU$, $\partial_1V=0$, the coefficient of $UV$ in
$[h,\partial_1h]=i[V+V^*,U-U^*]$ is $i(e^{-it}-1)\ne0$;
the other terms have different Fourier degrees. Thus a positive
right coefficient need not commute with its derivative. In particular
\[
 \partial_j(\kappa^2)=(\partial_j\kappa)\kappa
                         +\kappa(\partial_j\kappa),\qquad
 \partial_j(\kappa^{-1})=-\kappa^{-1}(\partial_j\kappa)\kappa^{-1}.
\]
The order in these identities must be retained.

Since $\kappa$ commutes with the left algebra,
$[D_\kappa,L_a]=\kappa[D_0,L_a]\kappa=\kappa^2c_jL_{\partial_ja}$.
This explains the first action in~\eqref{eq:conformalactions}, but does
not prescribe the second. Expanding on smooth vectors gives
\[
 D_\kappa=\kappa^2c_j\partial_j+\kappa c_j\partial_j\kappa,
 \qquad \sigma_1(D_\kappa)=i\kappa^2c(\xi),\quad
 \sigma_0(D_\kappa)=\kappa c_j\partial_j\kappa.
\]
Both symbols contain only unhatted generators.

\begin{proposition}[Blindness on the numerical frame]\label{prop:blind}
For the amplification~\eqref{eq:conformalactions}, every numerical-frame curvature component vanishes:
\begin{equation}
 \RR_{D_\kappa}(e^i,e^j,e^a,e^b)=0.
\end{equation}
The same holds if the first pair has commuting scalar coefficients and its hatted commutator is a pure hatted bivector.
\end{proposition}
\begin{proof}
Symbol composition differentiates coefficients and multiplies Clifford
words; neither operation introduces hatted generators. The inverse
recursion starts with the scalar principal inverse
$\kappa^{-4}|\xi|^{-2}$ of $D_\kappa^2$, so the same statement holds
at every order for $(D_\kappa^2)^{1-m}$.
For $i\ne j$, the inserted double commutator is
$2\kappa^8\hc_i\hc_j[c_a,c_b]$. Each product symbol therefore has
external matrix part $\hc_i\hc_jP(c)$, retaining the order of all
opposite-algebra coefficients. Its exterior trace is zero by the
nonempty-monomial rule of Section~\ref{sec:conventions}. Equivalently,
conjugation by $\hc_i$ changes the sign for even $P$, while odd $P$
gives an odd operator and hence zero ordinary trace. For $i=j$ the
first commutator itself is zero. In particular the degree $-n$ trace
vanishes before taking the algebra trace or integrating the covariable.
The same argument applies whenever the first hatted commutator has
no scalar component, as specified in the statement.
\end{proof}
The vanishing is not evidence that a conformal metric is flat. In the commutative specialization take $\kappa=e^f$, so the principal symbol corresponds to $g_\kappa=e^{-4f}g_0$. Let $f=\epsilon\cos x_1$, $\epsilon\ne0$. At $x_1=0$ its sectional curvature in the $(1,2)$ plane is
\begin{equation}
 K_{12}(g_\kappa)=-2\epsilon e^{4\epsilon}\ne0.
\end{equation}
This follows from $g_\kappa=e^{2\omega}g_0$, $\omega=-2f$, and
$K_{12}=-e^{-2\omega}(\omega_{11}+\omega_{22})$ at a critical point of $\omega$. Proposition~\ref{prop:blind} therefore proves a failure of classical conformal recovery for this naive lift. A genuine conformal Hodge construction has degree-dependent metric identifications and additional exterior-connection terms. The spinorial conformal rescaling cannot replace it.

\paragraph{Checking the classical curvature obstruction.}
In the commutative specialization $\kappa=e^f$, the squared principal
symbol is $e^{4f}|\xi|_{g_0}^2$, so the covariant metric is
$g_\kappa=e^{-4f}g_0$, not $e^{4f}g_0$. Put $\omega=-2f$.
The Christoffel symbols of $e^{2\omega}g_0$ are
\[
 \widetilde\Gamma^k_{ij}
 =\delta^k_j\omega_i+\delta^k_i\omega_j-\delta_{ij}\omega_k.
\]
At a critical point of $\omega$ they vanish, but their derivatives need
not vanish. In particular
$\widetilde\Gamma^1_{22}=-\omega_1$ and
$\widetilde\Gamma^1_{12}=\omega_2$, so
$(R(\partial_1,\partial_2)\partial_2)^1=-\omega_{11}-\omega_{22}$.
The sectional-curvature numerator has the extra metric factor
$e^{2\omega}$ and its denominator is $e^{4\omega}$, proving the formula
above. For $f=\epsilon\cos x_1$ at $x_1=0$,
$\omega=-2\epsilon$, $\omega_{11}=2\epsilon$, $\omega_{22}=0$,
which gives $-2\epsilon e^{4\epsilon}$.
Scalar weights supported near this point still give a pure hatted
bivector in the first pair and hence zero residue. Thus the failure
cannot be explained by cancellation in an unweighted integral.

The principal symbol determines a metric but not the exterior
connection and its lower-order curvature terms. For a genuine
conformal Hodge operator, the $L^2$ inner product on degree $r$ forms
changes by $e^{(n-2r)\omega}$, because the pointwise norm changes by
$e^{-2r\omega}$ and the volume by $e^{n\omega}$. This degree dependence
is absent from the uniform action choice~\eqref{eq:conformalactions}.

There is also a differential obstruction to simply taking $\dd_\kappa=\kappa\dd_0\kappa$:
\begin{equation}
 \dd_\kappa^2=\kappa[\dd_0,\kappa^2]\dd_0\kappa,
\end{equation}
which is generally nonzero. Thus~\eqref{eq:conformalactions} is being used as two-action data, not asserted to satisfy Definition~\ref{def:DR}.

Indeed
\[
 (\kappa\dd_0\kappa)^2
 =\kappa\dd_0\kappa^2\dd_0\kappa
 =\kappa[\dd_0,\kappa^2]\dd_0\kappa+\kappa^3\dd_0^2\kappa,
\]
and the last summand is zero. Already in the commutative case, on a
function $u$ this gives
\[
 (\kappa\dd\kappa)^2u
 =\kappa\dd(\kappa^2)\wedge\dd(\kappa u)
 =2\kappa^3\dd\kappa\wedge\dd u,
\]
which is nonzero whenever these two differentials are independent.
This disproves nilpotence for this candidate differential in general;
it does not rule out a different conformal Hodge construction.
Finally the numerical-frame conclusion must not be extended to all
noncommuting coefficients: for left coefficients $a,b$,
\[
 [a\kappa^2\hc_i,b\kappa^2\hc_j]
 =\kappa^4\left(\tfrac12(ab+ba)[\hc_i,\hc_j]
                         +\delta_{ij}[a,b]I\right).
\]
The second term has scalar hatted degree and survives the trace
selection rule. The next subsection retains precisely this possibility.

\subsection{Algebra-valued conformal inputs: an exact finite prescription}
The preceding zero statement must not be extended to arbitrary noncommuting coefficients. Write $u=u_ie^i$, and similarly for the other inputs, with coefficients in the left algebra. Taking the hatted scalar component gives
\begin{equation}\label{eq:conformalR}
 \RR_{D_\kappa}(u,v,w,z)=\frac1{n-2}
 \sum_{a,b}\tau\left(\Big(\sum_i[u_i,v_i]\Big)
                    (w_az_b-z_aw_b)\right)\mathsf H_{ab}(\kappa),
\end{equation}
under the left--right trace-factorization hypothesis used for strictly irrational tori in \cite{DSZtorsion}. Here $\mathsf H_{ab}$ is computed solely in the opposite coefficient algebra. The following recursion is an exact way of evaluating it without commuting $\kappa$ past its derivatives.

Let $q_0=|\xi|^2$ and put
\begin{align}
 d_1&=i\kappa^2c(\xi),& d_0&=\kappa c_j\partial_j\kappa,\\
 l_2&=\kappa^4q_0,&
 l_1&=d_1d_0+d_0d_1+\partial_{\xi_j}d_1D_{x_j}d_1,\\
 l_0&=d_0^2+\partial_{\xi_j}d_1D_{x_j}d_0.
\end{align}
All products are in the represented opposite algebra; this convention fixes the order of right multipliers. Compute $b_0,b_1,b_2$ by~\eqref{eq:inverse-rec}, with $b_j$ meaning the homogeneous term of degree $-2-j$. Set $S_{1,j}=b_j$. The following three recurrences suffice:
\begin{align}
 S_{a+1,0}&=S_{a,0}b_0,\\
 S_{a+1,1}&=S_{a,1}b_0+S_{a,0}b_1
                +\partial_{\xi_j}S_{a,0}D_{x_j}b_0,\\
 S_{a+1,2}&=S_{a,2}b_0+S_{a,1}b_1+S_{a,0}b_2\notag\\
 &\quad+\partial_{\xi_j}S_{a,1}D_{x_j}b_0
 +\partial_{\xi_j}S_{a,0}D_{x_j}b_1
 +\tfrac12\partial_{\xi_j}\partial_{\xi_l}S_{a,0}D_{x_j}D_{x_l}b_0.
\end{align}
For $a=m-1$, $S_{a,2}$ has degree $-n$. Set $q_0=1$ and apply the normalized $n$-dimensional moment rule~\eqref{eq:mom}, with $n-1$ replaced by $n$, to its remaining polynomial in $\xi$. Denoting that finite operation by $\mathcal M_n$, the answer is
\begin{equation}\label{eq:Hab}
 \mathsf H_{ab}(\kappa)=\nu_{n-1}\tau^o\left(
 \tr_\Lambda\big(c_ac_b\kappa^8\mathcal M_n(S_{m-1,2})\big)\right).
\end{equation}
Equations~\eqref{eq:conformalR}--\eqref{eq:Hab} are a finite ordered derivative formula for arbitrary $\kappa$ in every even dimension $n\geq4$. They are not claimed to be a minimal modular-functional-calculus expression. For constant $\kappa$, $d_0=l_1=l_0=0$ and all differentiated symbols vanish, so $S_{m-1,2}=0$ and the full functional vanishes, including algebra-valued inputs.

\subsection{The two new functionals on the other examples}
We use the same extra exterior data as in the preceding subsections. The two new orderings do not remove the need to specify the hatted action.

\paragraph{Gauge-coupled exterior data.}
For ordinary geometric one-forms tensored with the identity on the gauge fiber, the normal-coordinate calculation of Proposition~\ref{prop:universal} applies with the total connection. In both antisymmetrized insertions the scalar hatted contraction vanishes. The extra potential $c_ic_jF^A_{ij}/2$ and extra connection curvature contain no hatted generators, so their traces against $[\hc(X),\hc(Y)]$ vanish. The exterior-curvature part is the classical one, repeated over the $N^2$-dimensional gauge fiber. Consequently
\begin{equation}
 \RR_{D_A}^{(2)}(X,Y,Z,W)=\RR_{D_A}^{(3)}(X,Y,Z,W)
 =N^2C_m\int_M R(X,Y,Z,W)\dd V_g.
\end{equation}
This equality concerns the specified geometric submodule. Matrix-valued insertions must retain the ordered matrix products and the derivatives of the inserted endomorphisms; no extension of this displayed equality to that larger module is asserted.

\paragraph{Two sheets: geometric inputs.}
For the diagonal geometric submodule the extra mass has square $\mu^2\Id$ and anticommutes with the principal Clifford action. Its contribution to either first-pair antisymmetrization has zero hatted trace. Therefore
\begin{equation}
 \RR_{D_\mu^{(2)}}^{(2)}(X,Y,Z,W)
 =\RR_{D_\mu^{(2)}}^{(3)}(X,Y,Z,W)
 =2C_m\int_MR(X,Y,Z,W)\dd V_g.
\end{equation}

\paragraph{Two sheets: finite-direction inputs.}
Let $K_F=[H_1,H_2]$, $B_F=C_3C_4$, with $H_j,C_j$ as in~\eqref{eq:finiteR}. These are constant along $M$. The total one-form actions are $\gamma\otimes H_j$, $\gamma\otimes C_j$. Set $L=D_\Lambda^2$. The two relevant insertion products have second-order and zero-order parts
\begin{align}
 K_FD B_FD&=K_FB_FL+K_FD_FB_FD_F+\text{first-order terms},\\
 K_F\gamma C_3D\gamma C_4D
 &=-K_FB_FL+K_FC_3D_FC_4D_F+\text{first-order terms}.
\end{align}
The first-order terms have zero curvature-order residue. Since
\[
 L(L+\mu^2)^{-m}=(L+\mu^2)^{1-m}-\mu^2(L+\mu^2)^{-m},
\]
their traced residue coefficient is $-ks/12-m\mu^2$, rather than $-k(s/12+\mu^2)$. Keeping this extra mass term gives
\begin{align}\label{eq:finiteordered}
 \RR^{(2)}&=\frac{C_m}{k}\int_M\tr_2\!\left[
 K_F\left\{\left(-\frac{ks}{12}-m\mu^2\right)B_F+D_FB_FD_F\right\}\right]\dd V_g,\\
 \RR^{(3)}&=\frac{C_m}{k}\int_M\tr_2\!\left[
 K_F\left\{\left(-\frac{ks}{12}-m\mu^2\right)B_F-C_3D_FC_4D_F\right\}\right]\dd V_g.\notag
\end{align}
In the explicit test $C_1=C_3=\mu E_{12}$, $C_2=C_4=\mu E_{21}$, one has
\[
 \tr_2(K_FB_F)=-\mu^4,\quad
 \tr_2(K_FD_FB_FD_F)=\mu^6,\quad
 C_3D_FC_4D_F=0.
\]
Hence
\begin{align}
 \RR^{(2)}&=C_m\mu^4\int_M\left(\frac{s}{12}+\frac{m+1}{m-1}\mu^2\right)\dd V_g,\\
 \RR^{(3)}&=C_m\mu^4\int_M\left(\frac{s}{12}+\frac{m}{m-1}\mu^2\right)\dd V_g.
\end{align}
On a flat base these become $(m+1)C_m\mu^6\vol(M)/(m-1)$ and $mC_m\mu^6\vol(M)/(m-1)$, respectively. The original functional gives $C_m\mu^6\vol(M)$. Thus the finite-direction calculation explicitly distinguishes all three orderings.

\paragraph{Conformal torus: numerical-frame inputs.}
All symbols of $D_\kappa$ and its powers contain only unhatted Clifford generators. For distinct first-pair frame indices, the commutator is a pure hatted bivector, so its trace against either ordered insertion vanishes; for equal indices the commutator itself vanishes. Thus
\begin{equation}
 \RR^{(2)}_{D_\kappa}(e^i,e^j,e^a,e^b)
 =\RR^{(3)}_{D_\kappa}(e^i,e^j,e^a,e^b)=0.
\end{equation}
The same conformal classical-recovery obstruction described in Proposition~\ref{prop:blind} therefore applies to both new functionals.

\paragraph{Conformal torus: arbitrary coefficients.}
For completeness the following finite prescription includes derivatives of the inserted one-forms and keeps the order of $\kappa$. Use $d_1,d_0$ and $S_{m,j}$ from the preceding subsection, now with exponent $m$, not $m-1$. For arbitrary $A,B$ form
\begin{align}
 p_2&=Ad_1Bd_1,\\
 p_1&=A\{d_1Bd_0+d_0Bd_1+\partial_{\xi_j}d_1D_{x_j}(Bd_1)\},\\
 p_0&=A\{d_0Bd_0+\partial_{\xi_j}d_1D_{x_j}(Bd_0)\}.
\end{align}
Define the explicit homogeneous polynomial on the unit sphere
\begin{align}
 \mathcal Z(A,B)={}&p_2S_{m,2}+p_1S_{m,1}+p_0S_{m,0}\notag\\
 &+\partial_{\xi_j}p_2D_{x_j}S_{m,1}
 +\partial_{\xi_j}p_1D_{x_j}S_{m,0}
 +\tfrac12\partial_{\xi_j}\partial_{\xi_l}p_2D_{x_j}D_{x_l}S_{m,0}.
\end{align}
The two answers are
\begin{align}
 \RR^{(2)}&=\frac{\nu_{n-1}}k(\tau\otimes\tau^o)
 \tr_\Lambda\mathcal M_n\mathcal Z(H,c_\kappa(w)c_\kappa(z)),\\
 \RR^{(3)}&=-\frac{\nu_{n-1}}k(\tau\otimes\tau^o)
 \tr_\Lambda\mathcal M_n\mathcal Z(Hc_\kappa(w),c_\kappa(z)).
\end{align}
Here $H=[\hc_\kappa(u),\hc_\kappa(v)]$; the left--right trace hypothesis is the same as for~\eqref{eq:conformalR}. Every product and derivative is specified, and $\mathcal M_n$ is the finite spherical moment rule. For general noncentral $\kappa$ this is an ordered finite coefficient answer, not a claimed minimal modular formula. For constant $\kappa$, all coefficients of $L$ below its principal symbol vanish, so both functionals vanish even for arbitrary algebra-valued inputs.

\section{Conclusions and boundaries of the results}
The four classical mixed residues recover scalar, Ricci, Einstein and Riemann information with explicit coefficients. The mixed $c\,c\,\hc\,\hc$ part of the exterior Weitzenbock potential is essential. For warped collars the specified boundary extension has explicit geometric densities in every even dimension at least four. The five separate cases can be nonzero even when their sum vanishes. The first-pair antisymmetrization of the second or third insertion cancels the full boundary density.

The noncommutative curvature residue has a well-defined scope. It agrees with the weighted Riemann tensor on classical Hodge data, vanishes on the flat noncommutative torus, and has the explicit perturbation formula~\eqref{eq:torusR}. Nevertheless it can fail module balancing, pair interchange, and Bianchi. The conformal amplification illustrates a different limitation: vanishing frame components can reflect an inadequate choice of hatted action rather than flat geometry. The finite two-sheet computation demonstrates how an additional differential choice changes which questions can be answered.

A fully tensorial noncommutative Riemann theory would require an appropriate bimodule or braided symmetry framework and a comparison with compatible connections. The present definitions and counterexamples do not solve that general reconstruction problem. Likewise the boundary results do not cover arbitrary non-umbilic collars or all alternative truncation factorizations. These limitations specify the proved statements rather than leave implicit assumptions in the formulas.

\section*{Acknowledgements}
\phantomsection
\addcontentsline{toc}{section}{Acknowledgements}
The author was supported by Science and Technology Development Plan Project of Jilin Province, China: No.\ 20260102245JC and NSFC No.\ 11771070.

\section*{Data Availability}
\phantomsection
\addcontentsline{toc}{section}{Data Availability}
Data sharing is not applicable to this article as no new data were created or analyzed
in this study.

\section*{Declarations}
\phantomsection
\addcontentsline{toc}{section}{Declarations}
\noindent\textbf{Conflict of interest}
The authors have no conflict of interest to declare that are relevant
to the content of this article.

\end{document}